\documentclass[leqno, myheadings, twoside]{amsart}

\usepackage{amsmath, amsthm, amssymb, amscd, amsxtra,graphicx}
\usepackage{latexsym, amsfonts}
\usepackage{easybmat}
\usepackage{etex}
\usepackage{url}
\usepackage{texdraw}
\usepackage{epsfig}
\usepackage{tikz}
\allowdisplaybreaks[4]

\usepackage{geometry}
\makeatletter \@addtoreset{equation}{section}

\makeatletter \renewcommand{\@biblabel}[1]{#1.}
\theoremstyle{remark}

\begin{document}

\title [Determination of isometric real-analytic metric and spectral invariants] {Determination of isometric real-analytic metric and spectral invariants for elastic Dirichlet-to-Neumann map on Riemannian manifolds}
\author{Genqian Liu}


\maketitle Department of Mathematics, Beijing Institute of
Technology,  Beijing 100081,  People's Republic of China.
 \ \    E-mail address:  liugqz@bit.edu.cn

\vskip 0.46 true cm

\vskip 0.45 true cm

\vskip 15 true cm

\begin{abstract}   \ \ In this paper, the elastic Dirichlet-to-Neumann map $\Xi_g$ is studied for the stationary elasticity system in a compact Riemannian manifold $(\Omega,g)$ with smooth boundary $\partial \Omega$. By overcoming methodological difficulties, we explicitly get matrix-valued full symbol for the elastic Dirichlet-to-Neumann map $\Xi_g$.  We prove that for a strong convex or extendable real-analytic manifold with boundary, the elastic Dirichlet-to-Neumann map $\Xi_g$ uniquely determines the metric $g$ of $\Omega$ in the sense of isometry, thereby solving an open problem for the uniqueness of the metric under real-analytic setting. Furthermore, by calculating the symbol representation of the resolvent operator $(\Xi-\tau I)^{-1}$ we can explicitly obtain all coefficients $a_0, a_1 \cdots, a_{n-1}$ of the asymptotic expansion $\sum_{k=1}^\infty e^{-t \tau_k}\sim \sum_{m=0}^{n-1} a_m t^{m+1-n} +o(1)$ as $t\to 0^+$, where $\tau_k$ is the $k$-th eigenvalue of the elastic Dirichlet-to-Neumann map $\Xi_g$ (i.e., $k$-th elastic Steklov eigenvalue). These coefficients (spectral invariants) provide important geometric information for the manifold, which give an answer to another open problem for the elastic Steklov spectral asymptotics.  \end{abstract}

\vskip 1.49 true cm

\section{Introduction}

\vskip 0.45 true cm

Let $(\Omega,g)$ be a smooth compact Riemannian manifold of
dimension $n$ with smooth boundary $\partial\Omega$ of dimension $n-1$. The manifold $\Omega$ is considered here as an elastic, isotropic, homogenous medium with Lam\'{e} constants $\mu$ and $\lambda$. Assume that Lam\'{e} constants $\mu$ and $\lambda$ satisfy $\mu>0$ and $\lambda+\mu\ge 0$.
Under the assumption of no body forces acting on $\Omega$, the boundary value
problem for the displacement vector $u(x) = (u^1(x), \cdots, u^n(x))$ in $\Omega$ produced by a displacement $f$ on $\partial \Omega$ is given by
the elastic Lam\'{e} system:
\begin{eqnarray} \label{18/12/22-1} \;\;\;\;\;\; \left\{ \!\!\! \begin{array} {ll} {\mathcal{L}}_g u:=\!-\!\mu \nabla^*\nabla u \!+\! (\mu\!+\!\lambda)\mbox{grad} \, \mbox{div}\,u \!+\!\mu \,\mbox{Ric}(u)\!=\!0 \,\; &\mbox{in}\,\, \Omega,\\
u=f \;\; &\mbox{on}\;\, \partial \Omega, \end{array} \right. \end{eqnarray}
 where $-\nabla^*\nabla$ is the Bochner Laplacian (see (\ref{19.9.30-1}) in Section 2),  $\mbox{div}$ and $\mbox{grad}$ are  the usual divergence and gradient operators, and  \begin{eqnarray} \label{18/12/22} \mbox{Ric} (u)= \big(\sum\limits_{l=1}^n R^{\;1}_{\,l} u^l,  \sum\limits_{l=1}^n R^{\;2}_{\,l} u^l, \cdots, \sum\limits_{l=1}^n R^{\;n}_{\,l} u^l\big)\end{eqnarray} denotes the action of Ricci tensor $\mbox{R}_l^{\;j}:=\sum_{k=1}^n R^{k\,\,j}_{\,lk}$ on $u$.
For  $u\in H^2(\Omega)$,  $$\mbox{traction} \;u:=2\mu \,(\mbox{Def}\;{u})^{\#} {\nu} +\lambda (\mbox{div}\; {u}){\nu} \;\; \mbox{on}\;\; \partial \Omega$$ is said to be the elastic Neumann boundary condition (which just is the normal component of the stress at the boundary), where  $(\mbox{Def}\, u)_{jk}= \frac{1}{2} \big(u_{j;k} +u_{k;j}\big)$ is the strain tensor, $\#$ is the sharp operator (for a tensor) by raising index and $\nu=(\nu_1, \cdots, \nu_n)$ is the outward unit normal vector to $\partial \Omega$.

We define the elastic Dirichlet-to-Neumann map $\Xi_g$ as follows:
\begin{eqnarray} \label{2022.5.25-3} \Xi_g (f)= \big(2\mu \,(\mbox{Def}\;{u})^{\#} {\nu} +\lambda (\mbox{div}\; {u}){\nu}\big)\big|_{\partial \Omega},\;\, \;\forall f\in H^{3/2}(\partial \Omega),\end{eqnarray}
where $u$ satisfies  (\ref{18/12/22-1}) (i.e., ${\mathcal{L}}_g u=0$ in $\Omega$ and $u|_{\partial \Omega}=f$).
Clearly, the elstic Dirichlet-to-Neumann map is well defined for each $f \in H^{3/2}(\partial \Omega)$ with value in $H^{1/2}(\partial\Omega)$.
 $\Xi_g$ is also said to be the {\it displacement-to-traction map} on the elastic boundary.

Since the elastic Dirichlet-to-Neumann map $\Xi_g$ is an elliptic, self-adjoint pseudodifferential operator of order one (see section 3), there exists a sequence of eigenvalues
\begin{eqnarray} 0\le\tau_1\le \tau_2 \le \cdots \le \tau_k \le \cdots \nearrow \infty\end{eqnarray} such that $\Xi_g v_k= \tau_k v_k$, where $v_1, v_2, \cdots, v_k,\cdots$ is an orthogonal basis of eigenvectors in $(L^2(\partial \Omega))^n$=$L^2(\partial \Omega)\times \cdots \times L^2(\partial \Omega)$ corresponding to the eigenvalues $\tau_1, \tau_2, \cdots, \tau_k, \cdots$.
Clearly, all the eigenvalues $\{\tau_k\}$ of the elastic Dirichlet-to-Neumann map are just all the {\it elastic Steklov eigenvalues} $\{\tau_k\}$:
 \begin{eqnarray} \label{18/12/22-4} \;\;\;\;\, \left\{\!\! \begin{array} {ll} -\mu \,\nabla^*\nabla v_k\!+\! (\mu\!+\!\lambda)\,\mbox{grad} \, \mbox{div}\,v_k \!+\!\mu \,\mbox{Ric}(v_k)=0 \;\; &\mbox{in}\;\, \Omega,\\
\Xi_g\,  v_k=\tau_k v_k  \;\; &\mbox{on}\;\, \partial \Omega. \end{array} \right. \end{eqnarray}

Let us point out that if $\Omega$ is a bounded domain in ${\Bbb R}^n$ with the standard Euclidean metric, then the Lam\'{e} operator ${\mathcal{L}}_g$ reduces to the classical elastic operator $\mathcal{L}u= -\mu \Delta u +(\mu +\lambda) \, \mbox{grad}\, \mbox{div}\, u$ for $u\in (H^n(\Omega))^n$; and further if $\lambda+\mu=0$ and the traction is replaced by the normal derivative of $u$ in the definition of the elastic Dirichlet-to-Neumann map $\Xi_g$, then $\Xi_g$ becomes the classical Dirichlet-to-Neumann map $\Lambda_g$ associated with the Laplace operator.

  The elastic Lam\'{e} operator is a complicated ``exotic'' differential operator because its second-order terms are not merely the Laplacian or even the Laplace-Beltrami operator of a nontrivial metric (see  \cite{Avr7}, \cite{Avr10}, \cite{Full},  \cite{KLV}, \cite{Kupr}, \cite{KGB}, \cite{LaLi},  \cite{SV}  and \cite{So}). The study of the elastic Dirichlet-to-Neumann map is important in areas such as geophysical exploration (for example, one may explore the interior structure of the earth by measuring the displacements and tractions of the earth's surface by earthquakes (or artificial earthquakes)), materials characterization and acoustic
emission of many important materials and nondestructive testing (see \cite{NU1} and \cite{NU2}). However, the
commutativity properties (which are usually applied for scalar differential operators) are not valid for the elastic Dirichlet-to-Neumann map, so there are additional difficulties in the elastic case (see also, p.$\,$166 of \cite{Isak}).

\vskip 0.16 true cm

The following two well-known problems have been open for a long time.

\vskip 0.16 true cm

 Problem A:  \  Given a compact $n$-dimensional Riemannian manifold with boundary, $n\ge 2$, if its boundary metric is fixed, can one determine the metric of the whole manifold (in the sense of isometry) by the elastic Dirichlet-to-Neumann map?

\vskip 0.14 true cm

  Problem B:  \ \  What geometric information of the spectral asymptotics can be explicitly obtained for the Riemannian manifold by knowing all the elastic Steklov eigenvalues?

\vskip 0.2 true cm

Problem A is the Lee-Uhlmann type conjecture for the elastic Dirichlet-to-Neumann map. In \cite{LU}, J. Lee and G. Uhlmann raised and solved several important conjectures for the Dirichlet-to-Neumann map $\Lambda_g$ associated with Laplace operator. They proved that for a strongly convex real-analytic $n$-dimensional Riemannian manifold with connected real-analytic boundary or extendable real-analytic Riemannian manifold, if $\pi_1(\Omega, \partial \Omega)=0$ (this topological assumption means that every closed path in $\bar \Omega$ with base point in $\partial \Omega$ is homotopic to some path that lies entirely in $\partial \Omega$, see \cite{LU} and p.$\,$162 of \cite{Whi}), then the Dirichlet-to-Neumann map $\Lambda_g$ determines the metric $g$ in the sense of isometry for $n\ge 3$, and determines $g$ in the sense of conformal isometry for $n=2$, where $\Lambda_g: H^{1/2} (\partial \Omega)\to H^{-1/2} (\partial \Omega)$ is defined by $\Lambda_g f=\frac{\partial u}{\partial \nu}\big|_{\partial \Omega}$, and $u$ satisfies \begin{eqnarray} \label{18/12/22-5} \left\{ \begin{array}{ll} \Delta_g u=0 \quad &\mbox{in}\;\; \Omega,\\
u=f\quad &\mbox{on}\;\; \partial \Omega.\end{array}\right.\end{eqnarray}  For the case of electromagnetic field, similar problem and some related results can be seen in \cite{KeSU}, \cite{Liu1}, \cite{KV1}, \cite{KV2}, \cite{McD}. Problem A actually is an general version of the Calder\'{o}n problem (see \cite{Cald}) in the case of elasticity. It has attracted much attention during the past decades, and some progresses have been made (see, for example, \cite{EINT}, \cite{CHKU}, \cite{Rac}, \cite{Sto}  and \cite{Bel}). Note that when $\mu+\lambda=0$, Problem A is false in two dimensions  because of the conformal invariance
of the Dirichlet problem (cf. \cite{LU}).
\vskip 0.12 true cm

Problem B is similar to the famous Kac's problem in \cite{Kac} (see also \cite{Lo} or \cite{We1}) and Protter's propose (see p.$\,$ 187 of \cite{Pro}) for the classical eigenvalue problem (i.e., can one hear the shape of a medium by hearing all tones of vibration of this  medium? see \cite{We5}, \cite{We6}, \cite{Liu4}). For the eigenvalues $\{\varsigma_k\}$ of the Dirichlet-to-Neumann map associated with Laplace operator, the author of the present paper \cite{Liu2} calculated the first fourth coefficients of the heat trace asymptotic expansion (Polterovich and D. A. Sher calculated the first three coefficients in \cite{PS} by a different way). A generalization of the above result for the perturbed polyharmonic Steklov problem has been given by the author of the present paper in \cite{Liu3}. The asymptotic formula with a sharp remainder estimate for the counting function of the Steklov eigenvalues $\{\varsigma_k\}$ can be seen in \cite{Liu3} (see also,  \cite{Liu5}, \cite{Liu6}). Therefore, it is quite important to explicitly give the asymptotic expansion of the trace for the elastic Dirichlet-to-Neumann map.

In this paper, by overcoming methodological difficulties, we first solve open Problem A when $n\ge 2$ under the same assumptions as given by Lee and Uhlmann in \cite{LU}; furthermore we get the fundamental geometric quantities of asymptotic expansion of the trace of the elastic Dirichlet-to-Neumann map for  Problem B (In other words, in the sense of spectral invariants, we solve the open Problem B).

Our main results are the following:

\vskip 0.2 true cm

\noindent{\bf Theorem 1.1.} \ {\it Let $\bar \Omega$ be a compact, connected, real-analytic $n$-manifold with real-analytic boundary, $n\ge 2$, and assume that $\pi_1(\bar \Omega, \partial \Omega)=0$. Suppose that the Lam\'{e} constants $\mu>0$ and $\lambda$ satisfy $\lambda+\mu> 0$ and $(n-1) \lambda^3+(4n-2) \lambda^2\mu +(n+5)\lambda \mu^2 +(14-8n)\mu^3 \neq 0$. Let $g$ and $\tilde{g}$ be real analytic metric on $\bar \Omega$ such that $\Xi_g=\Xi_{\tilde g}$, and assume that one of the following conditions holds:

 (a) \ \ $\bar \Omega$ is strongly convex with respect to both $g$ and $\tilde {g}$;

 (b) \ \  either $g$ or $\tilde{g}$ is extendable (i.e., it can be extend to a complete real-analytic metric on a non-compact real-analytic manifold $\tilde{\Omega}$ (without boundary) containing $\bar \Omega$).

 Then there exists a real-analytic diffeomorphism $\varrho: \bar \Omega\to \bar \Omega$ with $\varrho\big|_{\partial \Omega}=\mbox{identity}$, such that
 $g=\varrho^* \tilde {g}$.}

\vskip 0.25 true cm

  Since the above condition $\pi_1(\bar \Omega, \partial \Omega)=0$ and extendable assumption (b) hold when $\Omega$ is a bounded simply-connected real-analytic domain in ${\Bbb R}^n$ and $g$ is the Euclidean metric, we immediately have the following corollary:

\vskip 0.2 true cm

 \noindent{\bf Corollary 1.2.} \ {\it Let $\Omega \subset {\Bbb R}^n$, $n\ge 2$, be a simply-connected bounded open set
with connected real-analytic boundary, and suppose $\tilde{g}$ is a real-analytic metric on $\bar \Omega$
such that $\Xi_{\tilde{g}} = \Xi_g$, where $g$ is the Euclidean metric. If the Lam\'{e} constants $\mu>0$ and $\lambda$ satisfy $\lambda+\mu> 0$ and  $(n-1) \lambda^3+(4n-2) \lambda^2\mu +(n+5)\lambda \mu^2 +(14-8n)\mu^3 \neq 0$, then there exists a real-analytic diffeomorphism $\varrho:\bar \Omega\to \bar \Omega$ with $\varrho\big|_{\partial \Omega}=\mbox{identity}$, such that $g=\varrho^*{\tilde{g}}$.}

\vskip 0.22 true cm

Clearly, the knowledge provided by the elastic Dirichlet-to-Neumann map (i.e., the set $\{(u\big|_{\partial \Omega}$, $(\mbox{traction}\, u)\big|_{\partial \Omega})\}$ of all the displacement and traction data) is equivalent to the information given by all the elastic Steklov eigenvalues and the corresponding eigenvectors. However, if we only know the information of all elastic Steklov eigenvalues, then we have the following:

\vskip 0.20 true cm

\noindent{\bf Theorem 1.3.} \ {\it  Suppose that $(\Omega,g)$ is an $n$-dimensional, smooth Riemannian
manifold with smooth boundary $\partial \Omega$.  Assume that $\mu>0$ and $\lambda+\mu> 0$. Let $\{\tau_k\}$ be the all eigenvalues of the elastic Dirichlet-to-Neumann map $\Xi_g$.
Then} \begin{eqnarray} \label{19.6.15-1} \;\;\;\;\;\sum_{k=0}^\infty e^{-t\tau_k} = \mbox{Tr}\; e^{-t \Xi_g} \sim  \sum_{m=0}^{n-1} a_m t^{-n+m+1} + o(1)\, \quad \mbox{as}\;\; t\to 0^+,\end{eqnarray}
{\it where $a_m$ are constants, which can be explicitly calculated by the procedure given in section 5 for $m<n$.
In particular, if $n\ge 2$, then
   \begin{eqnarray} \label{6.0.1}    && \sum_{k=1}^\infty e^{-t\tau_k} =t^{1-n}
\int_{\partial \Omega} a_0(x') \,dS(x') + t^{2-n}  \int_{\partial \Omega} a_1(x')\, dS(x') \\  &&\quad \quad \quad\quad  \quad + \left\{\begin{array}{ll} O(t^{3-n}) \quad \;\, \mbox{when}\;\; n>2,\\ O(t\log t) \quad \, \mbox{when} \;\; n=2,\end{array}\right. \quad \;\;\mbox{as}\;\; t\to 0^+.\nonumber\end{eqnarray}
Here}
   \begin{align}\label{109.6.14-1} &a_0(x')=\frac{( n-2 )!\,\mbox{vol}( \mathbb{S}^{n-2}) }{(2\pi)^{n-1}\mu ^{n-1}}\bigg\{ \frac{1+\big( \frac{\lambda +{3}\mu}{\lambda +\mu} \big) ^{n-1}}{2^{n-1}}+\left( n-2 \right) \bigg\},\\
\label{06-060} &a_1(x')= \frac{(n\!-\!2)!\,\mbox{vol}( \mathbb{S}^{n\!-\!2} )\big(\sum_\alpha\!{\kappa _{\alpha}}\big) }{(2\pi)^{n-1}}\bigg\{
\!\frac{n(2-n)}{2(n\!-\!1)\mu^{n-2}}\\
 &\qquad \qquad  \!-\! \bigg(\! 1\! +\! \frac{(\lambda^3 \!+\!23 \lambda^2 \mu \!+\!87 \lambda \mu^2 \!+\!97 \mu^3)}{ 4(n-1)(\lambda+3\mu)^3} \bigg) \frac{\mu}{(2\mu)^{n-1}} \nonumber \\
&\qquad \qquad  + \bigg( \frac{2\lambda^3 +9 \lambda^2\mu+ 10\lambda\mu^2 -\mu^3}{(\lambda+3\mu)^2} \nonumber \\
& \qquad \qquad + \frac{\mu (15\lambda^3 +105\lambda^2 \mu +233\lambda\mu^2 +175 \mu^3)}{4(n-1)(\lambda+3\mu)^3} \bigg) \bigg(\frac{\lambda+3\mu}{2\mu(\lambda+\mu)}\bigg)^{n-1}\bigg\}, \nonumber   \end{align}
 {\it where $\kappa_1(x'),\cdots, \kappa_{n-1}(x')$ are the principal curvatures of the boundary $\partial \Omega$ at $x'\in \partial \Omega$}, $\text{vol}\left( \mathbb{S}^{n-2} \right) = 2\pi^{(n-1)/2}/\Gamma(\frac{n-1}{2}) $ {\it is the volume of $(n-2)$-dimensional unit sphere ${\Bbb S}^{n-2}$ in ${\Bbb R}^{n-1}$, and $dS$  denotes the volume element on $\partial \Omega$.}

\vskip 0.25 true cm

The heat invariants $a_m$ are spectral invariants of the operator $\Xi_g$ that encode the information about the asymptotic properties of the spectrum. They are of great importance in spectral geometry and have extensive application in physics because they describe real physical phenomena.
Note that the eigenvalues of the elastic Dirichlet-to-Neumann map can be
measured experimentally. Therefore, Theorem 1.3 shows that the boundary volume ${\mbox{vol}}(\partial \Omega)$, the total mean curvature of the boundary surface $\partial \Omega$  are all spectral invariants and can also be obtained by elastic Steklov eigenvalues (i.e., physical quantities).
\vskip 0.25 true cm

The main ideas of this paper, which contain a lot of novel techniques and methods, are as follows. For a Riemannian manifold $(\Omega, g)$, we first derive a precise form of the elastic Lam\'{e} operator:
 ${\mathcal{L}}_g u= -\mu \, \nabla^* \nabla u+(\lambda+\mu) \,\mbox{grad}\, \mbox{div}\, u+\mu\,\mbox{Ric}(u)$. From this, we see the elastic Lam\'{e} operator ${\mathcal{L}}_g$  possesses an invariant form under a local coordinate change in a general $n$-dimensional Riemannian manifold. Secondly, in boundary normal coordinates, we obtain a local representation for the elastic Dirichlet-to-Neumann map:
 \begin{eqnarray}\!\!\!&&\label{18/12/23-6} \Xi_g (u\big|_{\partial \Omega})=\left\{ \left[
 \begin{BMAT}(@, 6pt, 6pt){c.c}{c.c}\big[\mu \delta_{jk} \big]_{(n-1)\times (n-1)}&\big[0\big]_{(n-1)\times 1}\\
\big[0\big]_{1\times (n-1)} & \lambda+2\mu\end{BMAT}\right] \bigg[-\frac{\partial }{\partial x_n}I_n\bigg]\right.\quad \\
\!\!\!\!&& \qquad \quad \quad \left. - \left[ \begin{BMAT}(@, 6pt, 6pt){c.c}{c.c} [0]_{(n-1)\times (n-1)}& \bigg[\mu\frac{\partial}{\partial x_j}\bigg]_{(n-1)\times 1}\\
\bigg[\lambda\big(\frac{\partial} {\partial x_k} +\sum\limits_{\alpha} \Gamma_{k\alpha}^\alpha\big)\bigg]_{1\times (n-1)} & \lambda\sum\limits_{\alpha}\Gamma_{n\alpha}^\alpha\end{BMAT}\right]\right\}\begin{bmatrix}u^1\\
  \vdots\\
 u^{n-1}\\
 u^n\end{bmatrix},   \nonumber \end{eqnarray}
 where $I_n$ denotes the $n\times n$ identity matrix. (\ref{18/12/23-6}) is a basic formula, and we can call it as a transform relation about the operator $-\frac{\partial }{\partial x_n} I_n$ and $\Xi_g$. Therefore, it suffices to give a local representation for the $-\frac{\partial }{\partial x_n}I_n$ on $\partial \Omega$ from the elastic Lam\'{e} equation ${\mathcal{L}}_g u=0$.
Unlike acoustic or electromagnetic waves, in elastic medium (see \cite{Gur}, \cite{LaLi} or \cite{Sl}), elastic waves have both longitudinal and shear waves (Note that the acoustic waves are a type of longitudinal waves that propagate by means of adiabatic compression and decompression;
the electric and magnetic fields of an electromagnetic wave are perpendicular to each other
and to the direction of the wave, see \cite{CH}, \cite{Coo} or \cite{Stra}).
 Unfortunately, many standard methods for the classic Dirichlet-to-Neumann map are no longer valid for the elastic Dirichlet-to-Neumann map.
Now, we should consider an equivalent elastic equation
    \begin{eqnarray*}  0=\begin{bmatrix} \frac{1}{\mu}I_{n-1}  &0 \\
   0& \frac{1}{\lambda +2\mu} \end{bmatrix} {\mathcal{L}}_g u=  \left\{ \bigg[\frac{\partial^2}{\partial x_n^2}I_n\bigg]+B\bigg[\frac{\partial}{\partial x_n}I_n\bigg] +{C} \right\}\begin{bmatrix}u^1\\
  \vdots\\
 u^{n-1}\\
 u^n\end{bmatrix},\end{eqnarray*}
     and look for the factorization $$\frac{\partial^2}{\partial x_n^2}I_n+B\left[\frac{\partial}{\partial x_n}I_n\right] +{C}= \left(\left[\frac{\partial }{\partial x_n}I_n\right] +B+Q\right)\left(\left[\frac{\partial }{\partial x_n}I_n\right] +Q\right),$$ where operator $Q$ will be determined late. Because the matrix $B$ is a differential operator of order one, we will encounter three huge difficulties:

      i) \ \  How to solve the unknown $q_1$ from the following matrix equation?
       \begin{eqnarray} \label{19.4.30-1} q_1^2+b_1 q_1 -c_2=0,\end{eqnarray}
       where $q_1$, $b_1$ and $c_2$ are the principal symbols of the differential (or pseudodifferential) operators $Q$, $B$ and $C$, respectively.

        Generally, the quadratic matrix equation of the form (\ref{19.4.30-1}) can not be exactly solved (in other words, there is not a formula of the solution represented by the coefficients of matrix equation (\ref{19.4.30-1})). However, in our setting we can get the exact solution by applying a method of algebra ring theory in this paper. More precisely, by observing the coefficients of matrix equation (\ref{19.4.30-1}) we define an invariant sub-algebra ring $\mathfrak{F}$ which is generated by coefficients matrices of equation (\ref{19.4.30-1}), and $\mathfrak{F}$ is acted (left-multiply) by another sub-algebra ring (so-called ``coefficient ring'') of special diagonal matrices. This implies that the $q_1$ has a special form (see section 3), and hence by solving a higher-order polynomial equations (a system of resultant equations) for the unknown constants in $q_1$, we obtain an exact solution $q_1$ (it is a surprising method and result$\,$!). This method is inspired by Galois group theory to solve the polynomial equation (see, for example, \cite{Art} or \cite{HME}).

       ii) \ \   How to solve Sylvester's equation $(q_1-b_1) q_{j-1} +q_{j-1} q_1=E_j\,$? where $q_{j-1}$ ($j\le 1$) are the remain symbols of $q$ (here $q\sim \sum_{j\le 1} q_j$),  and $E_j$ can be seen in section 3.

In mathematics (more precisely, in the field of control theory), a Sylvester equation is a matrix equation of the form (see \cite{Syl} and \cite{BaS}):
\begin{eqnarray}\label{19.6.22-1}  LX+XM=V. \end{eqnarray}
Then given matrices $L$, $M$, and $V$, the problem is to find the possible matrices $X$ that obey this equation. A celebrated result (see \cite{BaS} or \cite{BhR}) states that Sylvester's equation (\ref{19.6.22-1}) has a unique solution $X$ for all $V$ if and only if $L$ and $-M$ have no common eigenvalues. Generally, it is impossible to obtain an explicit solution of the Sylvester's equation. However, by putting Sylvester's equation into an equivalent $n^2\times n^2$-matrix equation and by applying our invariants sub-ring method mentioned above, we can get the inverse $U^{-1}$ of
$n^2\times n^2$-matrix $U:= (I_n\otimes L) + (M^t \otimes I_n)$ and further obtain the exact solution of $q_{j-1}$, $(j\le 1)$ (another surprise result!)
     Therefore, the elastic Dirichlet-to-Neumann map $\Xi_g$ on $\partial \Omega$ is obtained from transform relation (\ref{18/12/23-6}) because we have got the pseudodifferental operator $Q$ and $-\frac{\partial}{\partial x_n}  I=Q$ modulo a smoothing operator on $\partial \Omega$.
 By analysing the full symbol of $\Xi_g$, we find that the $\Xi_g$ determines the metric $g$ and all its tangent and normal derivatives of any order at every point $x_0\in \partial \Omega$. Combining the expansion of Taylor's series of $g_{jk}$ at every $x_0\in \partial \Omega$, the analytic of the Riemannian manifold $\Omega$ and some given topology assumptions, we can complete the proof of Theorem 1.1.

iii)  \ \  How to obtain the precise coefficients in the asymptotic expansion of the integral of the parabolic trace for the elastic Dirichlet-to-Neumann map?

The heat invariants have been studied for the Laplacian with manifolds without boundaries (see \cite{Gil} or \cite{Gil2}) as well as for the Dirichlet-to-Neumann map associated with the
 Laplacian (see \cite{Liu2} or \cite{PS}) by considering heat
trace asymptotics and applying symbol calculus (see \cite{See2} and \cite{Gr}). However, there is no any result in the context of the elastic Steklov problem because of the technical complexity and, most importantly, lack of symbol representation for the ``exotic'' elastic Dirichlet-to-Neumann map $\Xi_g$ (In fact, ``most of the studies in spectral geometry and spectral asymptotics are restricted to so-called Laplace type operators'', see p.$\,$120 of \cite{Avr10}).
In fact, there are essential difficulties for studying the  asymptotic expansion of the trace of elastic Dirichlet-to-Neumann map.  Firstly, we can show (see section 5) that $e^{-t \Xi_g}$ admits an asymptotic expansion
\begin{align*} \label{18/12/23-9}  & \sum\limits_{k=0}^\infty\! e^{-t\tau_k}\! = \!\int_{\partial \Omega} \!\mbox{Tr}\,(K(t, x',x'))dS(x')\qquad \qquad \qquad \qquad \qquad \qquad  \\
&\quad\quad \quad\quad \!\!= \!\int_{\partial \Omega} \!\bigg\{\frac{1}{(2\pi)^{n\!-\!1}} \!\int_{{\Bbb R}^{n\!-\!1}}\!\!e^{i\langle x'\!-\!x',\xi'\rangle}\!\bigg(\!\frac{i}{2\pi} \!\int_{\mathcal{C}}\! e^{-t \tau} (\Xi_g\! -\!\tau I)^{-1}d\tau \!\bigg)d\xi'\!\bigg\} dS(x')\nonumber \qquad \qquad \qquad  \qquad    \\
\end{align*}
where $\mathcal{C}$ is a contour around the positive real axis. Secondly, by introducing some novel techniques and by combining operator algebra method (see section 5), we get the full symbol $\sum_{m=0}^\infty \psi_{-1-m}(x',\xi',\tau)$ of pseudodifferential operator $(\Xi_g -\tau I)^{-1}$. The next key difficulty is how to  calculate the corresponding residues $\frac{i}{2\pi}\int_{\mathcal{C}} e^{-t \tau} \big( \mbox{Tr}\big( \psi_{-m}( x', \xi',\tau)\big)d\tau$ for all $m\ge 1$.
 Since \begin{align*}  &\quad\quad\quad  \quad\!\!   \sum\limits_{k=0}^\infty\! e^{-t\tau_k}\! =\!\int_{\partial \Omega}\!\bigg\{\!\frac{1}{(2\pi)^{n-1}}\! \int_{{\Bbb R}^{n-1}}\! \bigg(\frac{i}{2\pi}\int_{\mathcal{C}} \!e^{-t \tau} \!\big( \mbox{Tr}\big( \psi_{-1}( x', \xi',\tau)\big)  \\
&\quad \quad\quad\quad\, \;\; \quad \quad \quad\quad \!+\!\mbox{Tr}
\big( \psi_{-2} (\tau, x', \xi')\big)\!+\! \cdots \big) d\tau \!\bigg) d\xi'\!\bigg\} dS(x')\nonumber\,  \\
&\quad\quad \; \quad\;\quad\quad\quad\quad\! \sim \sum\limits_{m=0}^\infty\! a_m t^{-n+1+m} \!+\!\sum\limits_{l=1}^\infty \!b_l t^l \log t \quad \mbox{as}\,\, t\to 0^+, \;\quad\nonumber \end{align*}
by arduous calculations we can finally get all coefficients $a_m$, $1\le m\le n-1$, for $\Xi_g$. These coefficients $a_m$ explicitly give some important geometric information for the manifold $\Omega$.

\vskip 1.49 true cm

\section{Elastic Lam\'{e} operator and boundary conditions on Riemannian manifold}

\vskip 0.45 true cm

\noindent {\bf 2.1. Representation of elastic Lam\'{e} operator on Riemannian manifolds.}   Let $\Omega$ be an $n$-dimensional Riemannian  manifold (possibly with boundary), and let $\Omega$ be equipped with a smooth metric tensor $g =
\sum\limits_{j,k=1}^n g_{jk}\, dx_j \otimes dx_k$. Denote by $[g^{jk}]_{n\times n}$ the inverse of the matrix $[g_{jk}]_{n\times n}$ and set $|g|:= \mbox{det}\, [g_{jk}]_{n\times n}$. In
particular, $d\mbox{V}$, the volume element in $\Omega$ is locally given by $d\mbox{V} = \sqrt{|g|}\, dx_1\cdots dx_n$. By $T\Omega$ and $T^*\Omega$ we denote, respectively, the tangent and cotangent bundle on $\Omega$. Throughout, we shall also denote by $T\Omega$ global ($C^\infty$-smooth) sections in $T\Omega$ (i.e., $T\Omega \equiv C^\infty (\Omega, T\Omega)$);
similarly, $T^*\Omega \equiv C^\infty(\Omega, T^*\Omega)$. Recall first that
\begin{eqnarray} \label{9-2.1} \mbox{div} \,{X} :=\sum\limits_{j=1}^n \frac{1}{\sqrt{|g|}}\, \frac{\partial(\sqrt{|g|} \,{{X}}^j)}{\partial x_j}\quad \, \mbox{if}\;\; {{X}}=\sum\limits_{j=1}^n {{X}}^j \frac{\partial}{\partial x_j}\in T\Omega,\end{eqnarray}
and \begin{eqnarray} \label{9-2.2} \mbox{grad}\, v = \sum\limits_{j,k=1}^n \bigg(g^{jk} \frac{\partial v}{\partial x_k}\bigg)\frac{\partial}{\partial x_j}\quad \, \mbox{if}\;\; v\in C^\infty(\Omega), \end{eqnarray}
are, respectively, the usual divergence and gradient operators. Accordingly, the Laplace-Beltrami
operator $\Delta_g$ is just given by
 \begin{eqnarray} \label{18/12/22-2} \Delta_g:= \mbox{div}\; \mbox{grad} = \frac{1}{\sqrt{|g|}} \sum_{j,k=1}^n \frac{\partial}{\partial x_j} \bigg(\sqrt{|g|}\,g^{jk} \frac{\partial}{\partial x_k}\bigg). \end{eqnarray}
Next, let $\nabla$ be the associated Levi-Civita connection. For each ${X} \in T\Omega$, $\nabla X$ is the tensor of type $(0,2)$ defined by
\begin{eqnarray} \label{9-2.5} (\nabla {X})({Y},{Z}):= \langle \nabla_{{Z}} {X}, {Y}\rangle, \quad \; \forall\, {Y}, {Z}\in T\Omega. \end{eqnarray}
It is well-known that in a local coordinate system with the naturally associated frame field on the tangent bundle,
\begin{eqnarray*} \nabla_{\frac{\partial}{\partial x_k}} X = \sum\limits_{j=1}^n \big(\frac{\partial X^j}{\partial x_k} +\sum\limits_{l=1}^n \Gamma_{lk}^j X^l  \big)\frac{\partial }{\partial x_j}\quad \; \mbox{for}\;\; X=\sum\limits_{j=1}^n X^j \frac{\partial}{\partial x_j}, \end{eqnarray*}
 where $\Gamma_{lk}^j= \frac{1}{2} \sum_{m=1}^n g^{jm} \big( \frac{\partial g_{km}}{\partial x_l} +\frac{\partial g_{lm}}{\partial x_k} -\frac{\partial g_{lk}}{\partial x_m}\big)$ are the Christoffel symbols associated with the metric $g$ (see, for example, \cite{Ta2}). If we denote \begin{eqnarray*} {X^j}_{;k}= \frac{\partial X^j}{\partial x_k} +\sum\limits_{l=1}^n \Gamma_{lk}^j X^l,\end{eqnarray*}
 then     \begin{eqnarray*} \label{18/10/29}  \nabla_Y X = \sum\limits_{j,k=1}^{n} Y^k {X^j}_{;k} \,\frac{\partial}{\partial x_j} \;\; \mbox{for}\,\; Y=\sum\limits_{k=1}^n Y^{k} \frac{\partial}{\partial x_k}.\end{eqnarray*}
  The symmetric part of $\nabla {X}$ is $\mbox{Def}\, {X}$, the deformation of ${X}$, i.e.,
\begin{eqnarray}\; (\mbox{Def}\; {X})({Y},{Z}) \!=\!\frac{1}{2} \{ \langle \nabla_{{Y}} {X}, {Z}\rangle \!+\!\langle \nabla_{{Z}} {X}, {Y}\rangle \}, \quad \forall\, {Y}, {Z}\in T\Omega\end{eqnarray}
(whereas the antisymmetric part of $\nabla\, {X}$ is simply $d{X}$, i.e.,
\begin{eqnarray*} d{X}({Y}, {Z}) = \frac{1}{2} \{ \langle \nabla_{{Y}} {X}, {Z}\rangle -\langle \nabla_{{Z}} {X}, {Y}\rangle \}, \quad \, \forall \,{Y}, {Z}\in T\Omega.)\end{eqnarray*}
Thus, $\mbox{Def}\; {X}$ is a symmetric tensor field of type $(0,2)$. In coordinate notation,
 \begin{eqnarray} \label{9-2.7} (\mbox{Def}\; {X})_{jk} = \frac{1}{2} ({{X}}_{j;k} +{{X}}_{k;j}), \quad \, \forall\, j,k,\end{eqnarray}
 where, ${{X}}_{k;j}:= \frac{\partial {{X}}_k}{\partial {\mathbf{x}}_j} -\sum_{l=1}^n \Gamma^{l}_{kj} {{X}}_l$ for a vector field ${X} =\sum_{j=1}^n {{X}}^j \frac{\partial}{\partial x_j}$, and $X_k=\sum_{l=1}^n g_{kl} X^l$.
 The adjoint  ${\mbox{Def}}^*$ of $\mbox{Def}$ is defined in local coordinates by $({\mbox{Def}}^* {w})^j=-\sum_{k=1}^n {{w}}^{jk}_{\;\;\;\,;k}$ for each symmetric tensor field ${w}$ of type $(0,2)$.
In particular, if $\nu\in T\Omega$ is the outward unit normal to $\partial \Omega\hookrightarrow \Omega$, then the integration by parts formula (see formula (2.16) in \cite{DMM}, or formula (12.4) of p.$\,$463 in \cite{Ta1})
\begin{eqnarray} \label{18/113-4} \int_{\Omega} \langle \mbox{Def}\; u, w\rangle dV= \int_\Omega \langle u, \mbox{Def}^* w\rangle dV+ \int_{\partial \Omega} w (\nu, u) \,dS\end{eqnarray}
holds for any $u\in T\Omega$ and any symmetric tensor field $w$ of type $(0,2)$.
The Riemann curvature tensor $\mathcal{R}$ of $\Omega$ is given by
\begin{eqnarray} \label{9-2.10} \mathcal{R}({X}, {Y}){Z} = [\nabla_{{X}}, \nabla_{{Y}}]{Z} -\nabla_{[{X},{Y}]} {Z}, \quad \, \forall\, {X}, {Y}, {Z} \in T\Omega, \end{eqnarray}
where $[{X}, {Y}] := {X}{Y} - {Y} {X}$ is the usual commutator bracket. It is convenient to change this
into a $(0, 4)$-tensor by setting
\begin{eqnarray*} \mathcal{R}({X}, {Y}, {Z}, {W}) := \langle \mathcal{R}({X}, {Y}){Z}, {W}\rangle, \quad \; \forall\, {X}, {Y}, {Z}, {W} \in T\Omega.
\end{eqnarray*}
In other words, in a local coordinate system such as that discussed above,
 \begin{eqnarray*} &&R_{jklm}= \bigg\langle \mathcal{R}\left(\frac{\partial}{\partial x_l}, \frac{\partial}{\partial x_m}\right) \frac{\partial}{\partial x_k}, \frac{\partial}{\partial x_j}\bigg\rangle.\end{eqnarray*}
The Ricci curvature $\mbox{Ric}$ on $\Omega$ is a $(0, 2)$-tensor defined as a contraction of $\mathcal{R}$:
\begin{eqnarray*} \mbox{Ric} ({X}, {Y} )\!:=\! \sum_{j=1}^n \!\bigg\langle \!\!\mathcal{R}\bigg(\!\frac{\partial }{\partial x_j}, {Y}\!\bigg) {X}, \frac{\partial }{\partial x_j}\!\bigg\rangle\! =\! \sum_{j=1}^n\! \bigg\langle \!\!\mathcal{R}\bigg(\!{Y},\frac{\partial }{\partial x_j}\!\bigg)\!\frac{\partial }{\partial x_j}, {X}\!\!\bigg\rangle, \; \forall\, {X}, {Y}\! \in T\Omega.\end{eqnarray*}
  That is,  \begin{eqnarray}\label{19.10.3-2}   R_{jk} = \sum\limits_{l=1}^n R^l_{jlk}=\sum\limits_{l,m=1}^n g^{lm}R_{ljmk}.\end{eqnarray}
 Note that \begin{eqnarray} \label{19.6.25-1} R^j_{klm}=\frac{\partial \Gamma^j_{km}}{\partial x_l}- \frac{\partial \Gamma^j_{kl}}{\partial x_m} +\sum\limits_{s=1}^n\Gamma^j_{sl} \Gamma^s_{km}-\sum\limits_{s=1}^n \Gamma^j_{sm}\Gamma^s_{kl}.\end{eqnarray}

 \vskip 0.27 true cm

 \noindent{\bf Lemma 2.1.1.} \ {\it On a Riemannian manifold $\Omega$, modeling a homogeneous, linear, isotropic, elastic medium, the Lam\'{e} operator   ${\mathcal{L}}_g$ is given by}
  \begin{eqnarray} \label{18/11/4-4} \;\;\; \quad\; {\mathcal{L}}_g u \!\!& \!\!\!=\!\!&\!\!\! -\!\mu \, \nabla^* \nabla u \!+\!(\lambda\!+\!\mu)\, \mbox{grad}\,\mbox{div}\, u \!+\! \mu\,\mbox{Ric} (u) \quad \mbox{for}\;\, u\!=\!\sum\limits_{k=1}^n\! u^k\!\frac{\partial }{\partial x_k}, \end{eqnarray}
 {\it where} $-\nabla^*\nabla u$ {\it is the Bochner Laplacian of $u$ defined by}  \begin{eqnarray} \label{19.9.30-1} && -\!\nabla^*\nabla u\!=\! \sum\limits_{j=1}^n\!\bigg\{\! \Delta_g u^j  \!+\!2 \!\sum\limits_{k,s,l=1}^n\! \!\!g^{kl} \Gamma_{sk}^j \frac{\partial u^s}{\partial x_l} \!+\!\!
   \sum_{k,s,l=1}^n \!\!\bigg(\!g^{kl} \frac{\partial \Gamma^j_{sl}}{\partial x_k}\\
    &&\qquad \qquad \;\; \!+\!\sum\limits_{h=1}^n \!g^{kl} \Gamma_{hl}^j \Gamma_{sk}^h \!-\!\sum\limits_{h=1}^n\! g^{kl} \Gamma_{sh}^j \Gamma_{kl}^h \!\bigg)u^s\!\bigg\}\!\frac{\partial}{\partial x_j},\! \!\nonumber \end{eqnarray}
   and  \begin{eqnarray} \label{18/11/25} \mbox{Ric} (u)= \sum\limits_{j=1}^n  \big(\sum\limits_{l=1}^nR_l^j u^l\big)\frac{\partial }{\partial x_j}. \end{eqnarray}
  {\it In particular, ${\mathcal{L}}_g$ is strongly elliptic, formally self-adjoint.
  If ${u}\in T\Omega$ denotes the displacement, natural boundary conditions for ${\mathcal{L}}_g$ include prescribing $u\big|_{\partial \Omega}$, Dirichlet type, and}
  \begin{eqnarray}\label{18/11/6;5} \mbox{traction} \; {u}:= 2\mu \,(\mbox{Def}\;{u})^\# {\nu} +\lambda (\mbox{div}\; {u}){\nu} \quad \mbox{on}\; \, \partial \Omega, \end{eqnarray}
  {\it Neumann type.}

\vskip 0.16 true cm
 Here, we identity $(\mbox{Def}\,u)^\#{\nu}$ with the vector field uniquely determined by the requirement that $\langle (\mbox{Def}\; {u})^\#{\nu}, {X}\rangle = (\mbox{Def}\; {u}) ({\nu}, {X})$ for each ${X}\in T\Omega$.

  \vskip 0.28 true cm

 \noindent  {\it Proof.} \  Consider the total free elastic energy functional
\begin{eqnarray}\label{18/11/3-1}  \mathcal{E}(u) = -\frac{1}{2} \int_{\Omega} E(x, \nabla u(x)) dV_x, \quad \; u\in T\Omega,\end{eqnarray}
 ignoring at the moment the displacement boundary conditions. Obviously, equilibria states correspond to minimizers of the above variational integral.
It follows from \cite{DMM} that the elastic tensor in the case of linear, isotropic, elastic medium is given by
\begin{eqnarray} \label{18/11/3-2} \;\;\;\;E(x, \nabla u(x))= \lambda (\mbox{div}\, u)^2 (x) + 2\mu \langle (\mbox{Def}\, u)(x), (\mbox{Def}\, u)(x)\rangle.\end{eqnarray}
Thus, we will be leads to considering the variational integral
\begin{eqnarray}\;\; \;\;\mathcal{E} (u)= -\frac{1}{2} \!\int_{\Omega} \big( \lambda(\mbox{div}\,u)^2 \!+\!2\mu \langle \mbox{Def}\, u, \mbox{Def}\, u\rangle \big) dV, \; \; u\in T\Omega.\end{eqnarray}
To determine the associated Euler-Lagrange equation, for an arbitrary $v\in T\Omega$, smooth and compactly supported, we compute
\begin{eqnarray} \label{18/11/3-3}\,\;\,\; \frac{d}{dt} \mathcal{E}(u\!+\!tv) \big|_{t=0}\!\!\!\!&=\!\!\!& - \int_{\Omega} \big(\lambda\,(\mbox{div}\,u)(\mbox{div}\, v) \!+\!2\mu\;\!\langle \mbox{Def}\, u, \mbox{Def}\, v\rangle \big) dV\\  \!\!\!\!&=\!\!\!& \int_{\Omega} \langle (\lambda\, \mbox{grad}\, \mbox{div} \!-\!2\mu \, \mbox{Def}^* \mbox{Def})u, v\rangle dV\nonumber\\
\!\!&& \!- \!\int_{\partial \Omega} \langle 2\mu (\mbox{Def}\, u)\nu \!+\!\lambda (\mbox{div}\, u)\nu, v\rangle dS\nonumber\end{eqnarray}
after integrating by parts, based on (\ref{18/113-4}) and the usual divergence theorem. By taking $\frac{d}{dt}\, \mathcal{E} (u+tv)\big|_{t=0}=0$,
it follows that  \begin{eqnarray} \label{9-2.8} {\mathcal{L}}_g =-2\mu \,{\mbox{Def}}^* {\mbox{Def}} +\lambda\, \mbox{grad} \,\mbox{div}.\end{eqnarray}

Furthermore, from (\ref{9-2.7}) we get \begin{eqnarray*} (\mbox{Def}\, u)_{jk}= \frac{1}{2} \big(u_{j;k} +u_{k;j}\big)\quad \mbox{for}\;\; u=\sum_{j=1}^n u^j \frac{\partial}{\partial x_j}\in T\Omega.\end{eqnarray*}
   It follows that \begin{eqnarray}\label{18/11/4-1} && \big(\mbox{Def}^*(\mbox{Def}\, u)\big)^j=-\frac{1}{2} \sum\limits_{k=1}^n \big( u^{j;k}+u^{k;j}\big)_{;k}\\
  &&\qquad\qquad \qquad\;\;\, = -\frac{1}{2}\sum\limits_{k=1}^n {u^{j;k}}_{\!;k} -\frac{1}{2}\sum\limits_{k=1}^n {u^{k;j}}_{\!;k}, \nonumber\end{eqnarray}
   where $u^{j;k}=\sum_{m,l=1}^n g^{jm} g^{kl} u_{m;l}$ is the twice raising indices by the metric tensor for $u_{m;l}$.
     The first term on the right in (\ref{18/11/4-1}) is $\frac{1}{2} \nabla^* \nabla u$, while the second term can be written as (see, for example, p.$\,$562 of \cite{Ta3} or Exercise 3 on p.$\,$305 of \cite{Ta2})
     $$-\frac{1}{2}\sum\limits_{k=1}^n \big(u^k{}_{;k}{}^{;j} +\sum\limits_{l=1}^n R^{k\,\,j}_{\,lk} u^l\big)= - \frac{1}{2}(\mbox{grad}\, \mbox{div} \, {u} +\mbox{Ric}({u}))^j.$$
   Noting that $u^{j;k}= \sum\limits_{l=1}^n g^{kl} {u^j}_{\!;l}$ and $g^{kl}_{\;\;\; ;m}=0$, we have
   \begin{align*} & {u^{j;k}}_{\!;k}=\big(\sum\limits_{l=1}^n g^{kl} {u^j}_{\!;l}\big)_{\!;k}=
       \sum\limits_{l=1}^n g^{kl} {u^j}_{\!;l;k}\\
      &=\sum_{l=1}^n g^{kl} \bigg\{ \frac{\partial^2 u^j}{\partial x_k\partial x_l} +\sum\limits_{s=1}^n \frac{\partial u^s}{\partial x_l}
      \Gamma_{sk}^j +\sum\limits_{s=1}^n \frac{\partial u^s}{\partial x_k} \Gamma_{sl}^j-\sum\limits_{s=1}^n\frac{\partial u^j}{\partial x_s}\Gamma_{kl}^s\\
      & \quad +   u^s \big(\frac{\partial \Gamma_{sl}^j}{\partial x_k} +\sum\limits_{h} \Gamma_{sk}^h \Gamma_{hl}^j -\sum\limits_{h}
      \Gamma_{sh}^j \Gamma_{kl}^h \big)\bigg\}\\
   &=\sum\limits_{l=1}^n \bigg( g^{kl} \frac{\partial^2 u^j}{\partial x_k\partial x_l} -\sum\limits_{s=1}^n g^{kl} \Gamma_{kl}^s \frac{\partial u^j}{\partial x_s}\bigg) + \sum\limits_{s,l=1}^n \bigg(g^{kl} \Gamma^j_{sk} \frac{\partial u^s}{\partial x_l} +g^{kl} \Gamma_{sl}^j \frac{\partial u^s}{\partial x_k}\bigg)\\
   & \;\;\; \;+ \sum\limits_{s,l=1}^n \bigg( g^{kl}\frac{\partial \Gamma^j_{sl}}{\partial x_k}  +\sum\limits_{h=1}^n g^{kl} \Gamma^j_{hl} \Gamma_{sk}^h - \sum\limits_{h=1}^n g^{kl} \Gamma^j_{sh} \Gamma_{kl}^h \bigg) u^s
     \end{align*}
   so that
    \begin{align*}  \frac{1}{2}\big(\nabla^*\nabla u\big)^j&= -\frac{1}{2} \sum\limits_{k=1}^n {u^{j;k}}_{;k}=-\frac{1}{2} \bigg\{ \Delta_g u^j +2 \sum\limits_{k,s,l=1}^n g^{kl} \Gamma_{sk}^j \frac{\partial u^s}{\partial x_l}
  \\ &\;\;\, + \sum_{k,s,l=1}^n \bigg(g^{kl} \frac{\partial \Gamma^j_{sl}}{\partial x_k} +\sum\limits_{h=1}^n g^{kl} \Gamma_{hl}^j \Gamma_{sk}^h -\sum\limits_{h=1}^n g^{kl} \Gamma_{sh}^j \Gamma_{kl}^h \bigg)u^s\bigg\} \end{align*}
   because $\sum\limits_{k,l=1}^n \bigg( g^{kl} \frac{\partial^2 u^j}{\partial x_k\partial x_l} -\sum\limits_{s=1}^n g^{kl} \Gamma_{kl}^s \frac{\partial u^j}{\partial x_s}\bigg)=\Delta_g u^j$.
    Combining these and (\ref{9-2.8}), we get that for any $u=\sum\limits_{j=1}^n u^j \frac{\partial}{\partial x_j}\in T\Omega$,
    \begin{eqnarray*} \label{18/11/4-3}  \quad\; {\mathcal{L}}_g u \!\!& \!\!\!=\!\!&\!\! -\mu \, \nabla^* \nabla u +(\lambda+\mu)\, \mbox{grad}\,\mbox{div}\, u + \mu\,\mbox{Ric} (u)\\
    \!\!&\!\!\! =\!\!&\!\!\mu\sum\limits_{j=1}^n \bigg\{ \Delta_g u^j  +2 \sum\limits_{k,s,l=1}^n g^{kl} \Gamma_{sk}^j \frac{\partial u^s}{\partial x_l}\\
    \!\!\!\!\!\!\!&&\!\! +
   \sum_{k,s,l=1}^n \bigg(g^{kl} \frac{\partial \Gamma^j_{sl}}{\partial x_k} +\sum\limits_{h=1}^n g^{kl} \Gamma_{hl}^j \Gamma_{sk}^h -\sum\limits_{h=1}^n g^{kl} \Gamma_{sh}^j \Gamma_{kl}^h \bigg)u^s\bigg\}\frac{\partial}{\partial x_j}\\
  \!\!\!\!\!\!\!\!&&\!\! +(\lambda+\mu)\, \mbox{grad}\,\mbox{div}\, u + \mu\,\mbox{Ric} (u).\nonumber\end{eqnarray*}
          In particular, ${\mathcal{L}}_g$ is strongly elliptic, formally self-adjoint. It follows from Theorem 3.2 of \cite{DMM} that for the displacement ${u}\in T\Omega$ (i.e., the solution of ${\mathcal{L}}_g u=0$ in $\Omega$), the boundary condition of Dirichlet type is $u\big|_{\partial \Omega}$, and
   the boundary condition of  Neumann type is  \begin{eqnarray}\label{18/11/6;5} \mbox{traction} \; {u}:= 2\mu \,(\mbox{Def}\;{u})^\# {\nu} +\lambda (\mbox{div}\; {u}){\nu} \quad \mbox{on}\; \, \partial \Omega. \end{eqnarray}
       \qed
\vskip 0.28 true cm
\noindent{\bf Remark 2.1.2.}$\,$  {\it For an $n$-dimensional Riemannian manifold with metric $g$, the elastic ``strain tensor'' $\epsilon_{ml}$ is
 \begin{eqnarray} \label{18/11/6;2} \epsilon_{ml}=\frac{1}{2} \sum\limits_{s=1}^n \big(g_{ms} \,\frac{\partial u^s}{\partial x_l} +g_{sl} \,\frac{\partial u^s}{\partial x_m} + u^s \, \frac{g_{ml}}{\partial x_s}\big).\end{eqnarray}
   This can be written as \begin{eqnarray*} \epsilon_{ml}&=& \frac{1}{2} \sum\limits_{s=1}^n \bigg( g_{ms}\big(\frac{\partial u^s}{\partial x_l} +\sum\limits_{r=1}^n u^r \Gamma_{rl}^s\big) +g_{sl}\big( \frac{\partial u^s}{\partial x_m} + \sum\limits_{r=1}^n u^r\Gamma_{rm}^s\big)\bigg)\\
 &=& \frac{1}{2}\sum\limits_{s=1}^n \big(g_{ms}\, {u^s}_{;l} + g_{sl}\, {u^s}_{;m}\big)\nonumber\end{eqnarray*}
 because of the property $g_{ml;s}= \frac{\partial g_{ml}}{\partial x_s} -\sum\limits_{r=1}^n \big(\Gamma_{sl}^r g_{mr} + \Gamma_{sm}^{r} g_{rl}\big)=0$ on a Riemannian manifold.
  In physics, the ``stress tensor'' $\sigma^{jk}$ has a linear relationship with the strain tensor at each point, i.e.,  \begin{eqnarray} \label{18/11/6;1} \frac{\sigma^{jk}}{\sqrt{|g|}} =\sum\limits_{l,m=1}^n E^{jklm} \epsilon_{lm},\end{eqnarray}
where the elastic coefficient tensor $E^{jklm}$ of an isotropic elastic medium has the form
$$ E^{jklm} =\lambda g^{jk} g^{lm} +\mu g^{jl} g^{km} +\mu g^{jm} g^{kl},$$
 and $\lambda$ and $\mu$ are two Lam\'{e} constants. Thus,  \begin{eqnarray*} \label{18/11/3-6} \frac{\sigma^{jk}}{\sqrt{|g|}} =\sum\limits_{l,m=1}^n E^{jklm}\epsilon_{lm}= \sum\limits_{l=1}^n \big(\lambda g^{jk} {u^l}_{;l} +\mu g^{kl} {u^j}_{;l} +\mu g^{jl} {u^k}_{;l}\big).\end{eqnarray*}
It follows that the elastic Lam\'{e} equations on Riemannian manifolds can be locally represented by
 \begin{eqnarray} \label{18/11/3-9} \quad \quad  \quad ({\mathcal{L}}_g u)^j\!:=\! \sum\limits_{k=1}^n\!\bigg(\!\frac{\sigma^{jk}}{\sqrt{|g|}} \!\bigg)_{;k} \!=\!\sum\limits_{k,l=1}^n \!\big(\lambda g^{jk} {u^l}_{;l;k}\! +\!\mu g^{kl} {u^j}_{;l;k} \!+\!\mu g^{jl} {u^k}_{;l;k}\big),\end{eqnarray}
  which is just another equivalent local representation of (\ref{18/11/4-4}) for the operator ${\mathcal{L}}_g$.}

\vskip 0.39 true cm

\noindent {\bf 2.2. Neumann boundary condition in boundary normal coordinate system}.\   In what follows, we will simply denote  matrix
 \[ \begin{bmatrix} a_{11} & \cdots & a_{1n}\\
 \cdots & \cdots & \cdots \\
 a_{n1}& \cdots  & a_{nn}\end{bmatrix}\]  by $\big[a_{jk}\big]_{n\times n}$ (i.e., a matrix always is denoted in square brackets, $j$ always denotes the $j$-th row, and $k$ the $k$-th column). According to (\ref{18/11/6;5}), we have
 \begin{eqnarray} \label{9-3.15} \!\!\! \!\!\!\!\!\!  && \mbox{traction}\, {u}
:= 2\mu\, (\mbox{Def}\, {u})^\#{\nu} +{\lambda} (\mbox{div}\,u){\nu}\quad\quad \qquad \qquad \qquad \\
\!\!\! \!\!\!\!\!\!  &&   = {\mu}\, \big[{{u}}^{j}_{\;\,;k} +{{u}}_{k}^{\;\,;j}\big]_{n\times n} \,{\nu} +\lambda \big(\sum_{k=1}^n \frac{\partial u^k}{\partial x_k} +\sum_{k,l=1}^n \Gamma_{kl}^l u^k\big){\nu}
\nonumber\\
\!\!\! \!\!\!\!\!\!  && = \mu \left[  \frac{\partial u^j}{\partial x_k}+ \sum\limits_{l=1}^n\Gamma_{kl}^j u^l+ \sum\limits_{s,l=1}^n g^{js} g_{kl}
\Big( \frac{\partial u^l} {\partial x_s} + \sum\limits_{m=1}^n \Gamma^l_{sm} u^m\Big) \right]_{n\times n}
\left[\begin{array}{cc}
{{\nu}}_1\\ \vdots\\ {{\nu}}_n
\end{array}\right] \rule{0ex}{0.01em}\nonumber
\\
\!\!\! \!\!\!\!\!\!  &&\quad + \lambda \big(\sum_{k=1}^n \frac{\partial u^k}{\partial x_k} +\sum_{k,l=1}^n \Gamma_{kl}^l u^k\big)\left[\begin{array}{cc}
{{\nu}}_1\\ \vdots\\ {{\nu}}_n
\end{array}\right] \rule{0ex}{0.01em}\nonumber\\
 \!\!\! \!\!\!\!\!\!  &&  =\mu \left[\begin{array} {cc} \sum_{k=1}^n\Big( \frac{\partial u^1}{\partial x_k}+\sum\limits_{l=1}^n \Gamma_{kl}^1 u^l
 + \sum\limits_{s,l=1}^n g^{1s} g_{kl} (\frac{\partial u^l}{\partial x_s}
 +\sum\limits_{m=1}^n \Gamma_{sm}^l u^m\big)\Big) \nu_k\\
 \vdots\\
 \sum_{k=1}^n\Big( \frac{\partial u^n}{\partial x_k}+\sum\limits_{l=1}^n \Gamma_{kl}^n u^l
 + \sum\limits_{s,l=1}^n g^{ns} g_{kl} (\frac{\partial u^l}{\partial x_s}
 +\sum\limits_{m=1}^n \Gamma_{sm}^l u^m\big)\Big) \nu_k\end{array}\right] \rule{0ex}{0.01em}\nonumber\\
  \!\!\! \!\!\!\!\!\!  && \quad  + \lambda \big(\sum_{k=1}^n \frac{\partial u^k}{\partial x_k} +\sum_{k,l=1}^n \Gamma_{kl}^l u^k\big) \left[\begin{array}{cc}
{{\nu}}_1\\ \vdots\\ {{\nu}}_n
\end{array}\right] \rule{0ex}{0.01em}. \nonumber\end{eqnarray}
 For convenience, we will often use the following relations:
\[  \sum_{l=1}^n\Gamma_{kl}^l=\frac{1}{2}\sum_{l,j=1}^n g^{lj} \frac{\partial g_{lj}}{\partial x_k}  =\frac{1}{\sqrt{|g|}}\,\frac{\partial \sqrt{|g|}}{\partial x_k} = \frac{1}{2|g|} \, \frac{\partial |g|}{\partial x_k} =\frac{\partial \log \sqrt{|g|}}{\partial x_k}.\]

In order to describe the Dirichlet-to-Neumann map for the Lam\'{e} operator, we first recall the construction of usual geodesic coordinates with respect to the
boundary (see p.$\,$1101 of \cite{LU}). For each $x'\in \partial \Omega$, let $r_{x'}: [0, \epsilon)\to \bar \Omega$ denote the unit-speed geodesic
starting at $x'$ and normal to $\partial \Omega$. If $x':=\{x_1, \cdots, x_{n-1}\}$  are any local coordinates for
$\partial \Omega$ near $x_0\in \partial \Omega$, we can extend them smoothly to functions on a neighborhood
of $x_0$ in $\Omega$ by letting them be constant along each normal geodesic $r_{x'}$. If we then
define $x_n$ to be the parameter along each $r_{x'}$, it follows easily that $\{x_1, \cdots, x_{n}\}$
form coordinates for $\Omega$ in some neighborhood of $x_0$, which we call the boundary
normal coordinates determined by $\{x_1, \cdots, x_{n-1}\}$. In these coordinates $x_n>0$ in
 $\Omega$, and $\partial \Omega$ is locally characterized by $x_n= 0$. A standard computation shows
that the metric $g$ on $\bar \Omega$ then has the form
 (see p.$\,$1101 of \cite{LU} or p.$\,$532 of \cite{Ta2})
\begin{eqnarray}\!\!\!\!\!\!\!\!\!\!\!& \!\!\!\!\!\!\!\!\! \label{18/a-1} \;\qquad \quad\qquad \qquad \qquad \qquad \quad\;\\
\!\!\!\!\!\!& \big[g_{jk} (x',x_{n}) \big]_{n\times n} \!=\! \begin{bmatrix}
 g_{11} (x',x_{n}) & g_{12} (x',x_{n}) &\! \cdots \!& g_{1,n\!-\!1} (x',x_{n}) & 0\\
 \cdots\cdots& \cdots\cdots & \!\cdots \!&\cdots\cdots  & \cdots\\
 g_{n\!-\!1,1} (x',x_{n})  & g_{n\!-\!1,2} (x',x_{n}) &\! \cdots\! & g_{n\!-\!1,n\!-\!1} (x',x_{n}) & 0\\
 0& 0& 0& 0&1 \end{bmatrix}.  \nonumber \end{eqnarray}
       Furthermore, we can take a geodesic normal coordinate system for $(\partial \Omega, g|_{\partial \Omega})$ centered at $x_0=0$, with respect to $e_1, \cdots, e_{n-1}$, where  $e_1, \cdots, e_{n-1}$ are the principal curvature vectors. As Riemann showed, one has (see p.$\,$555 of \cite{Ta2}, or \cite{Spi2})
            \begin{eqnarray} \label{18/7/14/1} \begin{split}& g_{jk}(x_0)\!=\! \delta_{jk}, \; \; \frac{\partial g_{jk}}{\partial x_l}(x_0)
 \!=\!0  \,\;  \mbox{for all} \,\; 1\le j,k,l \le n\!-\!1,\\
 & \frac{1}{2}\,\frac{\partial g_{jk}}{\partial x_n} (x_0) =\kappa_k\delta_{jk}  \,\;  \mbox{for all} \;\, 1\le j,k \le n\!-\!1,\end{split}
 \end{eqnarray}
 where  \[ \delta_{jk}= \left\{\begin{array}{ll} 1 \;\; &\mbox{for}\;\, j=k\\
 0 \;\; &\mbox{for} \,\; j\ne k \end{array} \right.\]
 is the standard Kronecker symbol, and  $\kappa_1,\cdots, \kappa_{n-1}$ are the principal curvatures of $\partial \Omega$ at point $x_0=0$.
Under this normal coordinates, we take $-{\nu}(x)=[0,\cdots, 0,1]^t$, where $A^t$ denotes the transpose of the matrix (or vector) $A$.
By (\ref{18/a-1}) we immediately see that the inverse of metric tensor $g$ in the boundary normal coordinates has form:
    \begin{eqnarray*} g^{-1}(x',x_n) =\begin{bmatrix} g^{11}(x', x_n) & \cdots & g^{1,n-1} (x', x_n)& 0 \\
    \cdots & \cdots  & \cdots  & \cdots\\
    g^{n-1, 1}(x',x_n) & \cdots & g^{n-1,n-1}(x',x_n)& 0\\
    0&\cdots \cdots &0 &1\end{bmatrix}. \end{eqnarray*}
 Then, it follows from (\ref{9-3.15}) that
\begin{eqnarray} && \label{9.3-13} \\
&& \!\!\!\!\! - \mbox{traction}\, u =\mu \!\begin{bmatrix} \frac{\partial u^1}{\partial x_n}+\sum\limits_{l=1}^n \Gamma_{nl}^1 u^l
 + \sum\limits_{s,l=1}^n g^{1s} g_{nl} (\frac{\partial u^l}{\partial x_s}
 +\sum\limits_{m=1}^n \Gamma_{sm}^l u^m\big)\\  \vdots
\\  \frac{\partial u^n}{\partial x_n}+\sum\limits_{l=1}^n \Gamma_{nl}^n u^l
 + \sum\limits_{s,l=1}^n g^{ns} g_{nl} (\frac{\partial u^l}{\partial x_s}
 +\sum\limits_{m=1}^n \Gamma_{sm}^l u^m\big) \!\end{bmatrix} \nonumber \\
 && \quad +\lambda \big(\sum\limits_{k=1}^n \frac{\partial u^k}{\partial x_k} \!+\!\sum\limits_{k,l=1}^n\! \Gamma_{kl}^l u^k\big)\!\begin{bmatrix} 0\\ \vdots \\
 0\\ 1\end{bmatrix}
 = \left\{\mu\begin{bmatrix} h_{11} & h_{12} &\cdots & h_{1n} \\
h_{21} & h_{22} &\cdots & h_{2n}\\
\cdots & \cdots & \cdots &\cdots \\
h_{n1} &h_{n2} & \cdots & h_{nn}\end{bmatrix}\right.
  \nonumber \\
&&\quad \left. +\lambda \begin{bmatrix} 0& 0&\cdots & 0\\
\cdots & \cdots & \cdots &\cdots\\
0 & 0&\cdots &0 \\
 \frac{\partial}{\partial x_1}+\sum\limits_{l=1}^n \Gamma_{1l}^l & \frac{\partial}{\partial x_2}+\sum\limits_{l=1}^n \Gamma_{2l}^l& \cdots & \frac{\partial}{\partial x_n}+\sum\limits_{l=1}^n \Gamma_{nl}^l\end{bmatrix}\right\} \begin{small}\begin{bmatrix}u^1\\ u^2 \\
\vdots\\ u^n \end{bmatrix}\end{small}, \nonumber\end{eqnarray}
    where \begin{eqnarray} \label{18/10/14-1} \quad\;\;\;\quad \;\, \;h_{jk}\!=\! \delta_{jk}\frac{\partial}{\partial x_n}\!+\!\sum\limits_{s=1}^n g^{js} g_{nk} \frac{\partial}{\partial x_s}
 \!+\!  \Gamma_{nk}^j \!+\!\sum_{s,l=1}^n g^{js} g_{nl} \Gamma_{sk}^l, \quad 1\le j,k\le n.\end{eqnarray}

 \vskip 0.29 true cm

\noindent{\bf Remark 2.2.1.} \  {\it The first term of (\ref{9.3-13}) can also be obtained by another way. In fact, according to the definition of} $(\mbox{Def}\, u)^\#\nu$, {\it we have that}
  \begin{eqnarray} \label{2022.5.24-1}\langle (\mbox{Def}\; u)^\#\nu, X\rangle =(\mbox{Def}\; u) (\nu, X)  \quad \; \mbox{for all} \;\; X\in T\Omega,\end{eqnarray}
  {\it so that}
  \begin{align*} (\mbox{Def}\; u) (\nu, X)&=\frac{1}{2} \big( \langle \nabla_{\nu} u, X\rangle +\langle \nabla_X u, \nu\rangle \big) \\
  &= \frac{1}{2} \bigg( \Big\langle \sum\limits_{k,j=1}^n \nu^k \big(\frac{\partial u^j}{\partial x_k} +\sum\limits_{l=1}^n  \Gamma_{lk}^ju^l \big) \frac{\partial }{\partial x_j}, \sum_{m=1}^n X^m \frac{\partial }{\partial x_m}\Big\rangle \\
  &\quad + \Big\langle \sum\limits_{k,j=1}^n X^k \big(\frac{\partial u^j}{\partial x_k} +\sum\limits_{l=1}^n  \Gamma_{lk}^j u^l\big) \frac{\partial }{\partial x_j}, \sum_{m=1}^n \nu^m \frac{\partial }{\partial x_m}\Big\rangle\bigg)\\
  &= \frac{1}{2}\sum\limits_{j,k,m=1}^n \bigg(  \nu^k \big( \frac{\partial u^j}{\partial x_k} +\sum\limits_{l=1}^n \Gamma_{lk}^j u^l \big) g_{jm} X^m \\
  & \;\;\;\, +
  X^k \big( \frac{\partial u^j}{\partial x_k} +\sum\limits_{l=1}^n  \Gamma_{lk}^j u^l \big) g_{jm} \nu^m \bigg).\end{align*}
{\it Taking $-\nu (x)=[0,\cdots, 0,1]^t$, we get that for any $[X^1, X^2,\cdots, X^n]^t\in {\Bbb R}^n$},
\begin{eqnarray*} (\mbox{Def}\; u) (-\nu, X) = \frac{1}{2} \sum\limits_{m=1}^n  \bigg(\sum\limits_{j=1}^n\big(\frac{\partial u^j}{\partial x_n} +\sum_{l=1}^n \Gamma_{ln}^j u^l\big)g_{jm} +\frac{\partial u^n}{\partial x_m} +\sum_{l=1}^n  \Gamma_{lm}^nu^l\bigg)X^m.\end{eqnarray*}
{\it It follows from (\ref{2022.5.24-1}) that} \begin{eqnarray*}-\begin{small} 2\mu \,\big(\mbox{Def}\; u\big)^{\#}\,\begin{bmatrix} 0\\
 \vdots \\ 0\\
 1\end{bmatrix} = \mu  \begin{bmatrix} \frac{\partial u^1}{\partial x_n} +\sum_{l=1}^n  \Gamma_{ln}^1u^l + \sum\limits_{s=1}^n g^{1s}\big(\frac{\partial u^n}{\partial x_s} +\sum_{m=1}^n  \Gamma_{ms}^n u^m\big)\\
\vdots \\
 \frac{\partial u^n}{\partial x_n} +\sum_{l=1}^n  \Gamma_{ln}^n u^l + \sum\limits_{s=1}^n g^{ns}\big(\frac{\partial u^n}{\partial x_s} +\sum_{m=1}^n  \Gamma_{ms}^n u^m\big)\end{bmatrix},\end{small}\end{eqnarray*}
{\it which is just the first term after the first equal sign in (\ref{9.3-13}) because $g_{nn}=1$ and $g_{nl}=0$ for $1\le l\le n-1$.}

 \vskip 0.42 true cm

     Note that under boundary normal coordinates,  we have \begin{eqnarray} \label{19.10.5-1}\Gamma_{nk}^n &=&\frac{1}{2} \sum_{m=1}^n g^{nm}\bigg(\frac{\partial g_{nm}}{\partial x_k} +\frac{\partial g_{km}}{\partial x_n} -\frac{\partial g_{nk}}{\partial x_m}\bigg) \\ &=&\frac{1}{2} \bigg(\frac{\partial g_{nn}}{\partial x_k}+\frac{\partial g_{kn}}{\partial x_n} -\frac{\partial g_{nk}}{\partial x_n}\bigg)=0,\nonumber\end{eqnarray}
   \begin{eqnarray} \Gamma_{nn}^l =\frac{1}{2}\sum\limits_{m=1}^n g^{lm}\bigg(\frac{\partial g_{nm}}{\partial x_n} +\frac{\partial g_{nm}}{\partial x_n} -\frac{\partial g_{nn}}{\partial x_m}\bigg) =0\end{eqnarray}
 and  \begin{eqnarray}\label{18/10/14-2}  &&  \Gamma_{nk}^j +\sum_{s,l=1}^n g^{js} g_{nl} \Gamma_{sk}^l = \Gamma_{nk}^j +\sum_{s=1}^n g^{js} \Gamma_{sk}^n\\
 && \;\;\;\;\quad =\frac{1}{2}  \sum_{s=1}^{n} g^{js} \big( \frac{
  \partial g_{ks}}{\partial x_n} +\frac{\partial g_{ns}}{\partial x_k} -\frac{\partial g_{nk}}{\partial x_s}\big) \nonumber \\
 && \;\; \;\;\;\quad \,\;\;\; + \sum_{s=1}^n g^{js} \Big( \frac{1}{2}  \sum_{m=1}^n g^{nm} \big( \frac{\partial g_{km}}{\partial x_s} +\frac{\partial g_{sm}}{\partial x_k} -\frac{\partial g_{sk}}{\partial x_m}\big)\Big)\nonumber\\
 && \;\;\quad \quad =\frac{1}{2}\sum_{s=1}^{n} g^{js} \Big(\frac{\partial g_{ks}}{\partial x_n} +\frac{\partial g_{ns}}{\partial x_k} -\frac{\partial g_{nk}}{\partial x_s}\Big) +\frac{1}{2} \sum\limits_{s=1}^n g^{js}  \Big(\frac{\partial g_{kn}}{\partial x_s} +\frac{\partial g_{sn}}{\partial x_k} -\frac{\partial g_{sk}}{\partial x_n}\Big) =0.\nonumber\end{eqnarray}
 Combining these and (\ref{18/10/14-1}) we get
 \begin{eqnarray*}&& h_{jk} = \delta_{jk} \frac{\partial}{\partial x_n}, \quad  h_{jn}= \sum\limits_{l=1}^n g^{jl}\frac{\partial} {\partial x_l}, \quad  h_{nk}= 0 \quad  \mbox{for all}\;\; 1\le j,k\le n-1,\\
&&  h_{nn} =2\frac{\partial}{\partial x_n}. \end{eqnarray*}
 In what follows, we will let Greek indices run from $1$ to $n-1$, Roman indices from $1$ to $n$. Then $-\mbox{traction}\; u$ is reduced to
\begin{eqnarray*} \label{9.3-15}\!\!\!\! \!\!\!\!\!\!\!\!\! &&-\mbox{traction}\, u\\  [2.5mm] \!\! \!\!\!\!\!\!\!\! &&
= \begin{bmatrix} \mu \delta_{11} \frac{\partial}{\partial x_n} & \cdots &\mu \delta_{1,n-1}\frac{\partial}{\partial x_n} &\mu \sum\limits_{\alpha}g^{1\alpha} \frac{\partial}{\partial x_\alpha}\\
\cdots& \cdots& \cdots & \cdots \\
\mu \delta_{n-1,1} \frac{\partial}{\partial x_n} &\cdots &\mu \delta_{n-1,n-1}\frac{\partial}{\partial x_n} &\mu \sum\limits_{\alpha}g^{n-1,\alpha} \frac{\partial}{\partial x_\alpha} \\
  \lambda(\frac{\partial}{\partial x_1}\!+\!\sum\limits_{\alpha}\Gamma_{1\alpha}^\alpha) &  \cdots & \lambda(\frac{\partial}{\partial x_{n\!-\!1}}\!+\!\sum\limits_{\alpha }\Gamma_{(n\!-\!1)\alpha}^\alpha) & (\lambda\!+\!2\mu) \frac{\partial}{\partial x_n} \!+\!\lambda \sum\limits_{\alpha} \Gamma_{n\alpha}^\alpha \end{bmatrix}\!\!\begin{bmatrix}u^1\\
  \vdots\\
 u^{n\!-\!1}\\
 u^n\end{bmatrix}.\end{eqnarray*}
Throughout this paper, we will denote by
\begin{equation*}
  \begin{bmatrix}
\begin{BMAT}(@, 15pt, 15pt){c.c}{c.c}
   [a_{jk}]_{(n-1)\times (n-1)}&[b_j]_{(n-1)\times 1}\\
   [c_k]_{1\times (n-1)}& d
  \end{BMAT}
\end{bmatrix}
\end{equation*} the block matrix
\[ \begin{bmatrix}
\begin{BMAT}(@, 15pt, 15pt){ccc.c}{ccc.c}
    a_{11} & \cdots & a_{1,n-1} & b_{1}\\
 \cdots & \cdots &\cdots & \vdots \\
   a_{n-1,1} & \cdots & a_{n-1,n-1} & b_{n-1}\\
    c_{1} & \cdots & c_{n-1} & d
  \end{BMAT}
\end{bmatrix},  \]
 where $\big[a_{jk}\big]_{(n-1)\times (n-1)}$, $\big[ b_j\big]_{(n-1)\times 1}$ and $\big[c_k\big]_{1\times (n-1)}$ are the $(n-1)\times (n-1)$ matrix \begin{eqnarray*}\begin{bmatrix} a_{11}& a_{12} & \cdots &a_{1,n-1}\\
a_{21}& a_{22} & \cdots &a_{2,n-1}\\
\cdots & \cdots& \cdots & \cdots\\
a_{n-1,1}& a_{n-2,2} & \cdots &a_{n-1,n-1}\end{bmatrix}, \end{eqnarray*}
the $(n-1)\times 1$ matrix
\begin{eqnarray*}\begin{bmatrix} b_1 \\  b_2\\
\vdots \\ b_{n-1}\end{bmatrix} \end{eqnarray*}
and the $1\times (n-1)$ matrix \begin{eqnarray*}\begin{bmatrix} c_1 & c_2 &\cdots & c_{n-1} \end{bmatrix}, \end{eqnarray*} respectively. We also simply denote by $\Xi_g u$ the $(\mbox{traction}\, u)\big|_{\partial \Omega}$ as before.
 Thus we have obtained the following lemma:

\vskip 0.20 true cm

  \noindent{\bf Lemma 2.2.2.} \ {\it Let} ${\mathcal{L}}_g =-\mu \nabla^*\nabla +(\mu+\lambda)\mbox{grad}\, \mbox{div} +\mu\, \mbox{Ric}$ {\it be the elastic Lam\'{e} operator, and let $u\in T\Omega$ be the displacement.  Then, in boundary normal coordinates, the Neumann  boundary condition $(\mbox{traction}\; u)\big|_{\partial \Omega}$ (with respect to the outward unit normal vector of $\partial \Omega$) can be simply written as}
 \begin{equation*}    (\mbox{traction}\, u)\big|_{\partial \Omega}\!=\! -\!\begin{bmatrix}\!
\begin{BMAT}(@, 0pt, 0pt){c.c}{c.c}
  \bigg[ \mu \delta_{jk}\frac{\partial}{\partial x_n} \bigg]_{(n-1)\times (n-1)}& \bigg[\mu \sum\limits_{\alpha}g^{j\alpha}\frac{\partial}{\partial x_\alpha}\bigg]_{(n-1)\times 1}\\
   \bigg[\!\lambda(\frac{\partial}{\partial x_k}\!+\!\sum_{\alpha}\!\Gamma_{k\alpha}^\alpha)\bigg]_{1\times (n\!-\!1)} & (\lambda\!+\!2\mu) \frac{\partial}{\partial x_n} \!+\!\lambda \sum_{\alpha}\!\Gamma_{n\alpha}^\alpha
  \end{BMAT}\!
\end{bmatrix}\!\begin{bmatrix}u^1\\
  \vdots\\
 u^{n\!-\!1}\\
 u^n\end{bmatrix}.\end{equation*}
{\it Or equivalently},  \begin{eqnarray} &&  \label{18-9/22-4}\!\!\!  \end{eqnarray}  \begin{eqnarray}
 \!\!\!\!\! &&
 \Xi_g (u\big|_{\partial \Omega})=\left\{ \left[
 \begin{BMAT}(@, 6pt, 6pt){c.c}{c.c}\big[\mu \delta_{jk} \big]_{(n-1)\times (n-1)}&\big[0\big]_{(n-1)\times 1}\\
\big[0\big]_{1\times (n-1)}& \lambda+2\mu\end{BMAT}\right] \bigg[-\frac{\partial }{\partial x_n}I_n\bigg]\qquad  \qquad \qquad \qquad\quad  \qquad \qquad \right.\nonumber\\
\!\!\!\!\! &&  \left. \qquad  \quad \quad  \quad - \left[ \begin{BMAT}(@, 6pt, 6pt){c.c}{c.c} [0]_{(n-1)\times (n-1)}& \bigg[\mu\sum\limits_{\alpha} g^{j\alpha} \frac{\partial}{\partial x_\alpha}\bigg]_{(n-1)\times 1}\\
\bigg[\!\lambda(\frac{\partial} {\partial x_k} \!+\!\sum\limits_{\alpha} \Gamma_{k\alpha}^\alpha\big)\bigg]_{1\times (n\!-\!1)} & \lambda\sum\limits_{\alpha}\Gamma_{n\alpha}^\alpha\end{BMAT}\!\right]\!\right\}\begin{bmatrix}u^1\\
  \vdots\\
 u^{n\!-\!1}\\
 u^n\end{bmatrix}\nonumber  \end{eqnarray}
 {\it for $u$ satisfies \begin{eqnarray*} \left\{ \begin{array}{ll}{\mathcal{L}}_g u=0  &\mbox{in}\; \, \Omega,\\
  u=f  &\mbox{on}\;\, \partial \Omega.\end{array}\right.\end{eqnarray*}}
 We call (\ref{18-9/22-4}) the {\it transform relation} about operator $Q:=-\frac{\partial}{\partial x_n}\, I_n$ and $\Xi_g$. This is completely a new formula.

 \vskip 1.0 true cm

\section{Explicit expression of the Dirichlet-to-Neumann map}

\vskip 0.45 true cm

The purpose of this section is to give an expression for $-\frac{\partial }{\partial x_n} I_n$ (which is a pseudodifferential operator defined on $\partial \Omega$) by discussing the elastic Lam\'{e} equation ${\mathcal{L}}_g u=0$ in $\Omega$ and $u=f$ on $\partial \Omega$. So an explicit representation for $\Xi_g u$ will be obtained from such an expression and transform relation (\ref{18-9/22-4}).

In boundary normal coordinates, because of $\Gamma_{jn}^n=0$ and $g^{jn}=0$ for $1\le j\le n-1$, we have $$\;\mbox{div}\, u=  \sum\limits_{k=1}^n \big(\frac{\partial u^k}{\partial x_k} +\sum\limits_{\beta} \Gamma_{k\beta}^\beta u^k\big),$$ so that
\begin{eqnarray*} &&\!\!\!\!\!\!\!\!\! \mbox{grad}\, \mbox{div}\, u= \begin{bmatrix} \sum\limits_{m=1}^n g^{1m} \sum\limits_{k=1}^n \big( \frac{\partial^2 u^k}{\partial x_m\partial x_k} +\sum\limits_{\beta} \Gamma_{k\beta}^\beta \frac{\partial u^k}{\partial x_m} +\sum\limits_{\beta} \frac{\partial \Gamma_{k\beta}^\beta}{\partial x_m} u^k\big)\\
\vdots \\
 \sum\limits_{m=1}^n g^{nm} \sum\limits_{k=1}^n \big( \frac{\partial^2 u^k}{\partial x_m\partial x_k} +\sum\limits_{\beta} \Gamma_{k\beta}^\beta \frac{\partial u^k}{\partial x_m} +\sum\limits_{\beta} \frac{\partial \Gamma_{k\beta}^\beta}{\partial x_m} u^k\big)\end{bmatrix}\qquad\quad \quad \qquad \qquad  \qquad \qquad \qquad \\ [1mm]
&&\;\;\;\;\;\;\;\;\;\; \;= \begin{bmatrix} \sum\limits_{\alpha} g^{1\alpha} \sum\limits_{k=1}^n \big( \frac{\partial^2 u^k}{\partial x_\alpha\partial x_k} +\sum\limits_{\beta} \Gamma_{k\beta}^\beta \frac{\partial u^k}{\partial x_\alpha} +\sum\limits_{\beta} \frac{\partial \Gamma_{k\beta}^\beta}{\partial x_\alpha} u^k\big)\\
\vdots \\  \sum\limits_{\alpha} g^{n-1,\alpha} \sum\limits_{k=1}^n \big( \frac{\partial^2 u^k}{\partial x_\alpha\partial x_k} +\sum\limits_{\beta} \Gamma_{k\beta}^\beta \frac{\partial u^k}{\partial x_\alpha} +\sum\limits_{\beta} \frac{\partial \Gamma_{k\beta}^\beta}{\partial x_\alpha} u^k\big)\\
 \sum\limits_{k=1}^n \big(\frac{\partial^2 u^k}{\partial x_n\partial x_k} +\sum\limits_{\beta} \Gamma_{k\beta}^\beta \frac{\partial u^k}{\partial x_n} +\sum\limits_{\beta}\frac{\partial \Gamma_{k\beta}^\beta}{\partial x_n} u^k\big)\end{bmatrix}\qquad\quad   \qquad \quad\qquad \qquad \qquad \qquad \\ [1mm] \!\! \!\!\!\!\!\!\!&\!\!\!&\!\!\!\!\!\!\!\!\!\!\!\!= \!\!{ \begin{bmatrix} \!\sum\limits_{\alpha}\! g^{1\alpha}\! \big(\! \frac{\partial^2}{\partial x_\alpha\partial x_1}\! +\!\sum\limits_{\beta} \! \Gamma_{1\beta}^\beta \!\frac{\partial }{\partial x_\alpha}
 \!+\!\sum\limits_{\beta}\! \frac{\partial \Gamma_{1\beta}^\beta}{\partial x_\alpha} \!\big)& \!\cdots \!&  \sum\limits_{\alpha}\! g^{1\alpha} \!\big(\! \frac{\partial^2}{\partial x_\alpha\partial x_n}\! +\!\sum\limits_{\beta}\! \Gamma_{n\beta}^\beta\! \frac{\partial }{\partial x_\alpha}
 \!+\!\sum\limits_{\beta}\! \frac{\partial \Gamma_{n\beta}^\beta}{\partial x_\alpha}\!\big)\\
 \cdots \cdots\cdots& \!\cdots \!& \cdots\cdots\cdots \\
   \!\sum\limits_{\alpha} \!g^{n\!-\!1,\alpha} \!\big( \!\frac{\partial^2}{\partial x_\alpha\partial x_1} \!+\!\sum\limits_{\beta}\! \Gamma_{1\beta}^\beta \!\frac{\partial }{\partial x_\alpha}
 \!+\!\sum\limits_{\beta} \!\frac{\partial \Gamma_{1\beta}^\beta}{\partial x_\alpha}\!\big) & \!\cdots \!&  \sum\limits_{\alpha}\! g^{n\!-\!1,\alpha}\!\big( \! \frac{\partial^2}{\partial x_\alpha\partial x_n} \!+\!\sum\limits_{\beta}\! \Gamma_{n\beta}^\beta\! \frac{\partial }{\partial x_\alpha}
 \!+\!\sum\limits_{\beta} \frac{\partial \Gamma_{n\beta}^\beta}{\partial x_\alpha}\!\big)\\
 \!  \frac{\partial^2}{\partial x_n\partial x_1} \!+\!\sum\limits_{\beta}\!\Gamma_{1\beta}^\beta\!\frac{\partial }{\partial x_n}
 \!+\!\sum\limits_{\beta}\! \frac{\partial \Gamma_{1\beta}^\beta}{\partial x_n} &\!\cdots \!&  \frac{\partial^2}{\partial x_n^2} \!+\!\sum\limits_{\beta}\! \Gamma_{n\beta}^\beta\!\frac{\partial }{\partial x_n}
\! +\!\sum\limits_{\beta} \!\frac{\partial \Gamma_{n\beta}^\beta}{\partial x_n}\!\!\end{bmatrix}\!\!\begin{bmatrix} \!u^1\\
 \vdots \\
 u^{n\!-\!1} \\
 u^n \!\end{bmatrix}. }  \end{eqnarray*}
  According to Lemma 2.1.1, in the boundary normal coordinates, the elastic Lam\'{e} operator ${\mathcal{L}}_g u=-\mu \nabla^* \nabla u+(\mu+\lambda)\, \mbox{grad}\, \mbox{div}\,u +\mu \,\mbox{Ric}\,(u)$ can be written as
\begin{eqnarray} &&\label{9.4-1}\end{eqnarray}
 \begin{eqnarray*} \!\!\!\!\!&& \!\!\!\!\! {\mathcal{L}}_g u\!=\!\begin{small} \left\{\!\mu \bigg(\!\frac{\partial^2}{\partial x_n^2}\! +\!\sum\limits_{\beta} \!\Gamma_{n\beta}^\beta\! \frac{\partial }{\partial x_n}\! +\!\sum\limits_{\alpha, \beta}\! g^{\alpha \beta} \!\frac{\partial^2}{\partial x_\alpha\partial x_\beta}\! +\!\sum\limits_{\alpha, \beta}\!\big(\!\sum_{\gamma}g^{\alpha\beta}\! \Gamma_{\alpha \gamma}^\gamma \! +\!\frac{\partial g^{\alpha\beta}}{\partial x_\alpha}\!\big)\frac{\partial}{\partial x_\beta}\!\bigg)\!I_n \!\right.\end{small}\end{eqnarray*}
     \begin{eqnarray} \!\!\!&&   \,\, \;\;   +\mu \begin{small}\begin{bmatrix}\! 2\Gamma_{1n}^1& \cdots  & {2\Gamma_{(n-1)n}^1}_{{}_{{}_{}}} &0\\
   \cdots &\cdots & \cdots \!\!&\!\!\cdots \\
   2\Gamma_{1n}^{n-1}  & \cdots & {2\Gamma_{(n-1)n}^{n-1}}_{{}_{{}_{}}}  & 0\\
   0&\cdots & 0& 0\end{bmatrix}\end{small}\! \!\big(\frac{\partial }{\partial x_n}\! I_n\big)\!\nonumber \\
   &&\quad  + \mu\begin{small}\begin{bmatrix} 2\!\sum\limits_{\alpha,\beta}g^{\alpha\beta} {\Gamma_{1\alpha}^1 \frac{\partial }{\partial x_\beta}}_{{}_{{}_{}}} & \cdots &  2\!\sum\limits_{\alpha,\beta}g^{\alpha\beta} \Gamma_{(n-1)\alpha}^1 \frac{\partial }{\partial x_\beta}&  2\!\sum\limits_{\alpha,\beta}g^{\alpha\beta} \Gamma_{n\alpha}^1\frac{\partial }{\partial x_\beta}\!\!\\
   \cdots & \cdots& \cdots & \cdots\!\!\\
    2\!\sum\limits_{\alpha,\beta}g^{\alpha\beta} \Gamma_{1\alpha}^{n\!-\!1}\frac{\partial }{\partial x_\beta} & \cdots &  2\!\sum\limits_{\alpha,\beta}g^{\alpha\beta} \Gamma_{(n\!-\!1)\alpha}^{n\!-\!1}\frac{\partial }{\partial x_\beta} &  2\!\sum\limits_{\alpha,\beta}g^{\alpha\beta} {\Gamma_{n\alpha}^{n\!-\!1}\frac{\partial }{\partial x_\beta}}_{{}_{{}_{}}}\!\!\\
    2\!\sum\limits_{\alpha,\beta}g^{\alpha\beta} \Gamma_{1\alpha}^n \frac{\partial }{\partial x_\beta}&  \cdots & 2\!\sum\limits_{\alpha,\beta}g^{\alpha\beta} {\Gamma_{(n\!-\!1)\alpha}^n\frac{\partial }{\partial x_\beta}}_{{}_{{}_{}}} &  0\!\!\end{bmatrix}\end{small} \nonumber
    \end{eqnarray}
\begin{eqnarray} \!\!\!\! \!\!\!\!\!\!\!\! \!\!\!\!\!&&\!\!\!\!\!\! \!\!\!\!\!+\!(\!\lambda\!+\!\mu) \! {\begin{bmatrix}\!\! \sum\limits_{\alpha} \!g^{1\alpha} \big(\! \frac{\partial^2}{\partial x_\alpha\partial x_1} \!+\!\sum\limits_{\beta} \!\Gamma_{1\beta}^\beta \!\frac{\partial }{\partial x_\alpha}
\! +\!\sum\limits_{\beta} \!\frac{\partial \Gamma_{1\beta}^\beta}{\partial x_\alpha}\! \big)& \!\cdots \!&  \sum\limits_{\alpha} \!g^{1\alpha}\! \big( \! \frac{\partial^2}{\partial x_\alpha\partial x_n} \!+\!\sum\limits_{\beta} \!\Gamma_{n\beta}^\beta \!\frac{\partial }{\partial x_\alpha}
 \!+\!\sum\limits_{\beta} \!\frac{\partial \Gamma_{n\beta}^\beta}{\partial x_\alpha}\!\big)\\
 \cdots \cdots\cdots& \cdots & \cdots\cdots\cdots \\
   \sum\limits_{\alpha} \!g^{n\!-\!1,\alpha} \big( \frac{\partial^2}{\partial x_\alpha\partial x_1} \!+\!\sum\limits_{\beta} \!\Gamma_{1\beta}^\beta \frac{\partial }{\partial x_\alpha}
\! +\!\sum\limits_{\beta}\! \frac{\partial \Gamma_{1\beta}^\beta}{\partial x_\alpha}\big) & \!\cdots \!&  \sum\limits_{\alpha} \!g^{n\!-\!1,\alpha} \big( \!\frac{\partial^2}{\partial x_\alpha\partial x_n} \!+\!\sum\limits_{\beta} \!\Gamma_{n\beta}^\beta \frac{\partial }{\partial x_\alpha}
\! +\!\sum\limits_{\beta} \!\frac{\partial \Gamma_{n\beta}^\beta}{\partial x_\alpha}\big)\\
   \frac{\partial^2}{\partial x_n\partial x_1}\! +\!\sum\limits_{\beta}\! \Gamma_{1\beta}^\beta \frac{\partial }{\partial x_n}
\! +\!\sum\limits_{\beta} \!\frac{\partial \Gamma_{1\beta}^\beta}{\partial x_n} & \!\cdots\! &  \frac{\partial^2}{\partial x_n^2} \!+\!\sum\limits_{\beta} \!\Gamma_{n\beta}^\beta \frac{\partial }{\partial x_n}
\! +\!\sum\limits_{\beta} \!\frac{\partial \Gamma_{n\beta}^\beta}{\partial x_n}\!\!\end{bmatrix}}\nonumber
 \\  \!\!\! \!\!\!\!&&\! \!\!\!\!\left. + \mu\begin{small}  \begin{bmatrix} \rho_{11} +R_1^{\,1} & \cdots & \rho_{1n}+R_n^{\,1}\\
   \cdots& \cdots & \cdots \\
  \rho_{n-1,1}+R_1^{\,n-1} & \cdots & \rho_{n-1,n}+ R_n^{\,n-1}\\
  \rho_{n1}+R_1^{\,n} & \cdots & \rho_{nn}+ R_n^{\,n}\end{bmatrix}\end{small}\right\}
    \begin{bmatrix} u^1\\ \vdots\\ u^{n-1} \\ u^n\end{bmatrix},\quad \quad \nonumber \end{eqnarray}
where $I_n$ denotes $n\times n$ identity matrix, and $\rho_{js}=  \sum_{k,l=1}^n \big(g^{kl} \frac{\partial \Gamma^j_{sl}}{\partial x_k} +\sum\limits_{h=1}^n g^{kl} \Gamma_{hl}^j \Gamma_{sk}^h -\sum\limits_{h=1}^n g^{kl} \Gamma_{sh}^j \Gamma_{kl}^h \big)$. Here, for the sake of simplicity, we have again used the relationships that $\frac{1}{\sqrt{|g|}} \frac{\partial \sqrt{|g|}}{\partial x_n} =\frac{1}{2} \sum\limits_{\alpha, \beta} g^{\alpha\beta} \frac{\partial g_{\alpha\beta}}{\partial x_n}$ and $\frac{1}{\sqrt{|g|}} \frac{\partial \sqrt{|g|}}{\partial x_\alpha} =\sum\limits_{\gamma} \Gamma_{\alpha \gamma}^\gamma$.

 Since $\mu>0$ and $\lambda+2\mu>0$, we see that $u$ is a solution of the Lam\'{e} equations ${\mathcal{L}}_g u=0$ if and only if $u$ is a solution of \[  \left[\begin{BMAT}(@, 3pt, 3pt){c.c}{c.c}\bigg[\frac{1}{\mu}\delta_{jk}\bigg]_{(n-1)\times (n-1)} & [0]_{(n-1)\times 1} \\
[0]_{1\times (n-1)} & \frac{1}{\lambda+ 2\mu}  \end{BMAT}\right] {\mathcal{L}}_g u=0.\]
 We denote $\theta_\beta= \sum_{\alpha} \big(\sum_{\gamma}g^{\alpha \beta} \Gamma_{\alpha \gamma}^\gamma+\frac{\partial g^{\alpha\beta}}{\partial x_\alpha}\big)$.
By (\ref{9.4-1}) we have \begin{eqnarray} \label{19.3.20-21} \end{eqnarray}
\begin{eqnarray*} \!\!\!\!\!\!\!\!\;\;\,\;\;\,\;\;\left[\begin{BMAT}(@, 0pt, 0pt){c.c}{c.c}\!\bigg[\!\frac{1}{\mu}\delta_{jk}\bigg]_{(n-1)\times (n-1)} & [0]_{(n-1)\times 1} \\
[0]_{1\times (n-1)} & \frac{1}{\lambda+ 2\mu}  \end{BMAT}\right] {\mathcal{L}}_g u= \left\{\left(\frac{\partial^2}{\partial x_n^2} I_n\right)\right.\qquad \qquad \qquad \qquad \qquad\quad \qquad\quad\end{eqnarray*}
\begin{eqnarray*}
  +  {\begin{small}\left[\!
\begin{BMAT}(@, 0pt, 0pt){c.c}{c.c}\!\begin{small}\bigg[ \!\sum_{\beta}\! \big(\!\Gamma_{n\beta}^\beta\!\frac{\partial}{\partial x_n} \!+\!\sum_{\alpha}\! g^{\alpha \beta} \!\frac{\partial^2}{\partial x_\alpha\partial x_\beta}\! +\! \theta_\beta\! \frac{\partial}{\partial x_\beta}\big)\delta_{jk}\!\bigg]_{(n\!-\!1)\times (n\!-\!1)}\end{small} & \left[0\right]_{(n\!-\!1)\times 1} \\
\big[0\big]_{1\times (n\!-\!1)}& \frac{\mu}{\lambda\!+\!2\mu}\!\begin{small}\sum_{\beta} \! \big(\Gamma_{n\beta}^\beta\!\frac{\partial}{\partial x_n} \!+\! \sum_{\alpha}\!g^{\alpha \beta}\frac{\partial^2}{\partial x_\alpha\partial x_\beta}\! +\! \theta_\beta \frac{\partial}{\partial x_\beta}\big)\end{small}
   \end{BMAT}\!\right] \end{small}}\quad\nonumber\end{eqnarray*}
  \begin{eqnarray*}  + {\begin{small} \left[
\begin{BMAT}(@, 0pt, 0pt){c.c}{c.c}
   \bigg[2\Gamma_{kn}^j \frac{\partial }{\partial x_n}\bigg]_{(n\!-\!1)\times(n\!-\!1)}  & \big[0\big]_{(n\!-\!1)\times 1} \\
   \big[0\big]_{1\times (n\!-\!1)} &0 \end{BMAT}\right]\end{small}} \! +\!
   {\begin{small} \left[
\begin{BMAT}(@, 0pt, 0pt){c.c}{c.c}
   \bigg[2\sum\limits_{\alpha,\beta} g^{\alpha\beta} \Gamma_{k\alpha}^j \frac{\partial }{\partial x_\beta} \bigg]_{(n\!-\!1)\times (n-1)} &
   \bigg[ 2\sum\limits_{\alpha,\beta} g^{\alpha\beta} \Gamma_{n\alpha}^j \frac{\partial }{\partial x_\beta}\bigg]_{(n\!-\!1)\times 1}\\
   \bigg[\frac{2\mu}{\lambda\!+\!2\mu}\sum\limits_{\alpha,\beta} g^{\alpha\beta} \Gamma_{k\alpha}^n \frac{\partial }{\partial x_\beta}\bigg]_{1\times (n\!-\!1)} & 0\end{BMAT} \!\right]\end{small}} \end{eqnarray*}
    \begin{eqnarray*}  + {\begin{small}  \left[\!
\begin{BMAT}(@, 0pt, 0pt){c.c}{c.c}
  \!\bigg[\!\frac{\lambda\!+\!\mu}{\mu} \!\sum_{\alpha}\! g^{j\alpha} \!\big( \!\frac{\partial^2}{\partial x_\alpha\!\partial x_k }\! +\!
   \sum_{\beta} \!\big(\Gamma_{k\beta}^\beta \frac{\partial }{\partial x_\alpha}\! +\!\frac{\partial \Gamma_{k\beta}^\beta}{\partial x_\alpha} \big)\big) \!\bigg]_{(n\!-\!1)\times (n\!-\!1)} & \bigg[\!\frac{\lambda\!+\!\mu}{\mu} \!\sum_{\alpha}g^{j\alpha} \big( \!\frac{\partial^2}{\partial x_\alpha\!\partial x_n } \!+   \! \sum_{\beta}\! \big(\Gamma_{n\beta}^\beta \!\frac{\partial }{\partial x_\alpha} \!+\! \frac{\partial \Gamma_{n\beta}^\beta}{\partial x_\alpha}\! \big)\!\big)\!\bigg] \\
  \bigg[\! \frac{\lambda\!+\!\mu}{\lambda\!+\!2\mu} \big(\frac{\partial^2}{\partial x_n\!\partial x_k}\!+\!\sum_{\beta}\!\big(\Gamma_{k\beta}^\beta \frac{\partial}{\partial x_n} \!+\! \frac{\partial \Gamma_{k\beta}^\beta}{\partial x_n}\!\big)\!\big)\!\bigg]_{1\times (n\!-\!1)} &  \frac{\lambda\!+\!\mu}{\lambda\!+\!2\mu} \!\sum_\beta \!\big(\Gamma_{n\beta}^\beta\!\frac{\partial}{\partial x_n} \!+\! \frac{\partial  \Gamma_{n\beta}^\beta}{\partial x_n}\!\big) \!\end{BMAT} \!\right]\end{small}} \end{eqnarray*}
 \begin{eqnarray*} + { \begin{small} \left.
\left[\begin{BMAT}(@, 0pt, 0pt){c.c}{c.c}
   \!\left[\rho_{jk}+ R_k^{\,j}\right]_{(n-1)\times (n-1)} & \left[\!\rho_{jn}+ R_n^{\,j} \right]_{(n-1)\times 1} \\
 \!\left[\!\frac{\mu}{\lambda+2\mu}(\rho_{nk}+ R_k^{\,n})\!\right]_{1\times (n-1)} & \frac{\mu}{\lambda+2\mu} (\rho_{nn}+ R_n^{\,n})  \end{BMAT}
\right]\!\right\}\!
\end{small} \begin{bmatrix}u^1\\
  \vdots\\
 u^{n\!-\!1}\!\\
 u^n\end{bmatrix}}\!   \quad \qquad \qquad \qquad \quad \qquad \qquad   \\
  =\!\left\{\left(\!\frac{\partial^2}{\partial x_n^2}\, I_n\right) \!+\!B \left(\frac{\partial}{\partial x_n}I_n\!\right)\!+\!C\!\right\}\!
  \begin{bmatrix}u^1\\
  \vdots\\
 u^{n-1}\\
 u^n\end{bmatrix}\!=0,\qquad \qquad \quad\qquad \qquad \qquad \qquad \qquad\qquad\end{eqnarray*}
 where \begin{eqnarray} \label{19.3.19-5,} \end{eqnarray}
\begin{eqnarray*} \;\;\,\;\;\; \begin{small} B:=\left[
\begin{BMAT}(@, 16pt, 16pt){c.c}{c.c} \bigg[\frac{1}{2} \sum_{\alpha, \beta} g^{\alpha\beta} \frac{\partial g_{\alpha\beta}}{\partial x_n} \delta_{jk}+ 2\Gamma_{kn}^j\bigg]_{(n-1)\times (n-1)} & \bigg[\frac{\lambda+\mu}{\mu} \sum_{\alpha} g^{j\alpha}\frac{\partial}{\partial x_\alpha} \bigg]_{(n-1)\times 1}\\
\bigg[\frac{\lambda+\mu}{\lambda+2\mu} \big(\frac{\partial}{\partial x_k} +\sum_{\beta}\Gamma_{k\beta}^\beta\big)\bigg]_{1\times (n-1)} &
\frac{1}{2} \sum_{\alpha,\beta} g^{\alpha\beta} \frac{\partial g_{\alpha\beta}}{\partial x_n}  \end{BMAT} \right], \end{small}\qquad\quad\quad \quad \end{eqnarray*}
\begin{eqnarray}\label{19.3.19-6}\end{eqnarray}
\begin{eqnarray*}
C\!:= \!\!{\begin{small}\left[\!\!
\begin{BMAT}(@, 0pt, 0pt){c.c}{c.c} \begin{small}\bigg[\!\!\sum_{\alpha\!,\beta}\! \big(g^{\alpha \!\beta}\! \frac{\partial^2}{\partial x_\alpha\!\partial x_\beta}\! +\! \big( \!\sum_{\gamma}\! g^{\alpha\! \beta}\!\Gamma_{\alpha \!\gamma}^{\gamma} \!+\!\frac{\partial g_{\alpha\beta}}{\partial x_\alpha}\!\big) \!\frac{\partial}{\partial x_\beta}\big)\!\delta_{jk}\! \bigg]_{\!(\!n\!-\!1\!)\!\times \!(\!n\!-\!1\!)}\end{small}&  \big[0\big]_{(\!n\!-\!1\!)\times 1} \\
\big[0\big]_{1\times (\!n\!-\!1\!)} &\frac{\mu}{\lambda\!+\!2\mu}\!\sum_{\alpha,\beta}\!\big(\! g^{\alpha \beta}\! \frac{\partial^2}{\partial x_\alpha\!\partial x_\beta} \!+\!\! \big(\! \sum_{\gamma} \!g^{\alpha \beta}\!\Gamma_{\alpha \gamma}^{\gamma} \!\!+\!\!\frac{\partial g^{\alpha\beta}}{\partial x_\alpha}\!\big)\! \frac{\partial}{\partial x_\beta}\!\big)\end{BMAT}\!\! \right]\end{small}} \nonumber\end{eqnarray*}
 \[\;\;\;\;\;\;\;\; +{\begin{small} \left[\!
\begin{BMAT}(@, 0pt, 0pt){c.c}{c.c}\bigg[\!2\sum\limits_{\alpha,\beta}\!g^{\alpha\beta}\! \Gamma_{k\alpha}^j\! \frac{\partial }{\partial x_\beta}\bigg]_{(n\!-\!1)\times (n\!-\!1)} & \bigg[ 2\sum\limits_{\alpha,\beta} g^{\alpha\beta} \Gamma_{n\alpha}^j \frac{\partial}{\partial x_\beta}\bigg]_{(n\!-\!1)\times 1}\\
 \bigg[ \frac{2\mu}{\lambda\!+\!2\mu} \sum\limits_{\alpha,\beta} \!g^{\alpha\beta} \Gamma_{k\alpha}^n \frac{\partial }{\partial x_\beta}\bigg]_{1\times(n\!-\!1)} &0\end{BMAT}\right]\end{small}}\qquad \qquad \qquad \qquad \qquad\qquad \qquad \qquad \qquad \qquad \!\]
  \[\, +\!{\begin{small} \left[\!
\begin{BMAT}(@, 0pt, 0pt){c.c}{c.c}\!\bigg[\! \!\frac{\lambda\!+\!\mu}{\mu} \! \sum_{\alpha} \!g^{j\alpha} \! \big(\!\frac{\partial^2}{\partial x_\alpha\!\partial x_k} \!+\!\sum_{\beta} \!\big(\!\Gamma_{k\beta}^\beta \!\frac{\partial}{\partial x_\alpha} \! +\!\frac{\partial \Gamma_{k\beta}^\beta}{\partial x_\alpha}\!\big)\!\big) \!+\! \rho_{jk}\!+\!R_k^{\,j}\!\bigg]&\bigg[\!\frac{\lambda\!+\!\mu}{\mu}\! \sum_{\alpha,\beta} \!g^{j\alpha}\!  \big(\!\Gamma_{\beta n}^\beta\! \frac{\partial}{\partial x_\alpha} \! +\! \frac{\partial \Gamma_{n\beta}^\beta}{\partial x_\alpha}\big) \!+\! \rho_{jn}\!+\!R_n^{\,j}\bigg]_{\!(\!n\!-\!1\!)\times 1\!}\!\\
\bigg[\! \frac{\lambda\!+\!\mu}{\lambda\!+\!2\mu} \!\sum_\beta\! \frac{\partial \Gamma_{k\beta}^\beta}{\partial x_n} \!+\!\frac{\mu}{\lambda\!+\!2\mu} \!\big(\rho_{nk}\!+\!R_k^{\,n}\big)\bigg]_{\!1\times (\!n\!-\!1\!)\!} &{} \frac{\lambda\!+\!\mu}{\lambda\!+\!2\mu}\! \sum_\beta \!\frac{\partial \Gamma_{n\beta}^\beta}{\partial x_n} \!+\!\frac{\mu}{\lambda\!+\!2\mu} \big(\rho_{nn}\!+\!R_n^{\,n}\big)\end{BMAT} \!\right]   \end{small} }\!\]

\vskip 0.2 true cm

\noindent Throughout this paper, we denote $\sqrt{-1}=i$.
\vskip 0.28 true cm

 \noindent{\bf Proposition 3.1.} \ {\it There exists a pseudodifferential operator $Q(x, D_{x'})$ of order one in $x'$ depending smoothly on $x_n$ such that \begin{eqnarray} \label{19.3.19-1} \end{eqnarray}
\begin{eqnarray} \left[\begin{BMAT}(@, 3pt, 3pt){c.c}{c.c}\bigg[\frac{1}{\mu}\delta_{jk}\bigg]_{(n-1)\times (n-1)} & [0]_{(n-1)\times 1} \\
[0]_{1\times (n-1)} & \frac{1}{\lambda+ 2\mu}  \end{BMAT}\right] {\mathcal{L}}_g = \left( \frac{\partial}{\partial x_n} \,I_n + B -Q \right) \left(\frac{\partial }{\partial x_n}\,I_n  +Q\right), \nonumber\end{eqnarray}
 modulo a smoothing operator, where $D_{x'}=(D_{x_1}, \cdots, D_{x_{n-1}})$, $\,D_{x_j}=\frac{1}{i}\, \frac{\partial }{\partial x_j}$.}

 \vskip 0.20 true cm

 \noindent  {\it Proof.} \    We will divide this proof into several steps.

 Step 1.  \  It follows from (\ref{19.3.20-21}), (\ref{19.3.19-5,}) and (\ref{19.3.19-6}) that \begin{eqnarray*} \left[\begin{BMAT}(@, 3pt, 3pt){c.c}{c.c}\bigg[\frac{1}{\mu}\delta_{jk}\bigg]_{(n-1)\times (n-1)}  & [0]_{(n-1)\times 1} \\
[0]_{1\times (n-1)} & \frac{1}{\lambda+ 2\mu}  \end{BMAT}\right] {\mathcal{L}}_g= \left(\frac{\partial^2}{\partial x_n^2}\, I_n\right) +B \left(\frac{\partial}{\partial x_n}I_n\right)+C. \end{eqnarray*}
Let us assume that we have a factorization \begin{eqnarray*} \left[\begin{BMAT}(@, 3pt, 3pt){c.c}{c.c}\bigg[\frac{1}{\mu}\delta_{jk}\bigg]_{(n-1)\times (n-1)} & [0]_{(n-1)\times 1} \\
[0]_{1\times (n-1)}& \frac{1}{\lambda+ 2\mu}  \end{BMAT}\right] {\mathcal{L}}_g =\left(\frac{\partial}{\partial x_n} \,I_n +B - Q\right)\left( \frac{\partial}{\partial x_n}\,I_n+Q\right).\end{eqnarray*}
  Then we have
   \begin{eqnarray} \label{19.3.19-3}  0\!&=\!& - \left[\begin{BMAT}(@, 3pt, 3pt){c.c}{c.c}\bigg[\frac{1}{\mu}\delta_{jk}\bigg]_{(n-1)\times (n-1)} & [0]_{(n-1)\times 1} \\
[0]_{1\times (n-1)}& \frac{1}{\lambda+ 2\mu}  \end{BMAT}\right] {\mathcal{L}}_g\\
\!&& \!+ \left(\frac{\partial}{\partial x_n} I_n +B - Q\right)\left( \frac{\partial}{\partial x_n}I_n+Q\right)\nonumber\\
   \!  &=\!& -C +\left( \big(\frac{\partial}{\partial x_n} I_n\big)Q - Q\big(\frac{\partial}{\partial x_n}I_n\big) \right) +BQ-Q^2.\nonumber \end{eqnarray}
   Let $q(x', \xi')$, $b(x',\xi')$ and $c(x',\xi')$ be the full symbols of $Q$ and $B$ and $C$, respectively. Clearly,  $q(x', \xi') \sim \sum_{j\ge 0} q_{1-j} (x', \xi')$,   $\;b(x',\xi')=b_1(x',\xi')+b_0(x',\xi')$ and $c(x', \xi')= c_2(x',\xi') +c_1(x', \xi') +c_0(x', \xi')$, where \begin{eqnarray}
  \label{19.3.19-4'} \qquad \;\; b_1(x', \xi') \!=\!\left[\begin{BMAT}(@, 0pt, 1pt){c.c}{c.c}\left[0\right]_{(n-1)\times (n-1)} & \left[\frac{i(\lambda\!+\!\mu)} {\mu} \sum_{\alpha} g^{j\alpha} \xi_\alpha \right]_{(n-1)\times 1} \\
\left[ \frac{i(\lambda\!+\!\mu)}{\lambda+2\mu} \xi_k\right]_{1\times (n-1)} & 0 \end{BMAT}\right],\\
 \label{19.3.19-5}\;\;\;\;\;  b_0(x'\!, \xi') \!=\!\left[\!\begin{BMAT}(@, 0pt, 0pt){c.c}{c.c}\left[\!\frac{1}{2}\! \sum_{\alpha,\beta}\!g^{\alpha \beta} \! \frac{\partial g_{\alpha\beta}}{\partial x_n}\delta_{jk}\! +\! 2\Gamma_{kn}^j\!\right]_{(n\!-\!1)\times (n\!-\!1)} & [0]_{(n\!-\!1)\times 1}
 \\  \left[ \!\frac{\lambda\!+\!\mu}{\lambda\!+\!2 \mu} \sum_{\beta}\! \Gamma_{k\beta}^\beta\!\right]_{1\times (n\!-\!1)} & \frac{1}{2}\! \sum_{\alpha,\beta} \!g^{\alpha \beta}\! \frac{\partial g_{\alpha\beta}}{\partial x_n} \end{BMAT}\!\right];\end{eqnarray}
 \begin{eqnarray}
  \label{19.3.19-8} \!\!\!\!\!\!\!\!\!\!\!\!\!\!&&\quad\,\;\;\;\; \;\;{\begin{small} c_2(x'\!, \xi') \!=\! \left[\!\begin{BMAT}(@, 0pt, 0pt){c.c}{c.c}\!\left[\!-\!\sum_{\alpha,\beta}\!\! g^{\alpha\beta} \! \xi_\alpha\xi_\beta \delta_{jk}\!-\!\frac{\lambda\!+\!\mu}{\mu}\!
  \sum_{\alpha}\! \!g^{j\alpha} \xi_\alpha\xi_k \!\right]_{(\!n\!-\!1\!)\times (\!n\!-\!1\!)} & \left[0\right]_{(n\!-\!1)\times 1} \\
\left[ 0\right]_{1\times (n\!-\!1)}  & -\frac{\mu}{\lambda\!+\!2\mu}  \!\sum_{\alpha,\beta} \! g^{\alpha\beta}\! \xi_\alpha \xi_\beta \end{BMAT}\!\right]\end{small}},\end{eqnarray}
\begin{align}& \label{19.3.19-9} \\
& {\begin{small} c_1(x'\!, \!\xi')\! =\!\!\left[\!\begin{BMAT}(@, 0pt, 0pt){c.c}{c.c}\begin{small}\!\left[\! i\!\sum_{\beta} \! \theta_\beta \! \xi_\beta  \delta_{kj} \!+\! i \!\sum_{\alpha,\beta}\!\big(\!\frac{\lambda\!+\!\mu}{\mu} g^{j\alpha} \!\Gamma_{k\beta}^\beta \xi_\alpha \!+\! 2 g^{\alpha\beta} \!\Gamma_{k\alpha}^j \!\xi_\beta\!\big)\!\right]\end{small} & \begin{small}\!\left[\! i\!\sum_{\alpha,\beta}\!\big(\!\frac{ \lambda\!+\!\mu}{\mu} g^{j\alpha} \Gamma_{n\beta}^\beta \xi_\alpha\!+\!2  g^{\!\alpha\!\beta} \Gamma^j_{n\alpha} \xi_\beta\big)\!\right]_{\!(n\!-\!1)\times 1\!}\!\end{small}
 \\  \left[\! \frac{2i \mu}{\lambda\!+\!2\mu} \sum\limits_{\alpha,\beta}\! g^{\alpha\beta}\Gamma_{k\alpha}^n \xi_\beta\!\right]_{1\times (n\!-\!1)}  &  \frac{i\mu}{\lambda\!+\!2\mu} \sum_{\beta}\! \theta_\beta\xi_\beta \end{BMAT}\!\right]\end{small}}\!,\nonumber\end{align}
 \begin{eqnarray} && \label{19.3.19-9.}\end{eqnarray} \begin{eqnarray*}  c_0(x', \xi') \!=\!{\begin{small}\left[\!\begin{BMAT}(@, 0pt, 0pt){c.c}{c.c}\!\left[\! \frac{\lambda\!+\!\mu}{\mu}\! \sum_{\alpha}\! g^{j\alpha}  \frac{\partial \Gamma_{k\beta}^\beta }{\partial x_\alpha} \! +\!\rho_{jk}\!+\!R_k^{\,j}\!\right]_{\!(n\!-\!1)\times (n\!-\!1)\!} & \left[\!\frac{\lambda\!+\!\mu}{\mu} \! \sum_{\alpha, \beta}\! g^{j\alpha}  \frac{\partial \Gamma_{n\beta}^\beta}{\partial x_\alpha}  \!+\!\rho_{jn}\!+\!R_n^{\, j}\! \right]_{\!(n\!-\!1)\times 1\!} \\
  \!\left[\!\frac{\lambda\!+\! \mu}{\lambda\!+\!2\mu} \!\sum_{\beta} \!\frac{\partial \Gamma_{k\beta}^\beta}{ \partial x_n} \! +\!\frac{\mu}{\lambda\!+\!2\mu}\big(\rho_{kn}\!+\! R_k^{\,n}\big)\! \right]_{\!1\times (n\!-\!1)\!} &
    \frac{\lambda \!+\!\mu}{\lambda\!+\!2\mu}\! \sum_\beta \!\frac{\partial \Gamma_{n\beta}^\beta} {\partial x_n} \!+\!\frac{\mu}{\lambda\!+\!2\mu}\big(\rho_{nn}\!+\! R_n^{\, n}\big)  \end{BMAT}\!\right]\end{small}}.\end{eqnarray*}
 Note that for any smooth vector-valued function $v$,  \begin{eqnarray*} \left( \!\big(\frac{\partial}{\partial x_n} I_n\big)q \!-\! q\big(\frac{\partial}{\partial x_n}I_n\big)\! \right)v\!\!\!&=\!\!\!&\!
 \big(\frac{\partial}{\partial x_n} I_n\big)\big(qv)\! -\! q\big(\frac{\partial}{\partial x_n}I_n\big) v\\
 \!\!\!&\!=\!\!\!&\!\!\bigg(\frac{\partial q}{\partial x_n}\bigg) v \!+\! q \big(\frac{\partial }{\partial x_n} I_n\big) v \!-\! q \big(\frac{\partial }{\partial x_n} I_n\big)
 v\!=\! \bigg(\frac{\partial q}{\partial x_n}\!\bigg) v.\nonumber\end{eqnarray*}
 This implies that  \begin{eqnarray} \label{19.3.19-8} \big(\frac{\partial}{\partial x_n} I_n\big)q - q\big(\frac{\partial}{\partial x_n}I_n\big)
= \frac{\partial q}{\partial x_n},\end{eqnarray}
i.e., the symbol of $ (\frac{\partial }{\partial x_n} I_n)Q - Q(\frac{\partial }{\partial x_n} I_n)$ is $\frac{\partial q}{\partial x_n}$.
     Combining this, the right-hand side of (\ref{19.3.19-3}) and symbol formula for product of two pseudodifferential operators (see p.$\,37$ of \cite{Tre}, p.$\,$13 of \cite{Ta2} or \cite{KN})\, we get the full symbol equation:
   \begin{eqnarray} \label{19.3.19-4}\end{eqnarray} \begin{eqnarray} \sum_{\vartheta} \frac{(-i)^{|\vartheta|}}{\vartheta!} \big(\partial^{\vartheta}_{\xi'} q\big)\big(\partial^\vartheta_{x'}q\big)
 - \sum_{\vartheta} \frac{(-i)^{|\vartheta|}}{\vartheta!} \big( \partial^\vartheta_{\xi'}b\big) \big(\partial^\vartheta_{x'} q\big) -\frac{\partial q}{\partial x_n} +c=0.\nonumber\end{eqnarray}

\vskip 0.2 true cm

Step 2. \   Group the homogeneous terms of degree two in (\ref{19.3.19-4}) we obtain the matrix equation $q_1^2 -b_1q_1 +c_2=0$, i.e.,
\begin{eqnarray}\!\!\!\!\!\!\!\!\! \label{19.3.19-16} && q_1^2- \left[\begin{BMAT}(@, 3pt, 3pt){c.c}{c.c} [0]_{(n-1)\times (n-1)}  & \left[ \frac{i(\lambda+\mu)}{\mu} \sum_{\alpha} g^{j\alpha}\xi_\alpha\right]_{(n-1)\times 1} \\
 \left[ \frac{i(\lambda+\mu)}{\lambda+2\mu} \xi_k \right]_{1\times (n-1)} & 0 \end{BMAT} \right] q_1\\
\!\!\!\!\!\!\!\!&&\quad  -\!\left[\!\begin{BMAT}(@, 1pt, 1pt){c.c}{c.c}\left[ \!\sum_{\alpha, \beta}\! g^{\alpha\beta}\! \xi_\alpha \xi_\beta  \delta_{jk} \!+\!\frac{\lambda\!+\!\mu}{\mu} \!\sum_{\alpha}\! g^{j\alpha} \!\xi_\alpha \xi_k\! \right]_{(n\!-\!1)\times (n\!-\!1)} &  [0]_{(n-1)\times 1} \\
 [0]_{1\times (n-1)} & \frac{\mu}{\lambda \!+\!2\mu} \sum_{\alpha, \beta}\!g^{\alpha\beta} \!\xi_\alpha\xi_\beta \end{BMAT} \right]\!=\!0.\nonumber \end{eqnarray}
Our aim is to calculate the unknown $q_1$ by solving the matrix equation (\ref{19.3.19-16}).
 If $\lambda+\mu=0$ then $q_1= \pm \sqrt{\sum_{\alpha, \beta} g^{\alpha\beta}\xi_\alpha\xi_\beta}\, I_n$ (This essentially reduces to the case for
 Dirichlet-to-Neumann map, see \cite{LU}). Therefore, we will only discuss the case of $\lambda+\mu>0$.
  By observing the coefficient matrices of equation (\ref{19.3.19-16}), we see that    \begin{align*} &    {\begin{small}   \left[\!\begin{BMAT}(@, 0pt, 0pt){c.c}{c.c}  \left[\! \frac{\sum_\alpha \!g^{j\alpha} \xi_\alpha \xi_k}{\sqrt{\sum_{\alpha, \beta}\! g^{\alpha\beta} \xi_\alpha\xi_\beta}} \! \right]_{(n\!-\!1)\times (n\!-\!1)}  &[0]_{(n\!-\!1)\times 1} \\
  \left[0\right]_{1\times (n\!-\!1)} & \sqrt{\sum_{\alpha,\beta}\! g^{\alpha\beta} \xi_\alpha\xi_\beta}  \end{BMAT}\!\right]\end{small}} {\begin{small}\left[\begin{BMAT}(@, 0pt, 0pt){c.c}{c.c} \! \left[ \!\frac{\sum_\alpha \!g^{j\alpha} \xi_\alpha \xi_k}{\sqrt{\sum_{\alpha, \beta}\! g^{\alpha\beta} \xi_\alpha\xi_\beta}}  \right]_{(n\!-\!1)\times (n\!-\!1)} &[0]_{(n\!-\!1)\times 1} \\
  \left[0\right]_{1\times (n\!-\!1)} & \sqrt{\sum_{\alpha,\beta} \!g^{\alpha\beta} \xi_\alpha\xi_\beta}  \end{BMAT}\!\right]\end{small}}\end{align*}
      \begin{align*}& = {\begin{small} \left[\begin{BMAT}(@, 0pt, 0pt){c.c}{c.c} [0]_{(n-1)\times (n-1)} & \left[ \sum_\alpha g^{j\alpha} \xi_\alpha \right]_{(n-1)\times 1} \\
  \left[\xi_k\right]_{1\times (n-1)} & 0  \end{BMAT}\right] \end{small}}  {\begin{small}\left[\begin{BMAT}(@, 0pt, 0pt){c.c}{c.c} [0]_{(n-1)\times (n-1)} & \left[ \sum_\alpha g^{j\alpha} \xi_\alpha \right]_{(n-1)\times 1} \\
  \left[\xi_k\right]_{1\times (n-1)} & 0  \end{BMAT}\right]\end{small}}  \\  {}\\
  &     \,              = {\begin{small}  \left[\begin{BMAT}(@, 0pt, 0pt){c.c}{c.c}  \left[\sum_\alpha g^{j\alpha} \xi_\alpha \xi_k \right]_{(n-1)\times (n-1)}  &[0]_{(n-1)\times 1} \\
  \left[0\right]_{1\times (n-1)} &   \sum_{\alpha,\beta} g^{\alpha\beta} \xi_\alpha\xi_\beta  \end{BMAT}\right]\end{small}},\\
 {}\\ &          {\begin{small}\left[\begin{BMAT}(@, 0pt, 0pt){c.c}{c.c} \left[ \frac{\sum_\alpha g^{j\alpha} \xi_\alpha \xi_k}{\sqrt{\sum_{\alpha, \beta} g^{\alpha\beta} \xi_\alpha\xi_\beta}}  \right]_{(n-1)\times (n-1)} &[0]_{(n-1)\times 1} \\
  \left[0\right]_{1\times (n-1)} & \sqrt{\sum_{\alpha,\beta} g^{\alpha\beta} \xi_\alpha\xi_\beta}  \end{BMAT}\right]\end{small} } {\begin{small} \left[\begin{BMAT}(@, 0pt, 0pt){c.c}{c.c} [0]_{(n-1)\times (n-1)} & \left[ \sum_\alpha g^{j\alpha} \xi_\alpha \right]_{(n-1)\times 1} \\
  \left[\xi_k\right]_{1\times (n-1)} & 0  \end{BMAT}\right]\end{small}}\\
  &                  ={\begin{small}  \left[\begin{BMAT}(@, 0pt, 0pt){c.c}{c.c} [0]_{(n-1)\times (n-1)} &\left[ \sum_\alpha g^{j\alpha} \xi_\alpha \right]_{(n-1)\times 1} \\
  \left[\xi_k\right]_{1\times (n-1)} & 0  \end{BMAT}\right]\end{small}}  {\begin{small}   \left[\begin{BMAT}(@, 0pt, 0pt){c.c}{c.c}  \left[ \frac{\sum_\alpha g^{j\alpha} \xi_\alpha \xi_k}{\sqrt{\sum_{\alpha, \beta} g^{\alpha\beta} \xi_\alpha\xi_\beta}}  \right]_{(n-1)\times (n-1)}&[0]_{(n-1)\times 1} \\
  \left[0\right]_{1\times (n-1)} &\sqrt{\sum_{\alpha,\beta} g^{\alpha\beta} \xi_\alpha\xi_\beta}  \end{BMAT}\right]\end{small}}\\
  &         = {\begin{small} \left[\begin{BMAT}(@, 0pt, 0pt){c.c}{c.c} [0]_{(n-1)\times (n-1)} & \left[ \sqrt{\sum_{\alpha,\beta} g^{\alpha\beta} \xi_\alpha\xi_\beta}   \sum_\alpha g^{j\alpha} \xi_\alpha \right]_{(n-1)\times 1}\\
  \left[ \sqrt{\sum_{\alpha,\beta} g^{\alpha\beta} \xi_\alpha\xi_\beta}\,\, \xi_k\right]_{1\times (n-1)} & 0  \end{BMAT}\right]\end{small}}.\end{align*}
   Let $\mathcal{S}$ denote all the diagonal matrices of form:
 \begin{eqnarray*} && \begin{small} \left[\begin{BMAT}(@, 1pt, 1pt){c.c}{c.c} \left[d_1\delta_{jk}\right]_{(n-1)\times (n-1)} & \left[0\right]_{(n-1)\times 1} \\
 \left[ 0\right]_{1\times (n-1)} & d_2\end{BMAT} \right]\end{small}, \quad  d_1,\,d_2\in \mathbb{C}. \end{eqnarray*}
 Obviously, $\mathcal{S}$ is a commutative ring according to the addition and multiplication of matrices.
  Thus, the following three matrices play a key role:
\begin{eqnarray*} & \sqrt{\sum_{\alpha, \beta} g^{\alpha\beta} \xi_\alpha\xi_\beta}\,I_n, \;\,\,\;
    {\begin{small} \left[\begin{BMAT}(@, 0pt, 0pt){c.c}{c.c} \begin{small} \left[ \frac{\sum_\alpha g^{j\alpha} \xi_\alpha \xi_k}{\sqrt{\sum_{\alpha, \beta} g^{\alpha\beta} \xi_\alpha\xi_\beta}}  \right]_{(n-1)\times (n-1)} \end{small} &[0]_{(n-1)\times 1} \\
  \left[0\right]_{1\times (n-1)} &\begin{small} \sqrt{\sum_{\alpha,\beta} g^{\alpha\beta} \xi_\alpha\xi_\beta}\end{small}  \end{BMAT}\right]\end{small}}, \; \\
  &  \left[\begin{BMAT}(@, 0pt, 0pt){c.c}{c.c} [0]_{(n-1)\times (n-1)} & \begin{small}\left[ \sum_\alpha g^{j\alpha} \xi_\alpha \right]_{(n-1)\times 1}\end{small} \\
  \left[\xi_k\right]_{1\times (n-1)} & 0  \end{BMAT}\right].
\end{eqnarray*}
 The set $F$ of above three matrices can generate a matrix ring $\mathfrak{F}$ according to the following two operations: we first define a left multiplication operation between $\mathcal{S}$ and $\mathfrak{F}$: for every $V\in \mathcal{S}$ and $A\in \mathfrak{F}$, we have $VA\in \mathfrak{F}$ by the usual multiplication of matrix (The element of $\mathcal{S}$ is said to be the ``coefficient matrix'' of matrix ring $\mathfrak{F}$); we then define the addition and multiplication by using the usual matrix addition and multiplication of $\mathfrak{F}$. Clearly, $\mathfrak{F}$ is a three-dimensional matrix ring on ``coefficient ring'' $\mathcal{S}$, and  $F$ is a basis of $\mathfrak{F}$.
  This implies that the solution $q_1$ of  equation (\ref{19.3.19-16}) must have the following form:
   \begin{align} \label{19.3.21-1} &  q_1 =   s_1 \sqrt{\sum_{\alpha, \beta} g^{\alpha\beta} \xi_\alpha\xi_\beta} \; I_n\\
    & +  \begin{small} \left[\begin{BMAT}(@, 0pt, 0pt){c.c}{c.c}\left[ s_2\delta_{jk} \right]_{(n-1)\times (n-1)} & \left[0\right]_{(n-1)\times 1}\\
    \left[0\right]_{1\times (n-1)} & s_3 \end{BMAT} \right] \end{small} {\tiny\left[\begin{BMAT}(@, 0pt, 0pt){c.c}{c.c}  \begin{small}\left[ \frac{\sum_\alpha g^{j\alpha} \xi_\alpha \xi_k}{\sqrt{\sum_{\alpha, \beta} g^{\alpha\beta} \xi_\alpha\xi_\beta}}  \right]_{(n-1)\times(n-1)}\end{small}  &[0]_{(n-1)\times 1}  \\
  \left[0\right]_{1\times (n-1)} & \sqrt{\sum_{\alpha,\beta} g^{\alpha\beta} \xi_\alpha\xi_\beta}  \end{BMAT}\right]}\nonumber\\
  & +\begin{small} \left[\begin{BMAT}(@, 1pt, 1pt){c.c}{c.c} [s_4\delta_{jk}]_{(n-1)\times (n-1)} & \left[ 0 \right]_{(n-1)\times 1} \\
  \left[0\right]_{1\times (n-1)} & s_5  \end{BMAT}\right] \end{small}
  \begin{small} \left[\begin{BMAT}(@, 1pt, 1pt){c.c}{c.c} [0]_{(n-1)\times (n-1)} & \left[ \sum_\alpha g^{j\alpha} \xi_\alpha \right]_{(n-1)\times 1} \\
  \left[\xi_k\right]_{1\times (n-1)} & 0  \end{BMAT}\right], \end{small}\nonumber
\end{align}
where $s_1, \cdots, s_5$ are unknown constants which will be determined late.
  The above idea is similar to Galois group theory for solving the polynomial equation (see \cite{Art} or \cite{HME}).
 Substituting  (\ref{19.3.21-1}) into matrix equation (\ref{19.3.19-16}), we get
 \begin{eqnarray*} \label{19.4.5-1} && \!\!\!\!\!\! \left\{ \! s_1 \sqrt{\sum_{\alpha, \beta} g^{\alpha\beta} \xi_\alpha\xi_\beta} \; I_n +  \begin{small} \left[\begin{BMAT}(@, 1pt, 1pt){c.c}{c.c}  \begin{small}\left[ s_2\frac{\sum_\alpha g^{j\alpha} \xi_\alpha \xi_k}{\sqrt{\sum_{\alpha, \beta} g^{\alpha\beta} \xi_\alpha\xi_\beta}}  \right]_{(n-1)\times(n-1)} \end{small} &[0]_{(n-1)\times 1}  \\
  \left[0\right]_{1\times (n-1)} & s_3\sqrt{\sum_{\alpha,\beta} g^{\alpha\beta} \xi_\alpha\xi_\beta}  \end{BMAT}\right]\end{small} \right.\\
   &&\!\!\!\!\! \left.+\begin{small} \left[\begin{BMAT}(@, 1pt, 1pt){c.c}{c.c} [0]_{(n-1)\times (n-1)} & \left[ s_4\sum_\alpha g^{j\alpha} \xi_\alpha \right]_{(n-1)\times 1} \\
  \left[s_5\, \xi_k\right]_{1\times (n-1)} & 0  \end{BMAT}\right] \end{small}\right\}\left\{  s_1 \sqrt{\sum_{\alpha, \beta} g^{\alpha\beta} \xi_\alpha\xi_\beta} \; I_n \right.\\
[1.4mm]  && \!\!\! \!\!\left.+  \begin{small} \left[\begin{BMAT}(@, 1pt, 1pt){c.c}{c.c}  \begin{small}\left[ s_2\frac{\sum_\alpha g^{j\alpha} \xi_\alpha \xi_k}{\sqrt{\sum_{\alpha, \beta} g^{\alpha\beta} \xi_\alpha\xi_\beta}}  \right]_{(n-1)\times(n-1)} \end{small} &[0]_{(n-1)\times 1}  \\
  \left[0\right]_{1\times (n-1)} & s_3\sqrt{\sum_{\alpha,\beta} g^{\alpha\beta} \xi_\alpha\xi_\beta}  \end{BMAT}\right]\end{small} \right.\\
 [1.3mm] &&\!\!\!\!\! \left.+\begin{small} \left[\begin{BMAT}(@, 1pt, 1pt){c.c}{c.c} [0]_{(n-1)\times (n-1)} & \left[ s_4\sum_\alpha g^{j\alpha} \xi_\alpha \right]_{(n-1)\times 1} \\
  \left[s_5\, \xi_k\right]_{1\times (n-1)} & 0  \end{BMAT}\right] \end{small}\right\}\\
 [1.3mm]&& \!\!\!\!\!  - \begin{small} \left[\begin{BMAT}(@, 1pt, 1pt){c.c}{c.c} \left[0\right]_{(n-1)\times (n-1)}
   & \begin{small}\left[\frac{i(\lambda+\mu)}{\mu} \sum_{\alpha} g^{j\alpha} \xi_\alpha\right]_{(n-1)\times 1}\end{small} \\
    \left[ \frac{i(\lambda+\mu)}{\lambda+2\mu} \xi_k\right]_{1\times (n-1)} & 0\end{BMAT}\right] \end{small} \left\{
     s_1 \sqrt{\sum_{\alpha, \beta} g^{\alpha\beta} \xi_\alpha\xi_\beta} \; I_n \right.\\
  [1.3mm] &&\!\!\!\!\! \left.+  \begin{small} \left[\begin{BMAT}(@, 1pt, 1pt){c.c}{c.c}  \begin{small}\left[ s_2\frac{\sum_\alpha g^{j\alpha} \xi_\alpha \xi_k}{\sqrt{\sum_{\alpha, \beta} g^{\alpha\beta} \xi_\alpha\xi_\beta}}  \right]_{(n-1)\times(n-1)} \end{small} &[0]_{(n-1)\times 1}  \\
  \left[0\right]_{1\times (n-1)} & s_3\sqrt{\sum_{\alpha,\beta} g^{\alpha\beta} \xi_\alpha\xi_\beta}  \end{BMAT}\right]\end{small}  \right.\\
  [1.4mm]&&\!\!\! \!\left.+{ \left[\begin{BMAT}(@, 1pt, 1pt){c.c}{c.c} [0]_{(n-1)\times (n-1)} & \left[ s_4\sum_\alpha g^{j\alpha} \xi_\alpha \right]_{(n-1)\times 1} \\
  \left[s_5\, \xi_k\right]_{1\times (n-1)} & 0  \end{BMAT}\right] }\right\}\\
  [1.4mm] && \!\!\! \! - {\begin{small} \left[\begin{BMAT}(@, 1pt, 1pt){c.c}{c.c}  \left[\sum_{\alpha, \beta} g^{\alpha\beta} \xi_\alpha\xi_\beta \;\delta_{jk} + \frac{\lambda+\mu}{\mu}
   \sum_\alpha g^{j\alpha} \xi_\alpha \xi_k  \right]_{(n-1)\times(n-1)}  &[0]_{(n-1)\times 1}  \\
  \left[0\right]_{1\times (n-1)} &  \frac{ \mu}{\lambda+2\mu} \sum_{\alpha,\beta} g^{\alpha\beta} \xi_\alpha\xi_\beta \end{BMAT}\right]\end{small} }=0, \end{eqnarray*}
 i.e.,
    \begin{eqnarray*} \!\!&&\!\!\!\!\!\!\!\sqrt{\sum_{\alpha,\beta} g^{\alpha\beta} \xi_\alpha \xi_\beta} \left\{(s_1^2-1) \sqrt{\sum_{\alpha,\beta} g^{\alpha\beta} \xi_\alpha \xi_\beta}\; I_n   +\left(  \left[\begin{BMAT}(@, 1pt, 1pt){c.c}{c.c} \left[ 2s_1 s_2 \,\delta_{jk}\right]_{(n-1)\times (n-1)} & [0]_{(n-1)\times 1} \\
      [0]_{1\times (n-1)} &2s_1 s_3 \end{BMAT} \right]  \right.\right.\\ [1mm]
      \!\!\!&&\!\!\! \!\! \left.
      + \left[\begin{BMAT}(@, 1pt, 1pt){c.c}{c.c} \left[ s_2^2 \,\delta_{jk}\right]_{(n-1)\times (n-1)} & [0]_{(n-1)\times 1} \\
      [0]_{1\times (n-1)} & s_3^2 \end{BMAT} \right]
      +  \left[\begin{BMAT}(@, 1pt, 1pt){c.c}{c.c} \left[ s_4 s_5 \,\delta_{jk}\right]_{(n-1)\times (n-1)} & [0]_{(n-1)\times 1} \\
      [0]_{1\times (n-1)} & s_4 s_5 \end{BMAT} \right]\right.\\ [1.4mm]
      &&\!\!\!\!\!
     { \left.   -  \left[\begin{BMAT}(@, 1pt, 1pt){c.c}{c.c} \left[ \frac{i(\lambda+\mu)}{\mu} s_5 \,\delta_{jk}\right]_{(n-1)\times (n-1)} & [0]_{(n-1)\times 1} \\
      [0]_{1\times (n-1)} & \frac{i(\lambda+\mu)}{\lambda+2\mu} \,s_4 \end{BMAT} \right]  -
       \left[\begin{BMAT}(@, 1pt, 1pt){c.c}{c.c} \left[  \frac{\lambda+\mu}{\mu} \delta_{jk}\right]_{(n-1)\times (n-1)} & [0]_{(n-1)\times 1} \\
      [0]_{1\times (n-1)} & \frac{\mu}{\lambda+2\mu}-1 \end{BMAT} \right] \right)}\\ [1mm]
        \!\!\!&&\!\!\! \!\!\quad \times {  \left[\begin{BMAT}(@, 1pt, 1pt){c.c}{c.c}  \left[ \frac{\sum_\alpha g^{j\alpha} \xi_\alpha \xi_k}{\sqrt{\sum_{\alpha, \beta} g^{\alpha\beta} \xi_\alpha\xi_\beta}}  \right]_{(n-1)\times (n-1)}  &[0]_{(n-1)\times 1} \\
  \left[0\right]_{1\times (n-1)} & \sqrt{\sum_{\alpha,\beta} g^{\alpha\beta} \xi_\alpha\xi_\beta}  \end{BMAT}\right] }\\ [1mm]
   \!\!\!&&\!\!\! \!\!\!\;\;\! +\left(  \left[\begin{BMAT}(@, 1pt, 1pt){c.c}{c.c} \left[ 2s_1 s_4 \,\delta_{jk}\right]_{(n-1)\times (n-1)} & [0]_{(n-1)\times 1} \\
      [0]_{1\times (n-1)} &2s_1 s_5 \end{BMAT} \right]  \right.\\ [1.4mm]&& \!\!\!\!\left.
     + \left[\begin{BMAT}(@, 1pt, 1pt){c.c}{c.c} \left[ s_2s_4 \,\delta_{jk}\right]_{(n-1)\times (n-1)} & [0]_{(n-1)\times 1} \\
      [0]_{1\times (n-1)} & s_3 s_5 \end{BMAT} \right]
      +  \left[\begin{BMAT}(@, 1pt, 1pt){c.c}{c.c} \left[ s_3 s_4 \,\delta_{jk}\right]_{(n-1)\times (n-1)} & [0]_{(n-1)\times 1} \\
      [0]_{1\times (n-1)} & s_2 s_5 \end{BMAT} \right]  \right.
      \\ [1.4mm] \!\!\!&&\!\!\! \!\!
     \left. \left.\,  - { \left[\begin{BMAT}(@, 1pt, 1pt){c.c}{c.c} \left[ \frac{i(\lambda\!+\!\mu)}{\mu} s_1 \,\delta_{jk}\right]_{(n\!-\!1)\times (n\!-\!1)} & [0]_{(n\!-\!1)\times 1} \\
      [0]_{1\times (n\!-\!1)} & \frac{i(\lambda\!+\!\mu)}{\lambda\!+\!2\mu} \,s_1 \end{BMAT} \right]} \! -\!
               { \left[\begin{BMAT}(@, 1pt, 1pt){c.c}{c.c} \left[ \frac{i(\lambda\!+\!\mu)}{\mu} s_3 \,\delta_{jk}\right]_{(n\!-\!1)\times (n\!-\!1)} & [0]_{(n\!-\!1)\times 1} \\
      [0]_{1\times (n\!-\!1)} & \frac{i(\lambda\!+\!\mu)}{\lambda\!+\!2\mu} \,s_2 \end{BMAT} \right] } \right)\right.\qquad \quad\;\; \\
       \!\!\!&&\!\!\! \!\! \quad \left. \times  \left[\begin{BMAT}(@, 1pt, 1pt){c.c}{c.c} [0]_{(n-1)\times (n-1)} & \left[ \sum_\alpha g^{j\alpha} \xi_\alpha \right]_{(n-1)\times 1} \\
  \left[\xi_k\right]_{1\times (n-1)} & 0  \end{BMAT}\right] \right\}=0.\\
  && {} \\
  && {}\end{eqnarray*}
   By comparing the ``coefficient matrices'' in front of each element of $F$, we get
  \begin{eqnarray}  \label{19.3.22-1} \left\{ \begin{array}{ll} s_1^2-1=0, \\
2s_1 s_2 +s_2^2 +s_4s_5 -\frac{i(\lambda+\mu)}{\mu} s_5 -\frac{\lambda+\mu}{\mu} =0,\\
 2s_1s_3 +s_3^2 +s_4 s_5  - \frac{i(\lambda+\mu)}{\lambda+2\mu} s_4 +1-\frac{\mu}{\lambda+2\mu}  =0,\\
2s_1s_4 +s_2s_4 +s_3s_4 -\frac{i(\lambda+\mu)}{\mu} s_1 -\frac{i(\lambda+\mu)}{\mu}s_3 =0,\\
2s_1s_5 +s_3 s_5 +s_2s_5 -\frac{i(\lambda+\mu)}{\lambda+2\mu} s_1  -\frac{i(\lambda+\mu)}{\lambda+2\mu}s_2=0.\end{array}\right.\end{eqnarray}
Assume that $s_1+s_3\ne 0$ and $s_3-s_2\ne 0$. Then the equation system (\ref{19.3.22-1}) can be rewritten as
\begin{eqnarray}\! \! \!\!\label{19.4.5-3}\;\;\;\;\quad\;\;\;\;\;\;{  \left\{ \!\!\!\begin{array}{ll} s_1^2=1, \\ [2mm]
s_4= -i(s_1+s_3) \big( \frac{\lambda+2\mu}{\lambda+\mu} (2s_1+ s_2+s_3) +\frac{\lambda+3\mu}{\mu}(\frac{1}{s_3-s_2})\big),\\ [2mm]
s_5= -i(s_1+s_2) \big( \frac{\mu}{\lambda+\mu} (2s_1+ s_2+s_3) +\frac{\lambda+3\mu}{\lambda+2\mu}(\frac{1}{s_3-s_2})\big),\\
[3mm]
s_3^3\!+\!\left(4\,s_{1}\!+\!s_{2}\right)s_3^2\!+\!\big(4s_1^2 \!-\!s_2^2 \!+\!\frac{2(\lambda\!+\!\mu)}{\mu} \big)s_{3}
\!+\!\frac{(\lambda\!+\!\mu)(2\lambda s_1\!+ \!6\mu s_1\!+\!2\mu s_2)}{\mu(\lambda\!+\!2\mu)} \!-\!s_2(2s_1\!+\!s_2)^2\!=\!0,\\ [2mm]
 \big(\frac{\lambda\! +\!\mu }{\mu }\!+\!{s_{1}}^2\big)\,s_3^2\!+\!\big(2s_1^3 \!+\! \frac{2(\lambda\!+\!\mu)}{\mu}s_{1} \big)s_3 \!-\! s_1^2 s_{2}\,(2\,s_{1}\!+\!s_{2})\!+\!\frac{(\lambda\!+\!\mu)\big( (\lambda\!+\!2\mu)s_1^2 \!+\!\mu (s_1\!+\!s_2)^2\big)}{\mu (\lambda \!+\!2\mu)}
\!=\!0. \end{array}\right.}\end{eqnarray}
Because we have chosen the outward normal $\nu$ of $\partial \Omega$, we should take
\begin{eqnarray}\label{19.4.6-6} s_1=1, \quad s_2\ge 0, \quad s_1+s_3 > 0.\end{eqnarray}
Roughly speaking, such a choice means that the matrix $q_1$ is positive-definition .
The last two equations of (\ref{19.4.5-3}) can be solved by the well-known resultant method. Clearly, the resultant equation of the left-hand sides of last two equations in (\ref{19.4.5-3}) is
\begin{eqnarray*}\!\!\! &\!\!\! {\left| \begin{matrix}
	\!1&\!\!	\!	s_2+4&	\!\!	4-s_2^2 +\frac{2(\lambda\!+\!\mu)}{\mu}\!\!\!&\!	\!\!\!\!	\begin{small}\frac{2(\lambda\!+\!\mu)\big((\lambda+3\mu)\!+\!\mu s_2\big)}{\mu (\lambda+2\mu)}\!-\!s_2(2\!+\!s_2)^2\end{small} \!\!\!&	\!\!\!	0\\
	\!0\!&	\!\!\!	1\!&\!	\!	s_2\!+\!4\!\!&\!\!\!\!	\!	4\!-\!s_2^2 \!+\!\frac{2(\lambda\!+\!\mu)}{\mu}\!\!\!\!&\!	\!\!\!	\frac{2(\lambda\!+\!\mu)\big((\lambda+3\mu)\!+\!\mu s_2\big)}{\mu (\lambda\!+\!2\mu)}\!-\!s_2(2\!+\!s_2)^2\\
\!	\frac{\lambda \!+\!2\mu}{\mu}\!\!&	\!	\!\frac{2(\lambda\!+\!2\mu)}{\mu}\!\!&		\!\!\frac{(\lambda\!+\!\mu)( \mu s_2^2\!+\!2\mu s_2\!+\!\lambda\!+\!3\mu)}{\mu(\lambda\!+\!2\mu)} \!-\!s_2^2 \!-\!2s_2\!\!\!\! &\!	\!\!	\!0\!\!\!\!&\!	\!\!\!0\\
\!	\!0\!&	\!\!	\frac{\lambda \!+\!2\mu}{\mu}\!\!&	\!\!	\frac{2(\lambda\!+\!2\mu)}{\mu}\!\! \!&\!	\!\!	\!\frac{(\lambda\!+\!\mu)( \mu s_2^2 \!+\!2\mu s_2\!+\!\lambda\!+\!3\mu)}{\mu(\lambda\!+\!2\mu)} \!-\!s_2^2 \!\!-\!2s_2\!\!\! \!\!&\!	\!\!\!	0\\
\!	\!0\!&\!	\!0\!&\!	\!	\frac{\lambda \!+\!2\mu}{\mu}\! \!&	\!\!	\!\frac{2(\lambda\!+\!2\mu)}{\mu}\!\! \!\!&	\!\!\!	\!\!\frac{(\lambda\!+\!\mu)( \mu s_2^2 \!+\!2\mu s_2\!+\!\lambda\!+\!3\mu)}{\mu(\lambda\!+\!2\mu)}\! -\!s_2^2 \!-\!2s_2
\end{matrix}\! \right|}\\ {} \\ [3mm]
&\!\!\!\!\!\!\!\!\!\!\!\!\!\!\!\!\!\!\!\!\!\!\!\!\!\!\!\!\!\!\!\!\!\!=\frac{\left( s_2+1 \right) ^4\,\left( \lambda +\text{3\,}\mu \right) \,\left( \lambda +\mu \right) ^3\,\left( \lambda \,s_2-\mu -\lambda +\text{3\,}\mu \,s_2 \right) \,\left( \text{3\,}\lambda +\text{7\,}\mu +\lambda \,s_2+\text{3\,}\mu \,s_2 \right)}{\mu ^3\,\left( \lambda +\text{2\,}\mu \right) ^3}=0, \qquad \qquad \qquad \qquad \qquad \qquad \end{eqnarray*}
which has only three solutions $s_2=-1$, $s_2=-\frac{3\lambda+7\mu}{\lambda+3\mu}$ and $s_2=\frac{\lambda+\mu}{\lambda+3\mu}$.
 It is well-known (see \cite{GKZ}) that $(s_2,s_3)$ is a solution of the last two equations of (\ref{19.4.5-3})
 if and only if $s_2$ is a solution of the above resultant equation.
 Inserting these values of $s_2$ into (\ref{19.4.5-3}) we get that the solutions of  last two equations of (\ref{19.4.5-3}) are
\begin{eqnarray*} \left\{ \begin{array}{ll} s_2=-1 \\
s_3=-1, \end{array} \right.\qquad \;\end{eqnarray*}
\begin{eqnarray*} \mbox{or}\;\; \left\{ \begin{array}{ll} s_2= -\frac{3\lambda+7\mu}{\lambda+3\mu}\\
 s_3= -\frac{\lambda+5\mu}{\lambda+3\mu},\end{array}\right.\qquad \; \end{eqnarray*}
 \begin{eqnarray*} \mbox{or}\;\;\left\{ \begin{array}{ll}  s_2= \frac{\lambda+\mu}{\lambda+3\mu}\\
  s_3= -\frac{\lambda+\mu}{\lambda+3\mu}.\end{array}\right.\qquad \; \end{eqnarray*}
Therefore, according to (\ref{19.4.6-6}) we obtain the solution of last two equations of (\ref{19.4.5-3}):
\begin{eqnarray} \label{19.4.6-7} s_2=  \frac{\lambda+\mu}{\lambda+3\mu}\;\, \mbox{and}\;\, s_3= -\frac{\lambda+\mu}{\lambda+3\mu}.\end{eqnarray}
Substituting them into (\ref{19.4.5-3}) again, we get
\begin{eqnarray} \label{19.4.6-8} \left\{\begin{array}{ll} s_1=1,\quad s_2= \frac{\lambda+\mu}{\lambda+3\mu}, \quad
s_3=- \frac{\lambda+\mu}{\lambda+3\mu},\\
s_4= \frac{i(\lambda+\mu)}{\lambda+3\mu},\quad
s_5 = \frac{i(\lambda+\mu)}{\lambda+3\mu}.\end{array}\right.\end{eqnarray}
Thus,  \begin{eqnarray} \label{19.3.22-8}\end{eqnarray} \begin{eqnarray}   q_1(x',\xi')\! \!&\!\!=\!\!&\! \!{\begin{small} \left[\!\begin{BMAT}(@, 1pt, 1pt){c.c}{c.c} \left[\!s_1 \sqrt{\sum_{\alpha, \beta} \!g^{\alpha\beta} \xi_\alpha\xi_\beta}  \delta_{jk}\! +\! \frac{s_2\sum_\alpha \!g^{j\alpha} \xi_\alpha \xi_k}{\sqrt{\sum_{\alpha, \beta}\! g^{\alpha\beta} \xi_\alpha\xi_\beta}}  \!\right]_{(n\!-\!1)\times(n\!-\!1)} & \left[s_4\sum_\alpha\! g^{j\alpha} \xi_\alpha \right]_{(n\!-\!1)\times 1} \\
\left[s_5\, \xi_k\right]_{1\times (n\!-\!1)}
   &(s_1  \!+\!s_3)\sqrt{\sum_{\alpha,\beta}\! g^{\alpha\beta} \xi_\alpha\xi_\beta}  \end{BMAT}\right],\end{small} }\nonumber
\end{eqnarray}
where $s_1, \cdots, s_5$ are given by (\ref{19.4.6-8}).

\vskip 0.2 true cm

Step 3. \  The terms of degree one in (\ref{19.3.19-4}) are
\begin{eqnarray}  \label{19.3.22-10} &&  q_1 q_0 +q_0q_1  -i\sum_{l=1}^{n-1}  \frac{\partial q_1}{\partial \xi_l} \, \frac{\partial q_1}{\partial x_l}  - b_1 q_0 - b_0 q_1\\
&&  \qquad -\frac{1}{i} \sum_{l=1}^{n-1} \frac{\partial b_1}{\partial \xi_l} \, \frac{\partial q_1}{\partial x_l}   -\frac{\partial q_1}{\partial x_n} +c_1=0, \nonumber\end{eqnarray}
i.e., \begin{eqnarray} \label{19.4.17-1} (q_1- b_1) q_0 +q_0 q_1= E_1,\end{eqnarray}
where \begin{eqnarray} \label{19.3.23-1} E_1:= i\sum_{l=1}^{n-1} \frac{\partial q_1}{\partial \xi_l} \, \frac{\partial q_1}{\partial x_l} + b_0q_1 -i \sum_{l=1}^{n-1} \frac{\partial b_1}{\partial \xi_l}\, \frac{\partial q_1}{\partial x_l} +\frac{\partial q_1}{\partial x_n}  -c_1,\end{eqnarray}
 and $b_0$ and $c_1$ are given in (\ref{19.3.19-5}) and (\ref{19.3.19-9}).
More precisely, \begin{eqnarray}  \label{19.3.22-11}\end{eqnarray} \begin{eqnarray} \!\!\!\!\!\!\!\!&& \!\! {\begin{small}   \left[\begin{BMAT}(@, 1pt, 1pt){c.c}{c.c}\! \left[ s_1 \sqrt{\sum_{\alpha, \beta}\! g^{\alpha\beta} \xi_\alpha \xi_\beta} \, \delta_{jk}
 \!+\!s_2 \frac{\sum_{\alpha} g^{j\alpha} \xi_\alpha \xi_k}{\sqrt{\sum_{\alpha,\beta} g^{\alpha \beta} \xi_\alpha \xi_\beta}} \right]_{(n\!-\!1)\times(n\!-\!1)} &  \left[\big(s_4 \!-\!\frac{i(\lambda\!+\!\mu)}{\mu}\big) \sum_{\alpha} g^{j\alpha}\xi_\alpha\right]_{(n\!-\!1)\times 1} \\
 \left[\big(s_5\!-\! \frac{i(\lambda\!+\!\mu)}{\lambda\!+\!2\mu} \big) \xi_k\!\right] _{1\times(n-1)} & (s_1\!+\!s_3) \sqrt{\sum_{\alpha, \beta}\! g^{\alpha \beta} \xi_\alpha \xi_\beta }\end{BMAT} \right]q_0\end{small}}\nonumber \\
\!\! \!\!\!\!\!\!&&\!\! +
 q_0 {\begin{small} \left[\begin{BMAT}(@, 1pt, 1pt){c.c}{c.c} \left[ s_1 \sqrt{\sum_{\alpha, \beta} g^{\alpha\beta} \xi_\alpha \xi_\beta} \, \delta_{jk}
 +s_2 \frac{\sum_{\alpha} g^{j\alpha} \xi_\alpha \xi_k}{\sqrt{\sum_{\alpha,\beta} g^{\alpha \beta} \xi_\alpha \xi_\beta}} \right]_{(n-1)\times(n-1)} &   \left[s_4\sum_{\alpha} g^{j\alpha}\xi_\alpha\right]_{(n-1)\times 1} \\
 \left[s_5 \xi_k\right]_{1\times(n-1)} & (s_1+s_3) \sqrt{\sum_{\alpha, \beta} g^{\alpha \beta} \xi_\alpha \xi_\beta }\end{BMAT} \right]\end{small}}=E_1.\nonumber \end{eqnarray}
   Now, we calculate $q_0$ by solving Sylvester's matrix equation (\ref{19.3.22-11}).
It is well-known that Sylvester's equation of the form
$LX+XM=E$ can be put into the form (see \cite{BaS} or \cite{BhR})
\begin{eqnarray} \label{19.10.6-1} U(\mbox{vec}\,X)= \mbox{V}\end{eqnarray}  for larger matrices $U$ and $V$.
 Here  $\mbox{vec}\, X$ is a stack of all columns of matrix $X$ (see, for example,  Chapter 4 of \cite{HoJ}).
Indeed, $U=(I_n\otimes L)+(M^t \otimes I_n)$, and $V=\mbox{vec}\,E$, where $\otimes$ denotes the Kronecker product.
Thus, if we can obtain the inverse $U^{-1}$ of the matrix $U$, then we have $\mbox{vec}\,X=U^{-1} (\mbox{vec}\,V)$, and the corresponding solution $X$ will immediately be obtained.
From (\ref{19.3.22-11}) we see that
 \begin{eqnarray*}  \label{19.3.22-12} L\!\!&\!\!:=\!\!&\!\! {\begin{small} \left[\begin{BMAT}(@, 1pt, 1pt){c.c}{c.c} \left[ s_1 \sqrt{\sum_{\alpha, \beta} g^{\alpha\beta} \xi_\alpha \xi_\beta} \; \delta_{jk} \! +\!s_2 \frac{\sum_{\alpha} g^{j\alpha} \xi_\alpha \xi_k}{\sqrt{\sum_{\alpha,\beta}\! g^{\alpha \beta} \xi_\alpha \xi_\beta}} \!\right]_{(n\!-\!1)\times(n\!-\!1)} &  \left[\big(s_4 \!-\!\frac{i(\lambda\!+\!\mu)}{\mu}\big) \sum_{\alpha} \!g^{j\alpha}\xi_\alpha\!\right]_{(n\!-\!1)\times 1} \\
 \left[\big(s_5\!-\! \frac{i(\lambda\!+\!\mu)}{\lambda\!+\!2\mu} \big) \xi_k\right]_{1\times (n\!-\!1)} & (s_1\!+\!s_3) \sqrt{\sum_{\alpha, \beta}\!g^{\alpha \beta} \xi_\alpha \xi_\beta }\end{BMAT} \right]\end{small}}\end{eqnarray*} and
 \begin{eqnarray*} M^t :={\begin{small} \left[\begin{BMAT}(@, 1pt, 1pt){c.c}{c.c} \left[ s_1 \sqrt{\sum_{\alpha, \beta} g^{\alpha\beta} \xi_\alpha \xi_\beta} \; \delta_{jk}
 +s_2 \frac{\sum_{\alpha} g^{k\alpha} \xi_\alpha \xi_j}{\sqrt{\sum_{\alpha,\beta} g^{\alpha \beta} \xi_\alpha \xi_\beta}} \right]_{(n-1)\times (n-1)}  &  \left[ s_5\,\xi_j \right]_{(n-1)\times 1} \\
 \left[s_4 \sum_{\alpha} g^{k\alpha} \xi_\alpha \right]_{1\times (n-1)} & (s_1+s_3) \sqrt{\sum_{\alpha, \beta} g^{\alpha \beta} \xi_\alpha \xi_\beta}\end{BMAT} \right].\end{small}}\nonumber \end{eqnarray*}
Thus the  $U$ has the form:
 \begin{eqnarray} \label{19.3.23-15}&&\!\!\!\! U = (I_n\otimes L) + (M^t \otimes I_n) \\
&&\quad \!\!\! =  \left(2s_1 \sqrt{\sum_{\alpha,\beta}\! g^{\alpha\beta} \xi_\alpha\xi_\beta}\,  I_{n^2}\!\right) \!+\! I_n\otimes\left(\begin{small}\left[\begin{BMAT}(@, 1pt, 1pt){c.c}{c.c} \left[ s_2 \, \delta_{jk}\right]_{(n\!-\!1)\times (n\!-\!1)}& \left[0\right]_{(n\!-\!1)\times 1}\\
\left[0\right]_{1\times (n\!-\!1)}&s_3  \end{BMAT}\right]\end{small}A_1\!\right)\nonumber\\
[1mm]&& \qquad \!\! +I_n\otimes \left(\begin{small}\left[\begin{BMAT}(@, 1pt, 1pt){c.c}{c.c} \left[\big(s_4-\frac{i(\lambda+\mu)}{\mu}\big)\delta_{jk}\right]_{(n\!-\!1)\times (n\!-\!1)}& \left[0\right]_{(n\!-\!1)\times 1}\\
\left[0\right]_{1\times (n\!-\!1)}  & s_5 -\frac{i(\lambda+\mu)}{\lambda\!+\!2\mu}\end{BMAT}\!\right]\end{small}A_2\right)\nonumber \\
[1mm]&& \qquad \! + \left(\begin{small}\left[\begin{BMAT}(@, 1pt, 1pt){c.c}{c.c} \left[s_2\,\delta_{jk}\right]_{(n-1)\times (n-1)}& \left[0\right]_{(n-1)\times 1}\\
\left[0\right]_{1\times (n-1)}  & s_3\end{BMAT}\right]\end{small}A_1^t\right)\otimes I_n \nonumber\\
[1.5mm]&&  \qquad \! +
 \left(\begin{small}\left[\begin{BMAT}(@, 1pt, 1pt){c.c}{c.c} \left[s_5\delta_{jk}\right]_{(n\!-\!1)\times (n\!-\!1)}& \left[0\right]_{(n-1)\times 1}\\
\left[0\right]_{1\times (n\!-\!1)}  & s_4\end{BMAT}\!\right]\!\end{small}\!A_2^t\!\right)\!\otimes \!I_n,\nonumber\end{eqnarray}
where \begin{eqnarray} \label{19.3.24-1} A_1 &=&  \begin{small}\left[\begin{BMAT}(@, 1pt, 1pt){c.c}{c.c} \begin{small}\left[
 \frac{\sum_{\alpha} g^{j\alpha} \xi_\alpha \xi_k}{\sqrt{\sum_{\alpha, \beta} g^{\alpha \beta} \xi_\alpha \xi_\beta}}\right]_{(n-1)\times (n-1)} \end{small} & \left[0\right]_{(n-1)\times 1}\\
   \left[ 0\right]_{1\times (n-1)} & \sqrt{\sum_{\alpha, \beta} g^{\alpha \beta} \xi_\alpha \xi_\beta} \end{BMAT} \right]\end{small} \\
    A_2 &=&\begin{small}\left[\begin{BMAT}(@, 1pt, 1pt){c.c}{c.c} \left[
 0\right]_{(n-1)\times (n-1)}  & \left[ \sum_{\alpha} g^{j\alpha} \xi_\alpha \right]_{(n-1)\times 1}\\
   \left[ \xi_k \right]_{1\times (n-1)} &0 \end{BMAT} \right].\end{small}  \end{eqnarray}
We denote by  $\theta(d_1,d_2; d_3, d_4)$ the $n^2\times n^2$ diagonal matrix of the form
\begin{eqnarray*} \theta(d_1, d_2; d_3, d_4)\! =\! {\tiny \begin{bmatrix}  \begin{bmatrix} d_1&\\
 & \ddots && \\
  && d_1  \\
   && & d_2\end{bmatrix} \\
 & \ddots \\
 && \begin{bmatrix} d_1\\
 & \ddots && \\
  && d_1  \\
   && & d_2\end{bmatrix}\\
  &&&\!\!\!\begin{bmatrix} d_3&\\
 & \ddots && \\
  && d_3  \\
   && & d_4\end{bmatrix}\end{bmatrix}},\end{eqnarray*}
and denote by $\mathfrak{J}$ the all such diagonal matrices $\theta(d_1,d_2;d_3,d_4)$, where $d_1$, $d_2$, $d_3$, $d_4\in \mathbb{C}$. We call each element of the matrix ring $\mathfrak{J}$ as ``coefficient matrix''. Since $(I_n \otimes A_l)(A_m \otimes I_n)=(A_m \otimes I_n) (I_n \otimes A_l)$ for all $1\le l,m\le 2$, we see that the following nine matrices generate a matrix ring $\mathfrak{X}$ about the addition and multiply of matrices as well as the left-multiplication operation between $\mathfrak{J}$ and $\mathfrak{X}$:
\begin{eqnarray} \label{19.4.6-p2}\end{eqnarray}\begin{eqnarray}\!\!\!\!&&\!\!\!\! H=\begin{small} \bigg\{ \sqrt{\sum_{\alpha, \beta} g^{\alpha \beta}\; \xi_\alpha \xi_\beta}\; I_{n^2},\; \quad I_n\otimes A_1, \quad I_n\otimes A_2, \quad A_1^t\otimes I_n, \quad A_2^t\otimes I_n\end{small},  \nonumber\\
\label{19.4.6-p2}\!\!\!\!&&\!\!\!\! \frac{1}{\sqrt{\sum_{\alpha\beta} g^{\alpha \beta} \xi_\alpha \xi_\beta}} \big(I_n\!\otimes \!
\!A_1\big)\big(A_1^t\!\otimes \!I_n\big), \quad \frac{1}{\sqrt{\sum_{\alpha\beta} g^{\alpha \beta} \xi_\alpha \xi_\beta}} \big(I_n\otimes \!
A_1\big)\big(A_2^t\otimes I_n\big),\nonumber\\
\!\!\!\!&& \!\!\!\! \frac{1}{\sqrt{\sum_{\alpha\beta} g^{\alpha \beta} \xi_\alpha \xi_\beta}} \big(I_n\!\otimes \!
A_2\big)\big(A_1^t\!\otimes \!I_n\big),
\quad  \frac{1}{\sqrt{\sum_{\alpha\beta}\! g^{\alpha \beta} \xi_\alpha \xi_\beta}} \big(I_n\!\otimes \!
A_2\big)\big(A_2^t\!\otimes \!I_n\big)\!\bigg\}.\nonumber\end{eqnarray}
 More precisely, for any $\theta(d_1, d_2; d_3, d_4)\in \mathfrak{J}$ and any $M\in H$,  we have the left multiplication of matrices about $\theta(d_1, d_2; d_3, d_4)$ to $M$, and $\big(\theta(d_1, d_2; d_3, d_4)\big) M\in \mathfrak{X}$; and for any two elements of $\mathfrak{X}$ we have the usual matrix addition and matrix multiplication in $\mathfrak{X}$. Clearly, $\mathfrak{X}$ is a nine-dimensional matrix ring with respect to $\mathfrak{J}$, and $H$ is a basis of $\mathfrak{X}$.
This implies that we should look for the inverse $U^{-1}$ of the form:
\begin{eqnarray}\label{19.3.26-6}
&& \;\;\quad U^{-1} \!=\! \frac{1}{\sum\limits_{\alpha, \beta} \!g^{\alpha\beta} \!\xi_\alpha \xi_\beta} \!\left\{\! {\tilde{s}}_1 \sqrt{\sum_{\alpha,\beta} g^{\alpha\beta}\xi_\alpha\xi_\beta}  I_{n^2}\right.\\
&&\!\!\!\left.+ I_n\otimes
   \left[\!\begin{BMAT}(@, 0pt, 0pt){c.c}{c.c} \!\left[\! {\tilde{s}}_2 \! \frac{\sum_{\alpha} \!g^{j\alpha} \xi_\alpha \xi_k}{\sqrt{\sum_{\alpha,\beta} g^{\alpha \beta} \xi_\alpha \xi_\beta}} \!\right]_{(n\!-\!1)\times(n\!-\!1)} &  \left[0\right]_{(n-1)\times 1} \\
 \left[0\right]_{1\times (n\!-\!1)} & {\tilde{s}}_3 \sqrt{\sum_{\alpha, \beta}\! g^{\alpha \beta}\! \xi_\alpha \xi_\beta }\end{BMAT}\! \right]\right.\nonumber
\\ &&\!\!+   I_n\!\otimes\!\begin{small}\left[\!\begin{BMAT}(@, 0pt, 0pt){c.c}{c.c} \!\left[ 0\right]_{(n\!-\!1)\times (n\!-\!1)}& \left[\!{\tilde{s}}_4 \!\sum_{\alpha}\! g^{j\alpha}\! \xi_\alpha \!\right]_{(n\!-\!1)\times 1}\\
\left[{\tilde{s}}_5\xi_k\right]_{1\times (n\!-\!1)}&0  \end{BMAT}\right]\end{small}\nonumber
\\
&& \!+ \begin{small}\left[\!\begin{BMAT}(@, 0pt, 0pt){c.c}{c.c} \begin{small}\left[\!
 {\tilde{s}}_6\frac{\sum_{\alpha} g^{k\alpha} \xi_\alpha \xi_j}{\sqrt{\sum_{\alpha, \beta} g^{\alpha \beta} \xi_\alpha \xi_\beta}}\right]_{(n\!-\!1)\times (n\!-\!1)} \end{small} & \left[0\right]_{(n\!-\!1)\times 1}\\
   \left[ 0\right]_{1\times (n\!-\!1)} & {\tilde{s}}_7 \sqrt{\sum_{\alpha, \beta} g^{\alpha \beta} \xi_\alpha \xi_\beta}  \end{BMAT}\!\right] \end{small}\otimes I_n\nonumber\\
   && \!  +
   \begin{small}\left[\begin{BMAT}(@, 0pt, 0pt){c.c}{c.c} \left[ 0\right]_{(n-1)\times (n-1)}& \left[{\tilde{s}}_9 \xi_j \right]_{(n-1)\times 1}\\
  \left[{\tilde{s}}_8\sum_{\alpha} g^{k\alpha} \xi_\alpha \right]_{1\times (n-1)}&0  \end{BMAT}\right]\end{small}\otimes I_n \nonumber\\
  && +\frac{\theta({\tilde{s}}_{10},{\tilde{s}}_{11}; {\tilde{s}}_{12}, {\tilde{s}}_{13})}{\sqrt{\sum\limits_{\alpha, \beta} g^{\alpha\beta} \xi_\alpha \xi_\beta}} \left( I_n \otimes
    A_1\right) \left( A_1^t \otimes I_n\right)
    \nonumber\\
     && \, +  \frac{\theta({\tilde{s}}_{14},{\tilde{s}}_{15}; {\tilde{s}}_{16}, {\tilde{s}}_{17})}{\sqrt{\sum\limits_{\alpha, \beta} g^{\alpha\beta} \xi_\alpha \xi_\beta}} \left( I_n \otimes
    A_1\right) \left( A_2^t \otimes I_n\right) \nonumber\\
    && \,+ \frac{\theta({\tilde{s}}_{18},{\tilde{s}}_{19}; {\tilde{s}}_{20}, {\tilde{s}}_{21})}{\sqrt{\sum\limits_{\alpha, \beta} g^{\alpha\beta} \xi_\alpha \xi_\beta}} \left( I_n \otimes
    A_2\right) \left( A_1^t \otimes I_n\right)\nonumber\\
     && \left.\,+  \frac{\theta({\tilde{s}}_{22},{\tilde{s}}_{23}; {\tilde{s}}_{24}, {\tilde{s}}_{25})}{\sqrt{\sum\limits_{\alpha, \beta} g^{\alpha\beta} \xi_\alpha \xi_\beta}} \left( I_n \otimes
    A_2\right) \left( A_2^t \otimes I_n\right)\right\}\nonumber\\
 \!\!&&\!\!= \frac{1}{\sum\limits_{\alpha, \beta} g^{\alpha\beta} \xi_\alpha \xi_\beta} \left\{ {\tilde{s}}_1 \sqrt{\sum_{\alpha,\beta} g^{\alpha\beta}\xi_\alpha\xi_\beta} \, I_{n^2} +\theta({\tilde s}_2, {\tilde{s}}_3; {\tilde s}_2, {\tilde{s}}_3) \big(I_n\otimes A_1\big)
 \right.\nonumber\\
 &&\left.\,\, + \,\theta( {\tilde s}_4, {\tilde{s}}_5; {\tilde s}_4, {\tilde{s}}_5) \big(I_n\otimes A_2\big)\right.+ \theta({\tilde{s}}_6, {\tilde{s}}_6; {\tilde{s}}_7,{\tilde{s}}_7)\big(A_1\otimes I_n\big) \nonumber \\
   &&\,\,\,+  \theta({\tilde{s}}_9, {\tilde{s}}_9; {\tilde{s}}_8,{\tilde{s}}_8)\big(A_2^t\otimes I_n\big) +
  \frac{\theta({\tilde{s}}_{10},{\tilde{s}}_{11}; {\tilde{s}}_{12}, {\tilde{s}}_{13})}{\sqrt{\sum\limits_{\alpha, \beta} g^{\alpha\beta} \xi_\alpha \xi_\beta}} \left( I_n \otimes
    A_1\right) \left( A_1^t \otimes I_n\right)\quad
    \nonumber\\
     && \,\,\,+ \frac{\theta({\tilde{s}}_{14},{\tilde{s}}_{15}; {\tilde{s}}_{16}, {\tilde{s}}_{17})}{\sqrt{\sum\limits_{\alpha, \beta} g^{\alpha\beta} \xi_\alpha \xi_\beta}} \left( I_n \otimes
    A_1\right) \left( A_2^t \otimes I_n\right) \nonumber\\
    &&\, \,\,+ \frac{\theta({\tilde{s}}_{18},{\tilde{s}}_{19}; {\tilde{s}}_{20}, {\tilde{s}}_{21})}{\sqrt{\sum\limits_{\alpha, \beta} g^{\alpha\beta} \xi_\alpha \xi_\beta}} \left( I_n \otimes
    A_2\right) \left( A_1^t \otimes I_n\right)\quad\nonumber\\
    &&\left.\,\,\,+ \frac{\theta({\tilde{s}}_{22},{\tilde{s}}_{23}; {\tilde{s}}_{24}, {\tilde{s}}_{25})}{\sqrt{\sum\limits_{\alpha, \beta} g^{\alpha\beta} \xi_\alpha \xi_\beta}} \left( I_n \otimes
    A_2\right) \left( A_2^t \otimes I_n\right)\right],\nonumber
    \end{eqnarray}
where ${\tilde{s}}_1, \cdots, {\tilde{s}}_{25}$ are the undetermined constants.
Let us note that \begin{eqnarray} &&\begin{small}\!\!\!\! \!\left(\! I_n \!\otimes\! \begin{small} \left[\begin{BMAT}(@, 0pt, 0pt){c.c}{c.c} \begin{small}\left[a\!\frac{\sum_{\alpha} g^{j\alpha} \xi_\alpha \xi_k}{\sqrt{\sum_{\alpha,\beta} g^{\alpha \beta} \xi_\alpha \xi_\beta}} \!\right]_{(n\!-\!1)\times(n\!-\!1)}\end{small} &  \left[0\right]_{(n\!-\!1)\times 1} \\
 \left[0\right]_{1\times (n\!-\!1)} & b \sqrt{\sum_{\alpha, \beta} g^{\alpha \beta} \xi_\alpha \xi_\beta }\end{BMAT} \! \right]\end{small}\right)\end{small}\!\theta(d_1, d_2; d_3, d_4) \nonumber\\
 && \!=\!\theta(ad_1, b d_2; a d_3, b d_4)\left(I_n\!\otimes \! A_1\right),  \nonumber\end{eqnarray}
 \begin{eqnarray}
  &&\!\!\! \!\!\!\!\begin{small}\left( I_n \otimes \begin{small}\left[\begin{BMAT}(@, 1pt, 1pt){c.c}{c.c} \left[
 0\right]_{(n-1)\times (n-1)}  & \left[a \sum_{\alpha} g^{j\alpha} \xi_\alpha \right]_{(n-1)\times 1}\\
   \left[ b\xi_k \right]_{1\times (n-1)} &0 \end{BMAT} \right]\end{small}\right) \end{small}\theta(d_1, d_2; d_3, d_4) \nonumber \\
   && = \theta(ad_2, b d_1; a d_4, bd_3) \left(I_n\otimes A_2\right), \qquad \quad \qquad\nonumber
    \end{eqnarray}
     \begin{eqnarray}&&\!\!\!\! \!\begin{small}\left( \begin{small} \left[\begin{BMAT}(@, 0pt, 0pt){c.c}{c.c} \begin{small}\left[a\frac{\sum_{\alpha} g^{k\alpha} \xi_\alpha \xi_j}{\!\sqrt{\!\sum_{\alpha,\beta} g^{\alpha \beta} \xi_\alpha \xi_\beta}} \right]_{(n\!-\!1)\times(n\!-\!1)}\end{small} &  \left[0\right]_{(n\!-\!1)\times 1} \\
 \left[0\right]_{1\times (n\!-\!1)} & b \sqrt{\sum_{\alpha, \beta} g^{\alpha \beta} \xi_\alpha \xi_\beta }\end{BMAT}\! \right]\end{small}\!\otimes\! I_n \!\right)\end{small}\theta(d_1, d_2; d_3, d_4) \nonumber \\
 && \!=\!\theta(ad_1, a d_2; b d_3, b d_4)\left(\!A_1^t\!\otimes \!I_n\right),  \nonumber\end{eqnarray}
 \begin{eqnarray}
  && \!\!\!\!\begin{small}\left(\begin{small}\left[\begin{BMAT}(@, 1pt, 1pt){c.c}{c.c} \left[ 0\right]_{(n-1)\times (n-1)}  & \left[a \xi_j\right]_{(n-1)\times 1}\\
   \left[ b \sum_{\alpha} g^{k\alpha} \xi_\alpha \right]_{1 \times (n-1)} &0 \end{BMAT} \right]\end{small}\otimes I_n\right)\end{small} \theta(d_1, d_2; d_3, d_4)\nonumber\\
   && = \theta(ad_3, a d_4; b d_1, bd_2) \left( A_2^t\otimes I_n\right). \qquad \quad\quad \;\, \qquad\nonumber
    \end{eqnarray}
     By inserting (\ref{19.3.23-15}) and (\ref{19.3.26-6}) into
 \begin{eqnarray} \label{19.3.26-7} UU^{-1} =I_{n^2},\end{eqnarray}
 and by comparing the ``coefficients matrix'' in front of each elements of $H$,
we get the following system of linear equations with unknown constants ${\tilde{s}}_1,  \cdots, {\tilde{s}}_{25}$.
\begin{align*}
  & 2s_1 {\tilde{s}}_1=1, \\
&s_2 {\tilde{s}}_1+ \big(2s_1+s_2\big) {\tilde{s}}_2 +  \big(s_4- \frac{i(\lambda+\mu)}{\mu}\big) {\tilde{s}}_5=0,\\
 &s_3 {\tilde{s}}_1 +\big( 2s_1+s_3\big) {\tilde{s}}_3 + \big(s_5- \frac{i(\lambda+\mu)}{\lambda+2\mu}\big) {\tilde{s}}_4=0,\qquad\qquad\qquad\qquad \qquad\qquad \qquad \quad \quad\\
& \big(s_4- \frac{i(\lambda+\mu)}{\mu}\big) {\tilde{s}}_1  + \big(s_4- \frac{i(\lambda+\mu)}{\mu}\big) {\tilde{s}}_3 + \big(2s_1 +s_2\big) {\tilde{s}}_4 =0,\\
& \big(s_5-\frac{i(\lambda+\mu)}{\lambda+2\mu} \big){\tilde{s}}_1+\big(s_5-\frac{i(\lambda+\mu)}{\lambda+2\mu} \big){\tilde{s}}_2
+\big(2s_1+s_3\big) {\tilde{s}}_5=0,\\
&  s_2 {\tilde{s}}_1 +\big(2s_1+s_2\big)  {\tilde{s}}_6+s_5   {\tilde{s}}_8=0,\\
&  s_3   {\tilde{s}}_1
+ \big(2s_1+s_3\big)   {\tilde{s}}_7+s_4  {\tilde{s}}_9=0,\\
& s_5 {\tilde{s}}_1 +s_5 {\tilde{s}}_7
+ (2s_1+s_2\big)   {\tilde{s}}_9=0,\\
 & s_4 {\tilde{s}}_1 +s_4 {\tilde{s}}_6
+ \big(2s_1+s_3\big)   {\tilde{s}}_8=0,\\
&  s_2 {\tilde{s}}_2  + s_2 {\tilde {s}}_6 + 2\big(s_1+s_2 \big) {\tilde {s}}_{10} +s_5 {\tilde{s}}_{16} + \big( s_4- \frac{i(\lambda+\mu)}{\mu} \big) {\tilde {s}}_{19} =0,\\
&  s_2 {\tilde{s}}_3  + s_3 {\tilde{s}}_6  +(2s_1+s_2+s_3) {\tilde {s}}_{11}+s_5 {\tilde {s}}_{17} + \big(s_5- \frac{i(\lambda+\mu)}{\lambda+2\mu} \big) {\tilde{s}}_{18}=0,\quad\,\qquad\qquad\quad\;\\
& s_3 {\tilde s}_2 + s_2 {\tilde{s}}_7 + \big(2s_1+s_2+s_3\big) {\tilde{s}}_{12} +s_4{\tilde {s}}_{14} +\big(s_4-\frac{i(\lambda+\mu)}{\mu}\big) {\tilde{s}}_{21} =0,\\
 &s_3 {\tilde s}_3 + s_3 {\tilde{s}}_{7} + 2(s_1+s_3) {\tilde{s}}_{13} + s_4{\tilde {s}}_{15} +\big(s_5-\frac{i(\lambda+\mu)}{\lambda+2\mu}\big) {\tilde{s}}_{20} =0,\\
   &  \,s_5 {\tilde{s}}_2  + s_2 {\tilde {s}}_9 + s_5 {\tilde{s}}_{12} +2(s_1+s_2) {\tilde{s}}_{14} +  \big( s_4- \frac{i(\lambda+\mu)}{\mu} \big) {\tilde {s}}_{23} =0,\\
& \, s_5 {\tilde{s}}_3  + s_3 {\tilde{s}}_9  +s_5 {\tilde{s}}_{13} +(2s_1+s_2+s_3) {\tilde {s}}_{15}+  \big(s_5- \frac{i(\lambda+\mu)}{\lambda+2\mu} \big) {\tilde{s}}_{22}=0,\qquad\qquad\quad\;\\
& \,s_4{\tilde{s}}_{2} + s_2 {\tilde{s}}_{8}+ s_4{\tilde{s}}_{10}+ \big(2s_1+s_2+s_3\big) {\tilde{s}}_{16}  +\big(s_4-\frac{i(\lambda+\mu)}{\mu}\big) {\tilde{s}}_{25} =0,\\
 & \,s_4{\tilde{s}}_{3} +s_3 {\tilde s}_8+ s_4 {\tilde{s}}_{11} +2(s_1+s_3) {\tilde{s}}_{17} + \big(s_5-\frac{i(\lambda+\mu)}{\lambda+2\mu}\big) {\tilde{s}}_{24} =0,\\
      & s_2 {\tilde{s}}_4 \! +\! \big( s_4\!-\!\frac{i(\lambda\!+\!\mu)}{\mu}\big) {\tilde {s}}_6
\!    +\! \big( s_4\!-\!\frac{i(\lambda+\mu)}{\mu}\big) {\tilde {s}}_{11}\!+\! 2(s_1\!+\!s_2){\tilde{s}}_{18} \!+\!s_5 {\tilde{s}}_{24}=0,\\
     & s_2 {\tilde{s}}_5 \! +\! \big( s_5\!-\!\frac{i(\lambda\!+\!\mu)}{\lambda\!+\!2\mu}\big) {\tilde {s}}_6 \!+\!
    \big( s_5\!-\!\frac{i(\lambda\!+\!\mu)}{\lambda\!+\!2\mu}\big) {\tilde {s}}_{10}\!+\! (2s_1\!+\!s_2\!+\!s_3){\tilde{s}}_{19} \!+\!s_5{\tilde{s}}_{25}\!=\!0,\quad\;\qquad\qquad\quad\;
   \\ &  s_3 {\tilde{s}}_{4} \!+\!  \big( s_4\!-\! \frac{i(\lambda\!+\!\mu)}{\mu} \big) {\tilde {s}}_{7} \!+\! \big( s_4\!-\! \frac{i(\lambda\!+\!\mu)}{\mu} \big) {\tilde {s}}_{13}
\!+\!(2s_1\!+\!s_2\!+\!s_3) {\tilde{s}}_{20}\!+\!s_4{\tilde{s}}_{22}\!=\!0,\\
& s_3 {\tilde{s}}_{5} +  \big( s_5- \frac{i(\lambda+\mu)}{\lambda+2\mu} \big) {\tilde {s}}_{7} + \big( s_5- \frac{i(\lambda+\mu)}{\lambda+2\mu} \big) {\tilde {s}}_{12}
+2(s_1+s_3) {\tilde{s}}_{21}+s_4{\tilde{s}}_{23}=0,  \\
     &  s_5 {\tilde{s}}_4  + \big( s_4-\frac{i(\lambda+\mu)}{\mu} \big){\tilde {s}}_9
    + \big( s_4-\frac{i(\lambda+\mu)}{\mu}\big) {\tilde {s}}_{15} + s_5 {\tilde{s}}_{20} + 2(s_1+s_2){\tilde{s}}_{22} =0,\\
     &  s_5 {\tilde{s}}_5 \!+\! \big( s_5\!-\!\frac{i(\lambda\!+\!\mu)}{\lambda\!+\!2\mu}\big) {\tilde {s}}_9 \!+\!
    \big( s_5\!-\!\frac{i(\lambda\!+\!\mu)}{\lambda\!+\!2\mu} \big){\tilde {s}}_{14}\!+\! s_5 {\tilde{s}}_{21}\! +\!(2s_1\!+\!s_2\!+\!s_3){\tilde{s}}_{23} \!=\!0,\quad \;\;\quad
   \\ &   s_4 {\tilde{s}}_{4} \!+\! \big( s_4\!-\! \frac{i(\lambda\!+\!\mu)}{\mu} \big) {\tilde {s}}_{8} \!+\! \big( s_4\!-\! \frac{i(\lambda\!+\!\mu)}{\mu} \big) {\tilde {s}}_{17}\!+\!s_4{\tilde{s}}_{18}
\!+\!(2s_1\!+\!s_2\!+\!s_3) {\tilde{s}}_{24}\!=\!0,\\
&  s_4 {\tilde{s}}_5 + \big( s_5- \frac{i(\lambda+\mu)}{\lambda+2\mu} \big) {\tilde {s}}_{8}+   \big( s_5- \frac{i(\lambda+\mu)}{\lambda+2\mu} \big) {\tilde {s}}_{16}+s_4{\tilde{s}}_{19}+2(s_1+s_3) {\tilde{s}}_{25}=0.
    \end{align*}
This system has a unique solution:
\begin{align} & {\tilde{s}}_1=	\frac{1}{2},\nonumber\\
	&- {\tilde{s}}_2={\tilde{s}}_3=\frac{\lambda +\mu}{\text{4\,}( \lambda +\text{3}\mu )}, \quad\qquad\quad \qquad\;\nonumber\\
&	{\tilde{s}}_4= \frac{i(\lambda+2\mu)(\lambda+\mu)}{4\mu( \lambda+3\mu)},\qquad \quad \qquad\;\nonumber\\
&	 {\tilde{s}}_5=\frac{i\mu(\lambda+\mu)}{4(\lambda+2\mu)(\lambda+3\mu)},\qquad \quad \qquad\;\nonumber\\
&	-{\tilde{s}}_6= {\tilde{s}}_7= \frac{\lambda +\mu}{\text{4}( \lambda +\text{3}\mu )},\nonumber\\
&	{\tilde{s}}_8={\tilde{s}}_9=- \frac{i(\lambda +\mu)}{4( \lambda +3\mu)},\nonumber\\
\label{2019.10.24-2}   &	-{\tilde{s}}_{10}={\tilde{s}}_{11}={\tilde{s}}_{12}=-{\tilde{s}}_{13}=\frac{( \lambda +\mu ) ^2}{4( \lambda +\text{3}\mu ) ^2},\quad \quad\quad\qquad  \\
	&	-{\tilde{s}}_{14}={\tilde{s}}_{15}=-{\tilde{s}}_{16}={\tilde{s}}_{17}= \frac{i( \lambda +\mu ) ^2}{\text{4}( \lambda +\text{3\,}\mu) ^2},\nonumber\\
& {\tilde{s}}_{18}=-	{\tilde{s}}_{20}= \frac{\text{i\,}( \lambda +\text{2}\mu )( \lambda+\mu ) ^2}{\text{4}\mu ( \lambda +\text{3}\mu ) ^2},\quad\;\,\nonumber\\
& {\tilde{s}}_{19}=-{\tilde{s}}_{21}=	\frac{i\mu ( \lambda +\mu ) ^2}{\text{4}( \lambda +\text{2}\mu ) ( \lambda +\text{3}\mu ) ^2},\nonumber\\
& {\tilde{s}}_{22}={\tilde{s}}_{24}=	-\frac{( \lambda +\text{2}\mu ) ( \lambda +\mu ) ^2}{4\mu ( \lambda +\text{3}\mu ) ^2\,},\qquad \quad\quad\nonumber \\
& {\tilde{s}}_{23}=	{\tilde{s}}_{25}=-\frac{\mu ( \lambda +\mu ) ^2}{4( \lambda +\text{2}\mu ) ( \lambda +\text{3}\mu ) ^2},\nonumber\end{align}
from which we explicitly obtain the unique inverse $U^{-1}$ of $U$.
Noting that \begin{align*} & {\tiny \theta ({\tilde{s}}_{10}, {\tilde{s}}_{11}; {\tilde{s}}_{12}, {\tilde{s}}_{13})\!=\!\!\left(\!\!I_n \!\otimes\!\!  \left[\!\begin{BMAT}(@, 0pt, 0pt){c.c}{c.c} \left[{\tilde{s}}_{10}\, \delta_{jk}\!\right]_{(n\!-\!1)\times (n\!-\!1)} &
 \left[ 0\right]_{(n\!-\!1)\times 1}\\
 \left[ 0\right]_{1\times (n\!-\!1)} & {\tilde{s}}_{11}\end{BMAT} \!\right]\!\right)\!\!\left( \!\left[\!\begin{BMAT}(@, 0pt, 0pt){c.c}{c.c} \left[ \delta_{jk}\right]_{(n\!-\!1)\times (n\!-\!1)} &
 \left[ 0\right]_{(n\!-\!1)\times 1}\\
 \left[ 0\right]_{1\times (n\!-\!1)} & -1\end{BMAT} \!\right]\!\!\otimes \!I_n\!\!\right)\!,}\\
 & \theta ({\tilde{s}}_{14}, {\tilde{s}}_{15}; {\tilde{s}}_{16}, {\tilde{s}}_{17})= \left(I_n \otimes  \left[\begin{BMAT}(@, 0pt, 0pt){c.c}{c.c} \left[{\tilde{s}}_{14}\, \delta_{jk}\right]_{(n-1)\times (n-1)} &
 \left[ 0\right]_{(n-1)\times 1}\\
 \left[ 0\right]_{1\times (n-1)} & {\tilde{s}}_{15}\end{BMAT} \right]\right)\left( I_n \otimes I_n\right),\\
 &  \theta ({\tilde{s}}_{18}, {\tilde{s}}_{19}; {\tilde{s}}_{20}, {\tilde{s}}_{21})\!\!=\!\!{\tiny \left(\!\!I_n\! \otimes \!\! \left[\!\begin{BMAT}(@, 0pt, 0pt){c.c}{c.c} \left[{\tilde{s}}_{18}\, \delta_{jk}\right]_{(\!n\!-\!1)\times (n\!-\!1)} &
 \left[ 0\right]_{(n\!-\!1)\times 1}\\
 \left[ 0\right]_{1\times (n\!-\!1)} & {\tilde{s}}_{19}\end{BMAT} \right]\!\right)\! \!\left(\! \left[\begin{BMAT}(@, 0pt, 0pt){c.c}{c.c}\! \left[ \delta_{jk}\right]_{(n\!-\!1)\times (n\!-\!1)} &
 \left[ 0\right]_{(n\!-\!1)\times 1}\\
 \left[ 0\right]_{1\times (n\!-\!1)} & -1\end{BMAT} \right]\!\!\otimes \!I_n\!\!\right)},\\
&  \theta ({\tilde{s}}_{22}, {\tilde{s}}_{23}; {\tilde{s}}_{24}, {\tilde{s}}_{25})= \left(I_n \otimes  \left[\begin{BMAT}(@, 0pt, 0pt){c.c}{c.c} \left[{\tilde{s}}_{22}\, \delta_{jk}\right]_{(n\!-\!1)\times (n\!-\!1)} &
 \left[ 0\right]_{(n\!-\!1)\times 1}\\
 \left[ 0\right]_{1\times (n\!-\!1)} & {\tilde{s}}_{23}\end{BMAT} \!\right]\!\right)\!\left(I_n \otimes I_n\right),\end{align*}
  we find from this and (\ref{19.3.26-6}) that
 \begin{align} \label{19.6.13-1}&\begin{small}\mbox{vec}\, X\end{small}\!=\!{\tiny \frac{1}{ \sum\limits_{\alpha, \beta}\!g^{\alpha\beta}\! \xi_\alpha \xi_\beta}}\! {\tiny\left(\!\! I_n \!\otimes\!
\!  \left[\!\begin{BMAT}(@, 0pt, 0pt){c.c}{c.c}\! \left[{\tilde{s}}_1\sqrt{\sum\limits_{\alpha,\beta}\! g^{\alpha \beta} \!\xi_\alpha \xi_\beta}\!+\!  {\tilde{s}}_2  \frac{\sum\limits_{\alpha} g^{j\alpha} \xi_\alpha \xi_k}{\sqrt{\sum\limits_{\alpha,\beta} g^{\alpha \beta} \xi_\alpha \xi_\beta}} \right] & \! \left[{\tilde{s}}_4 \!\sum_{\alpha}\! g^{j\alpha} \!\xi_\alpha \!\right]_{(n\!-\!1)\times 1} \\
 \left[{\tilde{s}}_5 \xi_k \right]_{1\times (n\!-\!1)} & ({\tilde{s}}_1\!+\!{\tilde{s}}_3) \sqrt{\sum_{\alpha, \beta}\! g^{\alpha \beta}\!\xi_\alpha \xi_\beta }\end{BMAT}\! \right]\!\right)}\!(\mbox{vec}\,  E)\\
&\, +\frac{1}{\sum\limits_{\alpha, \beta} \!g^{\alpha \beta} \xi_\alpha \xi_\beta }{\tiny \left(\left[\begin{BMAT}(@, 0pt, 0pt){c.c}{c.c}  \left[
 {\tilde{s}}_6\frac{\sum\limits_{\alpha} g^{j\alpha} \xi_\alpha \xi_k}{\sqrt{\sum\limits_{\alpha, \beta} g^{\alpha \beta} \xi_\alpha \xi_\beta}}\right]_{(n-1)\times (n-1)}  & \left[{\tilde{s}}_9 \sum\limits_{\alpha} g^{j\alpha} \xi_\alpha\right]_{(n-1)\times 1}\\
   \left[{\tilde{s}}_8 \xi_k\right]_{1\times (n-1)} & {\tilde{s}}_7 \sqrt{\sum\limits_{\alpha, \beta} g^{\alpha \beta} \xi_\alpha \xi_\beta}  \end{BMAT}\right]  \!  \otimes I_n\!\right)}\!(\mbox{vec}\, E)
   \nonumber\\
   &\,  +{\tiny \frac{1}{\big(\sum\limits_{\alpha, \beta} g^{\alpha\beta} \xi_\alpha \xi_\beta\big)^{3/2}}\left\{\left(I_n\otimes \left[\begin{BMAT}(@, 0pt, 0pt){c.c}{c.c}\left[\frac{{\tilde{s}}_{10} \sum\limits_{\alpha} g^{j\alpha} \xi_\alpha \xi_k}{\sqrt{\sum\limits_{\alpha, \beta} g^{\alpha \beta} \xi_\alpha \xi_\beta}}\right]_{(n-1)\times (n-1)} &
     \left[{\tilde{s}}_{18} \sum_{\alpha} g^{j\alpha} \xi_\alpha\right]_{(n-1)\times 1} \\
     \left[ {\tilde{s}}_{19} \xi_k\right]_{1\times (n-1)} & {\tilde{s}}_{11} \sqrt{\sum\limits_{\alpha, \beta} g^{\alpha \beta} \xi_\alpha \xi_\beta}  \end{BMAT}\right] \right)   \right.}\nonumber \\
    &\, \quad \qquad \quad \quad  \qquad \quad  \times{\tiny \left(\left[\begin{BMAT}(@, 0pt, 0pt){c.c}{c.c} \left[ \frac{\sum\limits_{\alpha} g^{j\alpha} \xi_\alpha \xi_k}{\sqrt{\sum\limits_{\alpha, \beta} g^{\alpha \beta} \xi_\alpha \xi_\beta}}\right]_{(n-1)\times (n-1)} &
     \left[0\right]_{(n-1)\times 1} \\
     \left[ 0\right]_{1\times (n-1)} & - \sqrt{\sum_{\alpha, \beta} g^{\alpha \beta} \xi_\alpha \xi_\beta}  \end{BMAT}\right]\otimes I_n\right)} (\mbox{vec}\,E)\nonumber\\
      &\,  \left. +
   {\tiny  \frac{1}{\big(\sum\limits_{\alpha, \beta} g^{\alpha\beta} \xi_\alpha \xi_\beta\big)^{3/2}} \left(I_n\otimes \left[\begin{BMAT}(@, 0pt, 0pt){c.c}{c.c} \left[ \frac{{\tilde{s}}_{14}\sum_{\alpha} g^{j\alpha} \xi_\alpha \xi_k}{\sqrt{\sum\limits_{\alpha, \beta} g^{\alpha \beta} \xi_\alpha \xi_\beta}}\right]_{(n-1)\times (n-1)}&
     \left[{\tilde{s}}_{22} \sum_{\alpha} g^{j\alpha} \xi_\alpha\right]_{(n-1)\times 1} \\
     \left[ {\tilde{s}}_{23} \xi_k\right]_{1\times (n-1)} & {\tilde{s}}_{15} \sqrt{\sum\limits_{\alpha, \beta} g^{\alpha \beta} \xi_\alpha \xi_\beta}  \end{BMAT}\right] \right)} \right.\nonumber \\
     & \left.\, \quad \qquad \quad \quad  \qquad \quad\quad  \times {\begin{small} \left(   \left[\begin{BMAT}(@, 0pt, 0pt){c.c}{c.c} \left[ 0\right]_{(n-1)\times (n-1)} &
     \left[\xi_j \right]_{(n-1)\times 1} \\
     \left[\sum_{\alpha} g^{k\alpha} \xi_\alpha \right]_{1\times (n-1)} & 0  \end{BMAT}\right] \otimes I_n\right)\end{small}}(\mbox{vec}\, E)\right\}\qquad\quad\qquad \qquad  \qquad \qquad\qquad\qquad\nonumber\end{align}
so that  \begin{eqnarray} \label{19.3.28-1}{}\end{eqnarray}
\begin{eqnarray*} X\!\!\! \!\!&=\!\!\!\!\!& {\begin{small} \frac{1}{ \sum\limits_{\alpha, \beta} g^{\alpha\beta} \xi_\alpha \xi_\beta}\!
  \left[\!\begin{BMAT}(@, 0pt, 0pt){c.c}{c.c} \!\left[\!{\tilde{s}}_1\sqrt{\sum\limits_{\alpha,\beta} g^{\alpha \beta} \xi_\alpha \xi_\beta}\, \delta_{jk}\!+\!  {\tilde{s}}_2 \! \frac{\sum\limits_{\alpha}\! g^{j\alpha} \xi_\alpha \xi_k}{\sqrt{\sum\limits_{\alpha,\beta} \!g^{\alpha \beta} \xi_\alpha \xi_\beta}} \right]_{(n\!-\!1)\times(n\!-\!1)} &  \left[{\tilde{s}}_4 \sum_{\alpha} \!g^{j\alpha} \xi_\alpha \right]_{(n\!-\!1)\times 1} \\
 \left[{\tilde{s}}_5 \xi_k \right]_{1\times (n-1)} & ({\tilde{s}}_1\!+\!{\tilde{s}}_3) \sqrt{\sum_{\alpha, \beta} \!g^{\alpha \beta}\! \xi_\alpha \xi_\beta }\end{BMAT} \!\right]\end{small}} E\\
&&\!\!\! \!\! \!\!\!\!\! +{\begin{small} \frac{1}{\sum\limits_{\alpha, \beta}\! g^{\alpha \beta} \xi_\alpha \xi_\beta } \; E  \left[\begin{BMAT}(@, 0pt, 0pt){c.c}{c.c} \left[
 {\tilde{s}}_6\frac{\sum_{\alpha} \!g^{j\alpha} \xi_\alpha \xi_k}{\sqrt{\sum_{\alpha, \beta}\! g^{\alpha \beta} \xi_\alpha \xi_\beta}}\right]_{(n\!-\!1)\times (n-1)}  & \left[{\tilde{s}}_8 \sum_{\alpha}\! g^{j\alpha} \xi_\alpha\right]_{(n\!-\!1)\times 1}\\
   \left[{\tilde{s}}_9 \xi_k\right]_{1\times (n\!-\!1)} & {\tilde{s}}_7 \sqrt{\sum_{\alpha, \beta} \!g^{\alpha \beta} \xi_\alpha \xi_\beta}  \end{BMAT}\right] \end{small}}
   \nonumber\\
   &&\!\! \!\!\!\!\! \!\!\! +{\begin{small} \frac{1}{\big(\sum\limits_{\alpha, \beta} g^{\alpha\beta} \xi_\alpha \xi_\beta\big)^{3/2}}\left\{ \left[\begin{BMAT}(@, 0pt, 0pt){c.c}{c.c} \left[\frac{{\tilde{s}}_{10} \sum_{\alpha} g^{j\alpha} \xi_\alpha \xi_k}{\sqrt{\sum_{\alpha, \beta} g^{\alpha \beta} \xi_\alpha \xi_\beta}}\right]_{(n-1)\times (n-1)} &
     \left[{\tilde{s}}_{18} \sum_{\alpha} g^{j\alpha} \xi_\alpha\right]_{(n-1)\times 1} \\
     \left[ {\tilde{s}}_{19} \xi_k\right]_{1\times (n-1)} & {\tilde{s}}_{11} \sqrt{\sum_{\alpha, \beta} g^{\alpha \beta} \xi_\alpha \xi_\beta}  \end{BMAT}\right]  E  \right. \end{small}}\\
    &&  \qquad  \qquad \quad\,  \times \begin{small}\left[\begin{BMAT}(@, 0pt, 0pt){c.c}{c.c} \begin{small}\left[ \frac{\sum_{\alpha} g^{j\alpha} \xi_\alpha \xi_k}{\sqrt{\sum_{\alpha, \beta} g^{\alpha \beta} \xi_\alpha \xi_\beta}}\right]_{(n-1)\times (n-1)} \end{small}&
     \left[0\right]_{(n-1)\times 1} \\
     \left[ 0\right]_{1\times (n-1)} & - \sqrt{\sum_{\alpha, \beta} g^{\alpha \beta} \xi_\alpha \xi_\beta}  \end{BMAT}\right] \end{small}\nonumber\\
      &&\!\!\! \!\!\!\!\!\!\!\!\left. +
     {\begin{small} \left[\!\begin{BMAT}(@, 0pt, 0pt){c.c}{c.c} \left[ \!\frac{{\tilde{s}}_{14}\sum\limits_{\alpha}\! g^{j\alpha} \!\xi_\alpha \xi_k}{\sqrt{\sum\limits_{\alpha, \beta}\! g^{\alpha \beta} \xi_\alpha \xi_\beta}}\!\right]_{(n\!-\!1)\times (n\!-\!1)} &
     \left[{\tilde{s}}_{22} \sum\limits_{\alpha}\! g^{j\alpha} \xi_\alpha\right]_{(n\!-\!1)\times 1} \\
     \left[ {\tilde{s}}_{23} \xi_k\right]_{1\times (n\!-\!1)} & {\tilde{s}}_{15} \sqrt{\sum_{\alpha, \beta}\! g^{\alpha \beta}\! \xi_\alpha \xi_\beta}  \end{BMAT}\!\right] \! E\!   \left[\begin{BMAT}(@, 0pt, 0pt){c.c}{c.c} \!\left[ 0\right]_{(n\!-\!1)\times (n\!-\!1)}&
     \left[\sum_{\alpha} \!g^{j\alpha} \xi_\alpha \right]_{(n\!-\!1)\times 1} \\
     \left[ \xi_k\right]_{1\!\times \!(n\!-\!1)} & 0  \end{BMAT}\!\right]\end{small}}\!\right\}.\qquad\quad\qquad \qquad  \qquad \qquad\qquad\qquad\nonumber\end{eqnarray*}
Therefore, when replacing the matrix $E$ in (\ref{19.3.28-1}) by $E_1$, we immediately get $q_0$.

\vskip 0.2 true cm

Step 4. \  Furthermore, by considering the terms of degree zero in (\ref{19.3.19-4}), we have
\begin{eqnarray} \label{19.3.28-6}& \big(q_1-b_1\big)  q_{-1} +q_{-1}q_1 =E_0,\end{eqnarray}
where \begin{eqnarray} \label{19.6.1-1}&&\;\; E_0:=-q_0^2 +i\sum_{l=1}^{n-1} \big(\frac{\partial q_1}{\partial \xi_l} \frac{\partial q_0}{\partial x_l} +\frac{\partial q_0}{\partial \xi_l}\, \frac{\partial q_1}{\partial x_l}\big) \\
 && \qquad \quad + \frac{1}{2} \sum_{l,\gamma=1}^{n-1} \frac{\partial^2 q_1}{\partial \xi_l\partial \xi_\gamma}   \, \frac{\partial^2 q_1}{\partial x_l \partial x_\gamma}   +b_0 q_0-i\sum_{l=1}^{n-1} \frac{\partial b_1}{\partial \xi_l} \, \frac{\partial q_0}{\partial x_l} +\frac{\partial q_0}{\partial x_n} -c_0.\nonumber\end{eqnarray}
Generally, for $m\ge 1$ we get
\begin{eqnarray} \label{19.3.28-7} \qquad \quad (q_1-b_1)  q_{-m-1} +q_{-m-1}q_1 =E_{-m},\end{eqnarray}
 where \begin{eqnarray} \label{19.6.1-2}  &&  E_{-m}:= -\begin{small}\sum\limits_{\substack{-m\le j,k\le 1 \\|\vartheta| = j+k+m}} \frac{(-i)^{|\vartheta|}}{\vartheta!}\end{small} ( \partial_{\xi'}^{\vartheta} q_j ) ( \partial_{x'}^{\vartheta} q_k) +b_0 q_{-m}  \\
 &&\quad \qquad \;\,\; -i\sum_{l=1}^{n-1} \frac{\partial b_1}{\partial \xi_l}  \frac{\partial q_{-m}}{\partial x_l} +\frac{\partial q_{-m}}{\partial x_n}.\nonumber\end{eqnarray}
Replacing the matrices $E$ and $X$ by the above $E_{-m}$ and $q_{-m-1}$ in (\ref{19.3.28-1}), respectively, we explicitly get the all $q_{-m-1}$, $m\ge 0$. \qed

\vskip 0.18 true cm

We have obtained the full symbol $q(x, \xi')\sim \sum_{l\le 1} q_l$ of $Q$ from above proposition 3.1. This implies that modulo a smoothing operator, the pseudodifferential operator $Q$ have been obtained on $\partial \Omega$. Thus we have the following:
\vskip 0.19 true cm

\noindent{\bf Proposition 3.2.} \ {\it In the local boundary normal coordinates, the Dirichlet-to-Neumann map $\Xi_g$ can be represented as:
 \begin{eqnarray}\label{19.3.28-10}\end{eqnarray} \begin{eqnarray} &&\;\;{\small
 \Xi_g (u\big|_{\partial \Omega})\!= \!  \left\{\!\left[\!
 \begin{BMAT}(@, 6pt, 6pt){c.c}{c.c}\big[\mu \delta_{jk} \big]_{(n\!-\!1)\times (n\!-\!1)}&\big[0\big]_{(n\!-\!1)\times 1}\\
\big[0\big]_{1\times (n\!-\!1)}& \lambda\!+\!2\mu\end{BMAT}\!\right] Q \qquad  \qquad \qquad \qquad \right.}\nonumber\\
&&\left.\qquad \qquad \quad \;  -{\small \left[ \!\begin{BMAT}(@, 6pt, 6pt){c.c}{c.c} [0]_{(n\!-\!1)\times (n\!-\!1)}& \bigg[\mu\sum\limits_{\alpha}g^{j\alpha} \frac{\partial}{\partial x_\alpha}\bigg]_{(n\!-\!1)\times 1}\\
\bigg[\lambda\big(\frac{\partial} {\partial x_k} \!+\!\sum\limits_{\alpha}\!\Gamma_{k\alpha}^\alpha\big)\bigg]_{1\times (n\!-\!1)} & \lambda\sum\limits_{\alpha}\!\Gamma_{n\alpha}^\alpha\end{BMAT}\!\right]\!}\right\}\begin{bmatrix}u^1\\
  \vdots\\
 u^{n\!-\!1}\\
 u^n\end{bmatrix}\qquad \qquad \;\qquad \nonumber  \end{eqnarray}
modulo a  smoothing operator.}

\vskip 0.26 true cm

 \noindent  {\it Proof.} \  Let $(x', x_n)$ be local boundary normal coordinates, for $x_3\in [0,T]$. Since the principal symbol of ${\mathcal{L}}_g$   is a negative-definite matrix, the hyperplane $\{x_n=0\}$ is non-characteristic, and hence ${\mathcal{L}}_g$ is partially hypoelliptic with respect to this boundary (see p.$\,$107 of \cite{HormL}). Therefore, the solution $u$ of the elastic Lam\'{e} equations ${\mathcal{L}}_g u=0$ is smooth in normal variable, i.e., in boundary normal coordinates $(x',x_n)$ with $x_n\in [0,T]$, $u\in (C^\infty([0,T]; {\mathfrak{D}}' ({\Bbb R}^{n-1})))^n$ locally. From Proposition 3.1, we see that (\ref{18/12/22-1}) is locally equivalent to the following system of equations for $u,v\in (C^\infty([0,T]; {\mathfrak{D}}' ({\Bbb R}^{n-1})))^{n}$: \begin{eqnarray*} &&\bigg(\frac{\partial }{\partial x_n}\, I_n + Q\bigg) u=v, \quad \; u\big|_{x_n=0}=f,\\
 &&\bigg(\frac{\partial }{\partial x_n} I_n +B-Q\bigg)v =h\in (C^\infty ([0,T]\times \mathbb{R}^{n-1}))^n. \end{eqnarray*}
 Making the substitution $t=T-x_n$ for the second equation mentioned above (as done in \cite{LU}), we get a backwards generalized heat equation system:
 \begin{eqnarray*} \bigg(\frac{\partial}{\partial t} I_n\bigg) v -(-Q+B)v =-h.\end{eqnarray*}
    Since $u$ is smooth in the interior of $\Omega$ by interior regularity for elliptic operator ${\mathcal{L}}_g$, it follows that
    $v$ is smooth in the interior of $\Omega$, and hence
    $v\big|_{x_n=T}$ is smooth. In view of the principal symbol of $Q$ is strictly positive for any $\xi\ne 0$, we get that the solution operator for this heat equation system is smooth for $t>0$ (see p.$\,$134 of \cite{Tre}).
    Therefore,
    \begin{eqnarray*} \left(\frac{\partial }{\partial x_n}I_n\right) u +Qu= v\in (C^\infty ([0,T]\times \mathbb{R}^{n-1}))^n \end{eqnarray*}
    locally.  Setting $Rf= v\big|_{\partial \Omega} $, we immediately see that $R$ is a smoothing operator and
    \begin{eqnarray*} \bigg(\big(\frac{\partial}{\partial x_n} I_n \big)u\bigg)\bigg|_{\partial \Omega} = -Qu\big|_{\partial \Omega} +Rf. \end{eqnarray*}
Combining this and (\ref{18-9/22-4}),  we obtain (\ref{19.3.28-10}).  \qed

 \vskip 1.38 true cm

\section{Determining metric of manifold from the elastic Dirichlet-to-Neumann map}

\vskip 0.45 true cm

It follows from (\ref{19.3.28-10}) that the full symbol $\sigma(\Xi_g)$ of $\Xi_g$ is
\begin{eqnarray} \label{9.10-8} \sigma(\Xi_g) \sim p_1+p_0+p_{-1} +\cdots+ p_{-m}+\cdots, \end{eqnarray}
where \begin{eqnarray}\label{19.3.28-13} \end{eqnarray} \begin{eqnarray}  p_1=  \left[ \begin{BMAT}(@, 6pt, 6pt){c.c}{c.c} \big[\mu \delta_{jk}\big]_{(n-1)\times (n-1)} & \big[0\big]_{(n-1)\times 1} \\
\big[0\big]_{1\times (n-1)} &\lambda+2\mu\end{BMAT} \right] q_1- \left[   \begin{BMAT}(@, 6pt, 6pt){c.c}{c.c} \big[0\big]_{(n-1)\times (n-1)} & \bigg[i\mu \sum\limits_{\alpha} g^{j\alpha} \xi_\alpha\bigg]_{(n-1)\times 1} \\
\bigg[i\lambda  \xi_k\bigg]_{1\times (n-1)}& 0 \end{BMAT}\right],\nonumber\end{eqnarray}
\begin{eqnarray} \label{19.3.28-14}\end{eqnarray} \begin{eqnarray} p_0=\left[ \begin{BMAT}(@, 6pt, 6pt){c.c}{c.c} \big[ \mu \delta_{jk} \big]_{(n-1)\times (n-1)} &\big[0\big]_{(n-1)\times 1} \\
\big[0\big]_{1\times (n-1)} & \lambda+2\mu\end{BMAT}\right] q_0 -  \left[ \begin{BMAT}(@, 6pt, 6pt){c.c}{c.c} \bigg[ 0 \bigg]_{(n-1)\times (n-1)} &\big[0\big]_{(n-1)\times 1}\\
\bigg[ \lambda\sum\limits_{\alpha} \Gamma_{k\alpha}^\alpha \bigg]_{1\times (n-1)} & \lambda \sum\limits_{\alpha} \Gamma_{n\alpha}^\alpha \end{BMAT}\right],\nonumber\end{eqnarray}
\begin{eqnarray} \label{19.3.28-15}\end{eqnarray} \begin{eqnarray} p_{-m}= \left[\begin{BMAT}(@, 6pt, 6pt){c.c}{c.c}  \big[\mu \delta_{jk}\big]_{(n-1)\times (n-1)} & \big[0\big]_{(n-1)\times 1} \\
\big[0\big]_{1\times (n-1)} & \lambda +2\mu \end{BMAT}\right] q_{-m},  \quad\, \; m\ge 1. \nonumber \qquad \qquad \qquad \qquad\;\end{eqnarray}

\vskip 0.39 true cm

 \noindent{\bf Proposition 4.1.} \ {\it Suppose $\mbox{dim}\, M=n\ge 2$, and assume $\mu>0$ and $\lambda+\mu\ge 0$.
 Let $\{ x_1,\cdots, x_{n-1}\}$ be any local coordinates for an open set $W\subset \partial \Omega$,  and let $\{ p_j\}_{j\le 1} $ denote the full symbol of $\Xi_g$ in these coordinates. If $(n-1)\lambda^3+(4n-2) \lambda^2\mu +(n+5)\lambda \mu^2 +(14-8n)\mu^3 \neq 0$, then for any $x_0\in W$, the full Taylor series of $g$ at $x_0$ in boundary normal coordinates is given by explicit formula in terms of the matrix-valued functions $\{ p_j\}_{j\le 1}$
 and their tangential derivatives at $x_0$.}

\vskip 0.22 true cm

\noindent  {\it Proof.} \   Denote by $\{x_1,\cdots, x_n\}$ the boundary normal coordinates associated with $\{x_1,\cdots, x_{n-1}\}$ as in section 2. According to the form of metric (\ref{18/a-1}) which we have chosen,  we immediately see that it suffices to show that the matrix-valued functions $\{p_j\}$ determine the metric $[g_{\alpha\beta}]_{(n-1)\times (n-1)}$ and all its normal derivatives along $\partial \Omega$. Noticing that $\frac{\partial g_{\alpha\beta}}{\partial x_n} =-\sum\limits_{\rho, \gamma} g_{\alpha \rho} \frac{\partial g^{\rho \gamma}}{\partial x_n} g_{\gamma \beta},$ it is also enough to determine the inverse matrix $[g^{\alpha\beta}]_{(n-1)\times (n-1)}$ and all its normal derivatives.
In (\ref{19.3.28-13}), because the $(n,n)$-entry of the matrix in the second term  \begin{eqnarray*}\left[   \begin{BMAT}(@, 6pt, 6pt){c.c}{c.c} \big[0\big]_{(n-1)\times (n-1)} & \bigg[i\mu \sum\limits_{\alpha}g^{j\alpha}\xi_\alpha\bigg]_{(n-1)\times 1} \\ \bigg[i\lambda  \xi_k\bigg]_{1\times (n-1)}& 0 \end{BMAT}\right]\end{eqnarray*}
 is $0$, we get that $p_1$ uniquely determines the $(n,n)$-entry of the matrix  \begin{eqnarray*} \left[ \begin{BMAT}(@, 6pt, 6pt){c.c}{c.c} \big[\mu g_{jk}\big]_{(n-1)\times (n-1)} & \big[0\big]_{(n-1)\times 1} \\
\big[0\big]_{1\times (n-1)} &\lambda+2\mu\end{BMAT} \right] q_1,\end{eqnarray*}
where $q_1$ is given by (\ref{19.3.22-8}) and (\ref{19.4.6-8}). More precisely, $p_1$ uniquely determines
the $(n,n)$-entry of the matrix   \begin{eqnarray} \label{19.4.18-1.} &&  \left[ \begin{BMAT}(@, 6pt, 6pt){c.c}{c.c} \big[\mu \delta_{jk}\big]_{(n\!-\!1)\times (n\!-\!1)} & \big[0\big]_{(n\!-\!1)\times 1} \\
\big[0\big]_{1\times (n\!-\!1)} &\lambda\!+\!2\mu\end{BMAT} \right] q_1 \\
&&\; =
 {\begin{small} \left[\begin{BMAT}(@, 1pt, 1pt){c.c}{c.c} \left[\!s_1 \mu\sqrt{\sum_{\alpha, \beta} g^{\alpha\beta} \xi_\alpha\xi_\beta} \delta_{jk}
  \! +\! s_2 \mu \frac{ \sum_{\alpha} g^{j\alpha}\xi_\alpha \xi_k}{\sqrt{\sum_{\alpha, \beta} g^{\alpha\beta} \xi_\alpha\xi_\beta}} \! \right]_{(n-1)\times(n-1)} & \left[s_4\mu  \sum_\alpha g^{j\alpha} \xi_\alpha \right]_{(n-1)\times 1} \\
\left[s_5(\lambda+2\mu) \xi_k\right]_{1\times (n-1)}
   &(s_1 \! +\!s_3)(\lambda+2\mu)\sqrt{\sum_{\alpha,\beta} g^{\alpha\beta} \xi_\alpha\xi_\beta}  \end{BMAT}\right]. \end{small}}\nonumber
\end{eqnarray}
From the $(n,n)$-entry of the matrix on the right side of (\ref{19.4.18-1.}), we see that $p_1$ uniquely determines $\sqrt{\sum_{\alpha,\beta} g^{\alpha\beta} \xi_\alpha\xi_\beta}$, and hence $p_1$ uniquely determines $\sum_{\alpha,\beta} g^{\alpha\beta} \xi_\alpha\xi_\beta$ for all $(\xi_1, \cdots,\xi_{n-1})\in {\Bbb R}^{n-1}$. This implies that $p_1$ uniquely determines $g^{\alpha \beta}$ on $\partial \Omega$ for all $1\le \alpha, \beta\le n-1$.

\vskip 0.2 true cm

Next,  by (\ref{19.3.23-1}), (\ref{19.3.22-8}), (\ref{19.3.19-4'}), (\ref{19.3.19-5}) and (\ref{19.3.19-9}) we can rewrite $E_1$ as
\begin{align} \label{19.4.12,1}
 &     E_1=b_0q_1+\frac{\partial q_1}{\partial x_n} - c_1 +T_0^{(1)}(g_{\alpha\beta})\\
&       \!= \!{\tiny \left[\!\!\begin{BMAT}(@, 0pt, 0pt){c.c}{c.c}\bigg[\!\frac{1}{2} \!\sum\limits_{\alpha,\beta} \!g^{\alpha\!\beta} \frac{\partial g_{\alpha\!\beta}}{\partial x_n}\!\delta_{\!jk}\! +\!2\Gamma_{kn}^j \bigg]_{\!(n\!-\!1)\!\times\! (n\!-\!1)}\!
 & \big[0\big]_{(n\!-\!1)\!\times \!1}\\
 \!\big[0\big]_{1\!\times \!(n\!-\!1)\!} & \frac{1}{2} \!\sum_{\alpha,\beta} \!g^{\alpha\!\beta} \frac{\partial g_{\alpha\beta}}{\partial x_n}\end{BMAT}\! \!\right]}\! {\tiny\left[\!\!\begin{BMAT}(@, 0pt, 0pt){c.c}{c.c} \left[\!s_1 \sqrt{\sum\limits_{\alpha,\beta}\! g^{\alpha \!\beta} \!\xi_\alpha\!\xi_\beta}\!+ \! \frac{s_2\!\sum\limits_{\alpha}\! g^{\!j\!\alpha} \xi_\alpha\! \xi_k}{\sqrt{ \sum\limits_{\alpha,\beta}\! g^{\alpha\!\beta} \!\xi_\alpha\! \xi_\beta}} \!\right] & \left[s_4 \!\sum\limits_{\alpha} \!g^{j\!\alpha}\!\xi_\alpha\right]_{\!(n\!-\!1)\!\times\! 1} \\
  \left[s_5\xi_k\right]_{1\times(n-1)} & (\!s_1\!+\!s_3\!)\! \sqrt{\sum\limits_{\alpha,\beta}\! g^{\alpha\! \beta}\! \xi_\alpha\! \xi_\beta}\! \end{BMAT}\!\right]}\nonumber\\
&       + {\tiny \left[\!\!\begin{BMAT}(@, 0pt, 0pt){c.c}{c.c}\left[\! \frac{s_1 \sum\limits_{\alpha,\beta} \frac{\partial g^{\alpha\beta}}{\partial x_n} \xi_\alpha\xi_\beta}{2 \sqrt{\sum\limits_{\alpha,\beta} g^{\alpha\beta}
 \xi_\alpha\xi_\beta}} + \frac{s_2\! \sum\limits_{\alpha}\frac{\partial g^{j\alpha}}{\partial x_n} \xi_\alpha\xi_k}
 {\sqrt{\sum\limits_{\alpha,\beta} g^{\alpha\beta}  \xi_\alpha\xi_\beta}} - \frac{s_2 \big(\sum\limits_{\alpha} g^{j\alpha} \xi_\alpha\xi_k\big)\sum\limits_{\alpha,\beta} \frac{\partial g^{\alpha\beta}}{\partial x_n} \xi_\alpha\xi_\beta}{2 \big(\sum\limits_{\alpha,\beta} g^{\alpha\beta} \xi_\alpha\xi_\beta\big)^{3/2}}\right]_{\!(n\!-\!1)\times(n\!-\!1)\!}
& \bigg[s_4 \sum\limits_{\alpha} \frac{\partial g^{j\alpha}}{\partial x_n} \xi_\alpha\bigg]_{(n-1)\times 1}\\
\big[0\big]_{1\times (n-1)} & \frac{(s_1+s_3)\sum\limits_{\alpha,\beta} \frac{\partial g^{\alpha\beta}}{\partial x_n} \xi_\alpha\xi_\beta}{2\sqrt{\sum\limits_{\alpha,\beta} g^{\alpha\beta}\xi_\alpha\xi_\beta}}\end{BMAT}\right] }\nonumber
\\
  &    - {\tiny  \left[\begin{BMAT}(@, 0pt, 0pt){c.c}{c.c} \big[0\big]_{(n\!-\!1)\times (n\!-\!1)}  &  \bigg[\!\frac{i(\lambda\!+\!\mu)}{\mu}\! \sum\limits_{\alpha}\!
 g^{j\alpha} \!\sum\limits_{\beta}\! \Gamma_{n\beta}^\beta\! \xi_\alpha  \!+\!2i \!\sum\limits_{\alpha,\beta}\! g^{\alpha\beta}\!  \Gamma_{n\alpha}^j\! \xi_\beta \!\bigg]_{(n\!-\!1)\times 1} \\
 \left[\frac{2i\mu}{ \lambda\!+\!2\mu} \!\sum\limits_{\alpha,\beta} \!g^{\alpha\beta}\! \Gamma_{k\alpha}^n\! \xi_\beta \!\right]_{1\times (n\!-\!1)} &0\end{BMAT}\!\right]}\!+\!T_0^{(2)}(g_{\alpha\beta})\nonumber\\ \nonumber
 \\
& :=  E_1^{(1)}+  E_1^{(2)}+ E_1^{(3)}+E_1^{(4)}+T_0^{(2)}(g_{\alpha\beta}),\nonumber\end{align}
where \begin{align*}& {\tiny E_1^{(1)}\!=\!\frac{1}{2}\! \sum_{\alpha,\beta} \!g^{\alpha\beta}\! \frac{\partial g_{\alpha\beta}}{\partial x_n}\!
 \left[\begin{BMAT}(@, 1pt, 1pt){c.c}{c.c} \left[s_1 \sqrt{\sum_{\alpha,\beta}\! g^{\alpha \beta} \xi_\alpha\xi_\beta} \delta_{jk}\!+\! \frac{s_2 \sum_{\alpha} \!g^{j\alpha} \xi_\alpha \xi_k}{\sqrt{ \sum_{\alpha,\beta}\! g^{\alpha\beta} \xi_\alpha \xi_\beta}} \!\right]_{(n\!-\!1)\times(n\!-\!1)} & \left[s_4 \sum_{\alpha} \!g^{j\alpha}\xi_\alpha\!\right]_{(n\!-\!1)\times 1} \\
  \left[s_5 \xi_k\right]_{1\times (n\!-\!1)} & (s_1\!+\!s_3) \sqrt{\sum_{\alpha,\beta}\! g^{\alpha \beta} \xi_\alpha \xi_\beta}\! \end{BMAT}\!\right]}\nonumber
  \\
   &{\tiny E_1^{(2)}=   \left[\begin{BMAT}(@, 0pt, 0pt){c.c}{c.c} \left[2s_1 \sqrt{\sum_{\alpha,\beta} g^{\alpha \beta} \xi_\alpha\xi_\beta} \,  \,\Gamma_{kn}^j +\frac{2s_2\sum\limits_{\alpha,\beta} \Gamma_{\beta n}^j g^{\beta \alpha} \xi_\alpha\xi_k}{\sqrt{\sum\limits_{\alpha,\beta} g^{\alpha\beta}  \xi_\alpha\xi_\beta}}\right]_{(n-1)\times (n-1)}   & \left[2s_4 \sum\limits_{\alpha,\beta}\Gamma_{\beta n}^j g^{\beta \alpha} \xi_\alpha\right]_{(n-1)\times 1}\\
   \big[0\big]_{1\times (n-1)} & 0\end{BMAT}\right]}\nonumber\\
   &   {\tiny E_1^{(3)}\!=\!  \left[\!\!\begin{BMAT}(@, 0pt, 0pt){c.c}{c.c}\left[\! \frac{s_1 \sum\limits_{\alpha,\beta} \!\frac{\partial g^{\alpha\beta}}{\partial x_n} \!\xi_\alpha\xi_\beta}{2 \sqrt{\sum\limits_{\alpha,\beta}\! g^{\alpha\beta}
 \xi_\alpha\xi_\beta}} \!+\! \frac{s_2\! \sum\limits_{\alpha}\!\frac{\partial g^{j\alpha}}{\partial x_n} \xi_\alpha\xi_k}
 {\sqrt{\sum\limits_{\alpha,\beta} \!g^{\alpha\beta}  \xi_\alpha\xi_\beta}} \!-\! \frac{s_2 \big(\sum\limits_{\alpha,\beta} \frac{\partial g^{\alpha\beta}}{\partial x_n}\! \xi_\alpha\xi_\beta\big)\sum\limits_{\alpha} \!g^{j\alpha} \xi_\alpha\xi_k}{2 \big(\sum\limits_{\alpha,\beta} \!g^{\alpha\beta} \xi_\alpha\xi_\beta\big)^{3/2}}\!\right]_{\!(n\!-\!1)\times(n\!-\!1)\!}
& \bigg[s_4 \sum\limits_{\alpha} \!\frac{\partial g^{j\alpha}}{\partial x_n} \xi_\alpha\bigg]_{(n\!-\!1)\times 1}\\
\big[0\big]_{1\times (n\!-\!1)} & \frac{(s_1\!+\!s_3)\sum\limits_{\alpha,\beta} \!\frac{\partial g^{\alpha\beta}}{\partial x_n} \xi_\alpha\xi_\beta}{2\sqrt{\sum\limits_{\alpha,\beta}\! g^{\alpha\beta}\xi_\alpha\xi_\beta}}\end{BMAT}\!\right] }\nonumber
\\
  &   {\tiny  E_1^{(4)}=   -  \left[\begin{BMAT}(@, 0pt, 0pt){c.c}{c.c} \big[0\big]_{(n\!-\!1)\times (n\!-\!1)}  &  \bigg[\frac{i(\lambda\!+\!\mu)}{\mu}\! \sum\limits_{\alpha}\!
 g^{j\alpha} \sum\limits_{\beta} \Gamma_{n\beta}^\beta \xi_\alpha  +2i \sum\limits_{\alpha,\beta} g^{\alpha\beta}  \Gamma_{n\alpha}^j \xi_\beta \bigg]_{(n\!-\!1)\times 1} \\
 \left[\frac{2i\mu}{ \lambda+2\mu} \sum\limits_{\alpha,\beta} g^{\alpha\beta} \Gamma_{k\alpha}^n \xi_\beta \right]_{1\times (n-1)} &0\end{BMAT}\right]}, \nonumber\quad\;\qquad\end{align*}
 here and elsewhere each $T_0^{(s)}(g_{\alpha \beta})$ is a matrix expression involving only $g_{\alpha\beta}$, $g^{\alpha\beta}$, and their tangential derivatives along $\partial \Omega$.
  Replacing matrix $E$ in (\ref{19.3.28-1}) by each $E_1^{(l)}$, we get
 \begin{align*} &     q_0^{\!(1)\!}\!=\! \frac{1}{2}\!\sum\limits_{\alpha,\beta}\!  g^{\alpha\!\beta} \frac{\partial g_{\alpha\!\beta}}{\partial x_n}\! {\tiny  \left[\!\begin{BMAT}(@, 1pt, 1pt){c.c}{c.c} \left[ \;* \;\right]_{(\!n\!-\!1)\times (\!n\!-\!1)}& \left[\;\;*\;\;\right]_{(n\!-\!1)\times 1}\\
\left[\; *\;\right]_{1\times (n\!-\!1)} & s_4({\tilde{s}}_5\!-\!{\tilde{s}}_{19}) \!+\! (s_1\!+\!s_3) ({\tilde{s}}_1\!+\!{\tilde{s}}_3\!+\!{\tilde{s}}_7\!-\!{\tilde{s}}_{11}) \!+\!s_5({\tilde{s}}_8\!+\!{\tilde{s}}_{15})\!+\!(s_1\!+\!s_2){\tilde{s}}_{23}
\end{BMAT}\!\right]} \!+\!T_0^{\!(3)}\!(g_{\!\alpha\!\beta}),\\
&   {\tiny   q_0^{(2)}= \frac{\sum\limits_{\alpha,\beta} \frac{\partial g^{\alpha\beta}}{\partial x_n} \xi_\alpha\xi_\beta}{ \sum\limits_{\alpha,\beta} g^{\alpha\beta}\xi_\alpha\xi_\beta}   \left[\begin{BMAT}(@, 1pt, 1pt){c.c}{c.c} \left[ \;\,* \;\;\right]_{(n-1)\times (n-1)}& \left[\;\;*\;\;\right]_{(n-1)\times 1}\\
\left[\;\; *\;\;\right]_{1\times (n-1)} & s_4(-{\tilde{s}}_5+{\tilde{s}}_{19}) - (s_1+s_2){\tilde{s}}_{23}
\end{BMAT}\right]} +T_0^{(4)}(g_{\alpha\beta}),\\
&    {\tiny q_0^{(3)}\!=\! \frac{\sum\limits_{\alpha,\beta} \frac{\partial g^{\alpha\beta}}{\partial x_n} \xi_\alpha\xi_\beta}{ \sum\limits_{\alpha,\beta} g^{\alpha\beta}\xi_\alpha\xi_\beta}  \left[\begin{BMAT}(@, 1pt, 1pt){c.c}{c.c} \left[ \;\,* \;\;\right]_{(n\!-\!1)\times (n\!-\!1)}& \left[\;\;*\;\;\right]_{(n\!-\!1)\times 1}\\
\left[\;\; *\;\;\right]_{1\times (n\!-\!1)} & s_4({\tilde{s}}_5\!-\!{\tilde{s}}_{19})\! +\!\frac{1}{2} (s_1\!+\!s_3)( {\tilde{s}}_1 \!+ \! {\tilde{s}}_3\! +\! {\tilde{s}}_7\!-\!  {\tilde{s}}_{11})\! +\!\frac{1}{2}(s_1\!+\!s_2){\tilde{s}}_{23}
\end{BMAT}\right]} \!+\!T_0^{\!(5)}\!(g_{\alpha\beta}),\\
  &   {\tiny  q_0^{\!(4)}\!=\!\frac{1}{2}\!\sum\limits_{\alpha,\beta} \!g^{\alpha\!\beta} \frac{\partial g_{\alpha\!\beta}}{\partial x_n}\! \!\left[\begin{BMAT}(@, 0pt, 0pt){c.c}{c.c} \!\left[ \;\,* \;\;\right]_{(\!n\!-\!1)\times (\!n\!-\!1)}& \left[\;\;*\;\;\right]_{(\!n\!-\!1)\times \!1}\\
\left[\;\; *\;\;\right]_{1\times (\!n\!-\!1)} & \frac{i(\lambda\!+\!\mu)}{\mu} (-{\tilde{s}}_5 \!+\!{\tilde{s}}_{\!19}) \end{BMAT}\!\!\right]}\!\\
& \quad \;\,\quad  +{\tiny     \frac{\sum\limits_{\alpha,\beta} \!\frac{\partial g^{\alpha\!\beta}}{\partial x_n} \xi_\alpha\!\xi_\beta}{ \sum\limits_{\alpha,\beta}\! g^{\alpha\!\beta}\xi_\alpha\!\xi_\beta} \!  \left[\begin{BMAT}(@, 1pt, 1pt){c.c}{c.c}\! \left[ \;\,* \;\;\right]_{(n\!-\!1)\times (n\!-\!1)}& \left[\;\;*\;\;\right]_{(n\!-\!1)\times 1}\\
\left[\;\; *\;\;\right]_{1\times (n\!-\!1)}& i\big({\tilde{s}}_5 \!-\!{\tilde{s}}_{19} \big) \!-\!\frac{i\mu}{\lambda\!+\!2\mu} \big({\tilde{s}}_{8}\!+\! {\tilde{s}}_{15} \big)\!\end{BMAT}\!\right]}
 +T_0^{(6)}(g_{\alpha\beta}).\end{align*}
Therefore
 \begin{eqnarray}  \label{19.4.10-4}q_0\!=\! \sum\limits_{l=1}^4\! q_0^{(l)}\!=\!\begin{small} \left[\!\begin{BMAT}(@, 1pt, 1pt){c.c}{c.c} \left[ \;\,* \;\;\right]_{(n\!-\!1)\times (n\!-\!1)}& \left[\;\;*\;\;\right]_{(n\!-\!1)\times 1}\\
\left[\;\; *\;\;\right]_{1\times (n\!-\!1)} & \frac{ d_1}{2} \sum_{\alpha,\beta} \!g^{\alpha\beta} \frac{\partial g_{\alpha\beta}}{\partial x_n} \!+\! \frac{ d_2\sum_{\alpha, \beta}\!
 \frac{\partial g^{\alpha\beta}}{\partial x_n}  \xi_\alpha \xi_\beta}{\sum_{\alpha, \beta}\! g^{\alpha \beta} \xi_\alpha \xi_\beta} \end{BMAT}\!\right]\end{small} \!+\!T_0^{(7)}(g_{\alpha\beta}), \end{eqnarray}
where \begin{eqnarray*} \label{19.4.10-10} d_1\!\!\!\! &\!=\!\!
\bigg(\!\!s_4\!-\!\frac{i(\lambda\!+\!\mu)}{\mu}\!\bigg)({\tilde{s}}_5 \!-\!{\tilde{s}}_{19})\!+\! (s_1\!+\!s_3)({\tilde{s}}_1\!+\!{\tilde{s}}_3\!+\!{\tilde{s}}_7\!-\!{\tilde{s}}_{11})\!+\!s_5({\tilde{s}}_8\!+\! {\tilde{s}}_{15})  \!+\!(s_1\!+\!s_2){\tilde{s}}_{23}\\ \!\!& \!\! =\! \frac{{2}\mu \left( \lambda +{2}\mu \right)}{\left( \lambda +{3}\mu \right) ^2},\qquad \qquad \qquad \quad\qquad \qquad \qquad \quad \quad\qquad\qquad \qquad \quad \qquad \qquad\; \quad\nonumber\\
 d_2 \!\!\!\! \!\!&\!\! \!\!=\!\frac{1}{2}(s_1\!+\!s_3)({\tilde{s}}_1\!+\!{\tilde{s}}_3\!+\!{\tilde{s}}_7\!-\!{\tilde{s}}_{11}) \!+\!i\big({\tilde{s}}_5 \!-\!{\tilde{s}}_{19}\big) \!-\! \frac{1}{2} (s_1+s_2){\tilde{s}}_{23} -\frac{i\mu}{\lambda+2\mu}({\tilde{s}}_8 +{\tilde{s}}_{15}) \\
\!\!\!\!\!\!\! &\!\!\!\!\!\!\! =\!  \frac{\mu ( \lambda ^2+{3}\lambda \mu +{3}\mu ^2 )}{( \lambda \!+\!{2}\mu )( \lambda \!+\!{3}\mu ) ^2}.\quad\qquad\qquad \qquad \quad \qquad \quad\qquad\qquad \qquad\qquad \qquad \quad \quad \quad\end{eqnarray*}

It follows from the expression of $p_0$ in (\ref{19.3.28-14}) that, except for a known term (involving only $g^{\alpha\beta}$, $g_{\alpha\beta}$ and their tangential derivatives along $\partial \Omega$), the $(n,n)$-entry of matrix $p_0$ is
\begin{eqnarray}&& \label{19.4.10-6} (\lambda+2\mu) \left(\frac{d_1}{2}\sum_{\alpha,\beta} g^{\alpha\beta}\frac{\partial g_{\alpha\beta}}{\partial x_n}+ \frac{d_2\sum_{\alpha,\beta} \frac{\partial g^{\alpha\beta}}{\partial x_n}\xi_\alpha \xi_\beta}{\sum_{\alpha,\beta} g^{\alpha\beta} \xi_\alpha\xi_\beta}\right)- \frac{\lambda}{2} \sum_{\alpha,\beta} g^{\alpha \beta} \frac{\partial g_{\alpha\beta}}{\partial x_n}.\end{eqnarray}
Here we have used the fact that $\sum_{\alpha}\Gamma_{n\alpha}^\alpha =\frac{1}{2} \sum_{\alpha,\beta} g^{\alpha\beta} \frac{\partial g_{\alpha\beta}}{\partial x_n}$. Note that $\sum_{\alpha, \beta} g^{\alpha\beta} g_{\alpha\beta}=n-1$, and so
$\frac{1}{2} \sum_{\alpha, \beta} g^{\alpha \beta} \frac{\partial g_{\alpha \beta}}{\partial x_n}=-\frac{1}{2} \sum_{\alpha,\beta} g_{\alpha\beta}
\frac{\partial g^{\alpha \beta}}{\partial x_n}$.
If we set $h^{\alpha\beta}_1=\frac{\partial g^{\alpha\beta}}{\partial x_n}$, $h_1=\sum_{\alpha, \beta} g_{\alpha\beta} h^{\alpha\beta}_1$,
then (\ref{19.4.10-6}) (which is uniquely determined by $p_0$) can further be written as
  \begin{eqnarray}&& \label{19.4.10-7} \frac{\lambda-(\lambda+2\mu) d_1 }{2}  \sum_{\alpha, \beta} g_{\alpha \beta} \frac{\partial g^{\alpha\beta}}{\partial x_n} + (\lambda+2\mu)d_2\frac{ \sum_{\alpha, \beta}
 \frac{\partial g^{\alpha\beta}}{\partial x_n} \, \xi_\alpha \xi_\beta}{\sum_{\alpha, \beta} g^{\alpha \beta} \xi_\alpha \xi_\beta}\\
  &&  \; \; =\big(\sum_{\alpha,\beta}g^{\alpha\beta} \xi_\alpha \xi_\beta \big)^{-1} \sum_{\alpha, \beta} \left((\lambda+2\mu)d_2 h^{\alpha \beta}_1+  \frac{\lambda- (\lambda+2\mu) d_1 }{2} h_1g^{\alpha \beta}\right) \xi_\alpha \xi_\beta .\nonumber\end{eqnarray}
Evaluating this on all unit vectors $\xi'=(\xi_1, \cdots, \xi_{n-1})\in T^* (\partial \Omega)$ shows that $p_0$ and the values of $g_{\alpha\beta}$ along $\partial \Omega$ completely determine the matrix of quadratic form
\begin{eqnarray} \label{19.4.10-8} k^{\alpha\beta}_1=(\lambda+2\mu)d_2 h^{\alpha \beta}_1+  \frac{\lambda-(\lambda+2\mu) d_1 }{2} h_1g^{\alpha \beta}\end{eqnarray}
along $\partial \Omega$.
From (\ref{19.4.10-8}) we have \begin{eqnarray*}\sum_{\alpha, \beta} k^{\alpha\beta}_1 g_{\alpha \beta}&=& (\lambda+2\mu)d_2\sum_{\alpha, \beta} h^{\alpha \beta}_1g_{\alpha \beta}+  \frac{ \lambda-(\lambda+2\mu) d_1}{2} h_1\sum_{\alpha, \beta}g^{\alpha \beta}g_{\alpha \beta}\\
&=& h_1 \left((\lambda+2\mu)d_2+ (n-1) \frac{\lambda-(\lambda+2\mu) d_1 }{2}\right).\end{eqnarray*}
 Note that, under our assumption for any positive integer $n\ge 2$,
 \begin{align} \label{19.10.5-6}&  (\lambda+2\mu)d_2\!+\! (n\!-\!1) \big(\frac{\lambda\!-\!(\lambda\!+\!2\mu) d_1 }{2}\big)\!=\!(\lambda\!+\!2\mu)\big\{ d_2 \!+\!\frac{n\!-\!1}{2} \big(\frac{\lambda}{\lambda\!+\!2\mu}\! -\!d_1\big) \big\} \\
 &\qquad \quad \quad= \frac{(n-1) \lambda^3+(4n-2) \lambda^2\mu +(n+5)\lambda \mu^2 +(14-8n)\mu^3 }{2( \lambda +{3}\mu )^2}\ne 0.\nonumber
\end{align}
It follows that
 $$h_1= \frac{\sum_{\alpha,\beta} k^{\alpha \beta}_1 g_{\alpha\beta}}{(\lambda+2\mu)d_2+ (n-1) \big(\frac{\lambda-(\lambda+2\mu) d_1 }{2}\big)}.$$
  Thus, by (\ref{19.4.10-8}) we get \begin{eqnarray*} (\lambda\!+\!2\mu)d_2 h^{\alpha \beta}_1\!=\! k^{\alpha \beta}_1 \!+\!  \frac{(\lambda+2\mu) d_1 \!-\!\lambda}{2} \left(\frac{\sum_{\alpha,\beta}\!k^{\alpha \beta}_1 g_{\alpha\beta}}{(\lambda\!+\!2\mu)d_2\!+\!  (n\!-\!1)\big(\frac{\lambda\!-\!(\lambda\!+\!2\mu) d_1 }{2}\big)}\right)g^{\alpha \beta}.\end{eqnarray*}
From $\mu>0$ and $\lambda+\mu\ge 0$ we get  $d_2 =\frac{\mu [(\lambda+\mu)^2 + \mu(\lambda+2\mu)]} {( \lambda +{2}\mu )( \lambda +{3}\mu ) ^2}>0$.
 This shows that $p_0$ uniquely determines the $h^{\alpha \beta}_1=\frac{\partial g^{\alpha\beta}}{\partial x_n}$ along $\partial \Omega$.

In (\ref{19.3.28-15}), by putting $m=1$ we see that $p_{-1}$ uniquely determines $q_{-1}$. It follows from (\ref{19.3.28-7}) that
\begin{eqnarray} \label{19.3.29.22} (q_1- b_1) q_{-1}+q_{-1} q_1= E_0,\end{eqnarray}
so
\begin{eqnarray}\label{19.4.18-2,} \mbox{vec}\,q_{-1} = U^{-1} (\mbox{vec}\, E_0), \end{eqnarray}
 where $U^{-1}$ is given in section 3 and \begin{eqnarray} \label{19.3.29-14} E_0\!\!&\!\!=\!\!&\!\! -q_0^2 \!+\!i \sum_{l=1}^{n-1} \!\big(\frac{\partial q_1}{\partial \xi_l} \, \frac{\partial q_0}{\partial x_l}+ \frac{\partial q_0}{\partial \xi_l} \, \frac{\partial q_1}{\partial x_l}\big) \!+\!\frac{1}{2} \sum_{l,\gamma}^{n-1} \frac{\partial^2 q_1}{\partial \xi_l\partial \xi_\gamma}\, \frac{\partial^2q_1}{\partial x_l \partial x_\gamma} \\
 && +b_0q_0 -i \sum_{l=1}^{n-1} \frac{\partial b_1}{\partial \xi_l}\,\frac{\partial q_0}{\partial x_l} +\frac{\partial q_0}{\partial x_n} -c_0.\nonumber\end{eqnarray}
  (\ref{19.4.18-2,}) implies that $q_{-1}$ and $E_0$ are determined each other.
 Therefore, $p_{-1}$ uniquely determines $E_0$. Clearly, $E_0$ can be written as
 \begin{eqnarray}&& \label{19.10.3-3}\!\!\!\! E_0\!=\! \frac{\partial q_0}{\partial x_n} \!-\!{\scriptsize \left[\!\begin{BMAT}(@, 1pt, 1pt){c.c}{c.c} \left[
  \rho_{jk}\! +\! R^{j}_k\right]_{(n\!-\!1)\times (n\!-\!1)} &  \left[ \rho_{jn}\!+\!R^j_n\right]_{(n\!-\!1)\times 1}\\
  \frac{\mu}{\lambda\!+\!2\mu} (\rho_{nk}\!+\!R^n_k) & \frac{\lambda\!+\!\mu}{\lambda\!+\!2\mu}\sum\limits_{\beta} \!\frac{\partial \Gamma_{n\beta}^\beta}{\partial x_n} \!+\!
  \frac{\mu}{\lambda\!+\!2\mu} (\rho_{nn}\!+\! R^n_{n} )\end{BMAT}\!\right]}
 \!+\!T_{-1}^{(1)}(g_{\alpha\beta}), \end{eqnarray}
 where and elsewhere each $T_{-1}^{(s)}(g_{\alpha\beta})$ is a matrix expression involving only $g_{\alpha\beta}$, $g^{\alpha\beta}$, their
 normal derivatives $\frac{\partial g^{\alpha\beta}}{\partial x_n}$ and tangential derivatives along $\partial \Omega$.
  Recall that $\rho_{js}=  \sum_{k,l=1}^n \big(g^{kl} \frac{\partial \Gamma^j_{sl}}{\partial x_k} +\sum\limits_{h=1}^n g^{kl} \Gamma_{hl}^j \Gamma_{sk}^h -\sum\limits_{h=1}^n g^{kl} \Gamma_{sh}^j \Gamma_{kl}^h \big)$ for $s,l=1,\cdots,n$.
  It follows that \begin{eqnarray} \label{19.10.3-4} \rho_{nn}= \sum_{k,l=1}^n g^{kl} \frac{\partial \Gamma^n_{nl}}{\partial x_k}  +T^{(2)}_{-1}(g_{\alpha\beta})= 0+ T^{(2)}_{-1}(g_{\alpha\beta}) \end{eqnarray} because of $\Gamma_{nl}^n =0$ (see (\ref{19.10.5-1})).
 Also, from \begin{eqnarray*} R^{\,r}_j = \sum_{s,l=1}^n g^{rs}\bigg( \frac{\partial \Gamma_{js}^l}{\partial x_l} -\frac{\partial \Gamma_{ls}^l}{\partial x_j} +\sum_{m=1}^n \left(\Gamma_{js}^m \Gamma_{lm}^l -\Gamma_{ls}^m \Gamma_{jm}^l\right)\bigg), \quad 1\le j,r\le n,\end{eqnarray*}
 we see that
\begin{align} \label{19.06.25-2} R_n^{\,n}&= \sum_{l=1}^n \bigg( -\frac{\partial \Gamma_{nl}^l}{\partial x_n} -\sum_{m=1}^n \Gamma_{ln}^m \Gamma_{nm}^l\bigg)\\
&= -\frac{\partial }{\partial x_n} \bigg(\frac{1}{2}\sum\limits_{\alpha,\beta} g^{\alpha\beta} \frac{\partial g_{\alpha\beta}}{\partial x_n}\bigg) + T_{-1}^{(3)}(g_{\alpha\beta})\nonumber\\
& =-\frac{1}{2}\sum\limits_{\alpha,\beta} g^{\alpha\beta} \frac{\partial^2 g_{\alpha\beta}}{\partial x_n^2} + T_{-1}^{(4)}(g_{\alpha\beta}).\nonumber\end{align}
 In view of (\ref{19.4.10-4}), we have
 \begin{eqnarray*} \frac{\partial q_0}{\partial x_n}\! =\! \frac{\partial }{\partial x_n}\! \left\{\!\begin{small} \left[\!\begin{BMAT}(@, 0pt, 0pt){c.c}{c.c} \left[ \;\,* \;\;\right]_{(n\!-\!1)\times (n\!-\!1)}& \left[\;\;*\;\;\right]_{(n\!-\!1)\times 1}\\
\left[\;\; *\;\;\right]_{1\times (n\!-\!1)} & \frac{ d_1}{2} \sum_{\alpha,\beta}\! g^{\alpha\beta} \frac{\partial g_{\alpha\beta}}{\partial x_n} \!+\! \frac{ d_2\sum_{\alpha, \beta}\!
 \frac{\partial g^{\alpha\beta}}{\partial x_n}  \xi_\alpha \xi_\beta}{\sum_{\alpha, \beta}\! g^{\alpha \beta} \xi_\alpha \xi_\beta} \end{BMAT}\right]\end{small} \!\!+\!T_0^{(2)}(g_{\alpha\beta})\right\}\\
 \!\!\!\!=  \begin{small} \left[\!\begin{BMAT}(@, 1pt, 1pt){c.c}{c.c} \left[ \;\,* \;\;\right]_{(n\!-\!1)\times (n\!-\!1)}& \left[\;\;*\;\;\right]_{(n\!-\!1)\times 1}\\
\left[\;\; *\;\;\right]_{1\times (n\!-\!1)} & \frac{ d_1}{2} \sum_{\alpha,\beta}\! g^{\alpha\beta} \frac{\partial^2 g_{\alpha\beta}}{\partial x_n^2} \!+\! \frac{ d_2\sum_{\alpha, \beta}
 \frac{\partial^2 g^{\alpha\beta}}{\partial x_n^2}  \xi_\alpha \xi_\beta}{\sum_{\alpha, \beta} \!g^{\alpha \beta} \xi_\alpha \xi_\beta} \end{BMAT}\right]\end{small} \!+\!T_{-1}^{(5)}(g_{\alpha\beta}). \qquad\end{eqnarray*}
Combining  this, (\ref{19.10.3-3}), (\ref{19.10.3-4}) and (\ref{19.06.25-2})  we get that the $(n,n)$-entry of the matrix $E_0$ is
\begin{eqnarray}&& \!\!\!\!\!\!\!\!\!\!\! \label{19.10.3-5}
 \!\!\! \!\!\! E_0^{nn}\!:=\! \frac{d_1}{2} \!\sum\limits_{\alpha,\beta}\! g^{\alpha\beta}\! \frac{\partial^2 g_{\alpha\beta}}{\partial x_n^2}\! +\! \frac{d_2\!\sum\limits_{\alpha, \beta}
 \frac{\partial^2 g^{\alpha\beta}}{\partial x_n^2} \! \xi_\alpha \xi_\beta}{\sum\limits_{\alpha, \beta}\! g^{\alpha \beta} \xi_\alpha \xi_\beta}
 \!-\!\frac{\lambda}{2(\lambda\!+\!2\mu)}\! \sum\limits_{\alpha,\beta}\! g^{\alpha\beta} \frac{\partial^2 g_{\alpha\beta}}{\partial x_n^2}
  \!+\!T_{-\!1}^{(6)}(g_{\alpha\beta})\\
&&\!\;\; = \frac{1}{2}\!\bigg( d_1\!-\!\frac{\lambda}{\lambda+2\mu}\bigg) \sum\limits_{\alpha,\beta} \!g^{\alpha\beta}\, \frac{\partial^2 g_{\alpha\beta}}{\partial x_n^2} + \frac{ d_2\sum\limits_{\alpha, \beta}\!
 \frac{\partial^2 g^{\alpha\beta}}{\partial x_n^2}  \xi_\alpha \xi_\beta}{\sum_{\alpha, \beta}\! g^{\alpha \beta} \xi_\alpha \xi_\beta}
   \!+\!T_{-1}^{(6)}(g_{\alpha\beta}).\nonumber\end{eqnarray}
 By applying $\sum_{\alpha, \beta} g^{\alpha\beta} g_{\alpha\beta}=n-1$ again, we have
$$ \sum_{\alpha,\beta} \left(\frac{\partial^2 g^{\alpha\beta}}{\partial x_n^2} g_{\alpha\beta} + 2\frac{\partial g^{\alpha\beta}}{\partial x_n}\,\frac{\partial g_{\alpha\beta}}{\partial x_n}+g^{\alpha\beta}\,\frac{\partial^2 g_{\alpha\beta}}{\partial x_n^2} \right)=0,$$
which implies \begin{eqnarray}\label{19.10.3-8}\sum_{\alpha,\beta} g^{\alpha\beta}\,\frac{\partial^2 g_{\alpha\beta}}{\partial x_n^2}
 =- \sum_{\alpha,\beta}  g_{\alpha\beta}\frac{\partial^2 g^{\alpha\beta}}{\partial x_n^2} +T_{-1}^{(7)}(g_{\alpha\beta})
.\end{eqnarray}
Hence
\begin{align}\label{19.10.3.-9} E_0^{nn}& =\frac{1}{2}\bigg(\frac{\lambda}{\lambda+2\mu}-  d_1\bigg) \sum_{\alpha,\beta} g_{\alpha\beta}\, \frac{\partial^2 g^{\alpha\beta}}{\partial x_n^2} + \frac{ d_2\sum_{\alpha, \beta}
 \frac{\partial^2 g^{\alpha\beta}}{\partial x_n^2} \, \xi_\alpha \xi_\beta}{\sum_{\alpha, \beta} g^{\alpha \beta} \xi_\alpha \xi_\beta}
   +T_{-1}^{(8)}(g_{\alpha\beta})\\
  & = \frac{1}{\sum_{\alpha, \beta}\! g^{\alpha \beta} \xi_\alpha \xi_\beta} \sum_{\alpha, \beta} \!\bigg\{\!\frac{1}{2}\!\bigg(\!\frac{\lambda}{\lambda\!+\!2\mu}\!-\!  d_1\bigg)h_2  g^{\alpha\beta} \!+\!d_2
 \frac{\partial^2 g^{\alpha\beta}}{\partial x_n^2} \!\bigg\} \xi_\alpha \xi_\beta \!+\! T_{-1}^{(8)}(g_{\alpha\beta}),\nonumber\end{align}
where $h_{2}:= \sum_{\alpha,\beta} g_{\alpha\beta} \, \frac{\partial^2 g^{\alpha \beta}}{\partial x_n^2}$.
 Because $E_0^{nn}$ is determined by $p_{-1}$, we see that $k^{\alpha\beta}_{2}:=\frac{1}{2}\left(\frac{\lambda}{\lambda+2\mu}-  d_1\right)h_2 \, g^{\alpha\beta} +d_2
 \frac{\partial^2 g^{\alpha\beta}}{\partial x_n^2}$, for all  $\alpha, \beta$, are also determined by $p_{-1}$.
 From (\ref{19.10.5-6}) and \begin{eqnarray*} && \sum_{\alpha,\beta} k^{\alpha\beta}_{2}g_{\alpha\beta} =\frac{1}{2}\bigg(\frac{\lambda}{\lambda+2\mu}-  d_1\bigg)h_{2}\sum_{\alpha, \beta}g^{\alpha\beta}  g_{\alpha \beta} +d_2\sum_{\alpha,\beta}  \frac{\partial^2 g^{\alpha\beta}}{\partial x_n^2}g_{\alpha\beta} \\
 &&\qquad \qquad \quad =\bigg\{ \frac{1}{2}(n-1) \bigg(\frac{\lambda}{\lambda+2\mu}-  d_1\bigg)+d_2\bigg\} h_2,\end{eqnarray*}
 we get \begin{eqnarray*}h_{2}= \frac{\sum_{\alpha,\beta} k^{\alpha\beta}_2g_{\alpha\beta}}{\frac{1}{2}(n-1) \big(\frac{\lambda}{\lambda+2\mu}-  d_1\big)+d_2}  \end{eqnarray*}  for $(n-1) \lambda^3+(4n-2) \lambda^2\mu +(n+5)\lambda \mu^2 +(14-8n)\mu^3 \neq 0$.
Thus \begin{eqnarray*} \frac{\partial^2  g^{\alpha\beta}}{\partial x_n^2} = \frac{1}{d_2} \bigg\{ k^{\alpha\beta}_{2} -
\frac{1}{2}\left(\frac{\lambda}{\lambda+2\mu}-  d_1\right)h_2 \, g^{\alpha\beta}\bigg\}
.\end{eqnarray*}
Since the right-hand side of the above equality is determined by $p_{-1}$, we get that $\frac{\partial^2 g^{\alpha\beta}}{\partial x_n^2}$,  $1\le \alpha, \beta\le n-1$, can be uniquely determined by $p_{-1}$ if $\mu>0$, $\lambda+\mu\ge 0$ and $(n-1) \lambda^3+(4n-2) \lambda^2\mu +(n+5)\lambda \mu^2 +(14-8n)\mu^3 \neq 0$.

 Finally, from the symbol equation we have
\begin{eqnarray*} p_{-m-1}=    \left[\begin{BMAT}(@, 1pt, 1pt){c.c}{c.c} \left[ \mu \delta_{jk}\right]_{(n-1)\times(n-1)} & \left[0\right]_{(n-1)\times 1}\\
\left[0\right]_{1\times (n-1)} & \lambda+2\mu\end{BMAT}\right] q_{-m-1},\quad m\ge 1, \end{eqnarray*}
so that $p_{-m-1}$ uniquely determines $q_{-m-1}$. Since $q_{-m-1}$ satisfies (\ref{19.3.28-7}) we have $\mbox{vec}\, E_{-m}=U\big(\mbox{vec}\, q_{-m-1}\big)$. This implies that $q_{-m-1}$
uniquely determines $\mbox{vec}\, E_{-m}$, and hence $p_{-m-1}$ uniquely determines $E_{-m}$.
Since $q_1,q_0,\cdots, q_{-m}$ have been recovered by $p_1, p_0, \cdots, p_{-m}$ from the previous discussion, if follows from (\ref{19.6.1-2}) that $E_{-m}$ can be written as \begin{eqnarray} \label{19.4.13;5}E_{-m}= \frac{\partial q_{-m}}{\partial x_n}+ T_{-m-1}^{(1)}(g_{\alpha \beta}),\end{eqnarray}
where and elsewhere $T_{-m-1}^{(s)}(g_{\alpha\beta})$ involves only the boundary values of $g_{\alpha\beta}$, $g^{\alpha\beta}$, and their normal derivatives of order at most $|-m-1|$. According to  (\ref{19.10.6-1}) we get that $p_{-m-1}$ uniquely determines
$U (\mbox{vec}\, E_{-m})$, so $p_{-m-1}$ uniquely determines $U \big(\mbox{vec}\, \frac{\partial q_{-m-1}}{\partial x_n} \big)+ U \big(\mbox{vec}\, T^{(1)}_{-m-1}(g_{\alpha\beta})\big)$,  where $U$ is given by (\ref{19.3.23-15}) and determined by $p_1$, and $\mbox{vec}\, X$ is a stack of all columns of matrix $X$.
Since $U(\mbox{vec}\,q_{-m}\big) =\mbox{vec}\, E_{-m+1}$, we have  \begin{align}  &\label{19.10.7-1} \; U \big(\mbox{vec}\, \frac{\partial q_{-m}}{\partial x_n} \big)+ U \big(\mbox{vec}\, T^{(1)}_{-m-1}(g_{\alpha\beta})\big)\\ & \qquad =
  \frac{\partial }{\partial x_n} \big(U \big(\mbox{vec}\,q_{-m}\big)\big) -  \frac{\partial U}{\partial x_n} \,\big(\mbox{vec}\,q_{-m}\big)\big) +  U \big(\mbox{vec}\, T^{(1)}_{-m-1}(g_{\alpha\beta})\big)\nonumber\\
 &\qquad   =  \frac{\partial  }{\partial x_n}\big(\mbox{vec}\, E_{-m+1}\big) +T^{(2)}_{-m-1}(g_{\alpha\beta}).\nonumber \end{align}
  Proceeding now by induction, let $m\ge 1$, and suppose we have show that, when $-m+1\le l+1\le 0$,
the $(n,n)$-entry of $E_{l+1}$ (i.e., the $n^2$-th entry of $\mbox{vec}\, E_{l+1}$)
 is \begin{eqnarray}&&\label{19.4.18-3..}
  \bigg(\!\sqrt{\sum\limits_{\alpha,\beta}\! g^{\alpha\beta}\xi_\alpha\xi_\beta}\bigg)^{\!l\!+\!1}\!\bigg(\! \frac{1}{2}\big( d_1\!-\!\frac{\lambda}{\lambda\!+\!2\mu}\big)\! \sum_{\alpha,\beta} \!g^{\alpha\beta} \frac{\partial^{|l|+1} g_{\alpha\beta}}{\partial x_n^{|l|+1}}\!\\
  && \quad \quad + \!\frac{ d_2\sum_{\alpha, \beta}
 \frac{\partial^{|l|+1} g^{\alpha\beta}}{\partial x_n^{|l|+1} } \xi_\alpha \xi_\beta}{\sum_{\alpha, \beta}\! g^{\alpha \beta} \xi_\alpha \xi_\beta}\!\bigg)
  \! +\!T_{l-1}^{(3)}(g_{\alpha\beta}),\nonumber \end{eqnarray}
where $T_{l-1}^{(3)}(g_{\alpha\beta})$ involves only the boundary values of $g_{\alpha\beta}$, $g^{\alpha\beta}$, and their normal derivatives of order at most $|l|$ (Note that (\ref{19.4.18-3..}) holds for $l+1=0$, see (\ref{19.10.3-5})).
Combining above discussion and (\ref{19.10.7-1}), we get that $p_{-m-1}$ uniquely determines the $n^2$-th entry of $\frac{\partial}{\partial x_n} (\mbox{vec}\,  E_{-m+1})$. However by (\ref{19.4.18-3..}),
we see that the the $n^2$-th entry of $\frac{\partial}{\partial x_n} (\mbox{vec}\,  E_{-m+1})$ equals to \begin{eqnarray*}\label{19.5.29-1}&& \frac{\partial}{\partial x_n} \bigg\{\bigg(\sqrt{\sum\limits_{\alpha,\beta} g^{\alpha\beta}\xi_\alpha\xi_\beta}\bigg)^{-m+1}\bigg( \frac{1}{2}\big( d_1-\frac{\lambda}{\lambda+2\mu}\big) \sum_{\alpha,\beta} g^{\alpha\beta}\, \frac{\partial^{m+1} g_{\alpha\beta}}{\partial x_n^{m+1}} \\
&&  \quad \,+ \frac{ d_2\sum_{\alpha, \beta}
 \frac{\partial^{m+1} g^{\alpha\beta}}{\partial x_n^{m+1} }\, \xi_\alpha \xi_\beta}{\sum_{\alpha, \beta} g^{\alpha \beta} \xi_\alpha \xi_\beta}\bigg)\bigg\}
   +T_{-m-1}^{(4)}(g_{\alpha\beta})\\
 &&  = \bigg(\sqrt{\sum\limits_{\alpha,\beta} g^{\alpha\beta}\xi_\alpha\xi_\beta}\bigg)^{-m+1}\bigg( \frac{1}{2}\big( d_1-\frac{\lambda}{\lambda+2\mu}\big) \sum_{\alpha,\beta} g^{\alpha\beta}\, \frac{\partial^{m+2} g_{\alpha\beta}}{\partial x_n^{m+2}} \\
 &&\quad\,  + \frac{ d_2\sum_{\alpha, \beta}
 \frac{\partial^{m+2} g^{\alpha\beta}}{\partial x_n^{m+2} }\, \xi_\alpha \xi_\beta}{\sum_{\alpha, \beta} g^{\alpha \beta} \xi_\alpha \xi_\beta}\bigg)
   +T_{-m-1}^{(5)}(g_{\alpha\beta})\quad\;  \; \\
  && = \bigg(\sqrt{\sum\limits_{\alpha,\beta} g^{\alpha\beta}\xi_\alpha\xi_\beta}\bigg)^{-m+1}\bigg( \frac{1}{2}\big(- d_1+\frac{\lambda}{\lambda+2\mu}\big) \sum_{\alpha,\beta} g_{\alpha\beta}\, \frac{\partial^{m+2} g^{\alpha\beta}}{\partial x_n^{m+2}} \\
  && \quad\, + \frac{ d_2\sum_{\alpha, \beta}
 \frac{\partial^{m+2} g^{\alpha\beta}}{\partial x_n^{m+2} }\, \xi_\alpha \xi_\beta}{\sum_{\alpha, \beta} g^{\alpha \beta} \xi_\alpha \xi_\beta}\bigg)
   +T_{-m-1}^{(6)}(g_{\alpha\beta}). \end{eqnarray*}
 Thus, $p_{-m-1}$ uniquely determines the following $k_{m+2}^{\alpha \beta}$:
 $$k^{\alpha\beta}_{m+2}:= \frac{1}{2}\big(- d_1+\frac{\lambda}{\lambda+2\mu}\big) \big(\sum_{\gamma,\sigma} g_{\gamma\sigma}\, \frac{\partial^{m+2} g^{\gamma\sigma}}{\partial x_n^{m+2}}\big)g^{\alpha\beta} +  d_2
 \frac{\partial^{m+2} g^{\alpha\beta}}{\partial x_n^{m+2} }.$$
Since \begin{align*} \sum\limits_{\alpha,\beta} \!k^{\alpha\beta}_{m+2}\,g_{\alpha\beta}:=\!&\, \frac{1}{2}\!\big(\!\!- \!d_1\!+\!\frac{\lambda}{\lambda\!+\!2\mu}\big) \!\big(\!\sum\limits_{\gamma,\sigma}\! g_{\gamma\sigma} \frac{\partial^{m\!+\!2} \!g^{\gamma\sigma}}{\partial x_n^{m\!+\!2}}\!\big)\!\!\sum\limits_{\alpha,\beta}\!g^{\alpha\beta}g_{\alpha\beta} \!+\!  d_2\!
\sum\limits_{\alpha, \beta}g_{\alpha\beta}\! \frac{\partial^{m\!+\!2} g^{\alpha\beta}}{\partial x_n^{m+2} }\\
\,= &\,\bigg(\frac{n-1}{2}\big(-\! d_1\!+\!\frac{\lambda}{\lambda+2\mu}\big)  \!+\!  d_2\bigg) \sum\limits_{\alpha,\beta} g_{\alpha\beta} \frac{\partial^{m+2} g^{\alpha\beta}}{\partial x_n^{m+2}},\end{align*}
we get \begin{eqnarray*} \sum_{\alpha,\beta} g_{\alpha\beta}\, \frac{\partial^{m+2} g^{\alpha\beta}}{\partial x_n^{m+2}}= \frac{\sum\limits_{\alpha,\beta} k^{\alpha\beta}_{m+2}g_{\alpha\beta}}{\frac{n-1}{2}\big(- d_1+\frac{\lambda}{\lambda+2\mu}\big)  +  d_2},\end{eqnarray*}
so that \begin{eqnarray}&&\label{19.10.7-6}\frac{\partial^{m+2} g^{\alpha\beta}}{\partial x_n^{m+2} }\!=\! \frac{1}{d_2}\! \bigg\{k^{\alpha\beta}_{m+2}\!-\! \frac{1}{2}\!\big(\!-\! d_1\!+\!\frac{\lambda}{\lambda\!+\!2\mu}\big)\bigg(\!\frac{\sum\limits_{\alpha,\beta} \!k^{\alpha\beta}_{m+2}g_{\alpha\beta}}{\frac{n\!-\!1}{2}\big(\!-\! d_1\!+\!\frac{\lambda}{\lambda\!+\!2\mu}\big) \! +\!  d_2}\bigg)
g^{\alpha\beta}\!\bigg\}. \end{eqnarray}
This shows that $p_{-m-1}$ (together with $p_1, \cdots, p_{-m}$) determines the $(m+2)$-order normal derivatives of $g^{\alpha\beta}$ if $\mu>0$, $\lambda+\mu\ge 0$ and $(n-1) \lambda^3+(4n-2) \lambda^2\mu +(n+5)\lambda \mu^2 +(14-8n)\mu^3 \neq 0$. Therefore, the first $(m+2)$-order normal derivatives of $g^{\alpha \beta}$ can be determined by $p_1, \cdots, p_{-m-1}$ if $\mu>0$, $\lambda+\mu\ge 0$ and $(n-1) \lambda^3+(4n-2) \lambda^2\mu +(n+5)\lambda \mu^2 +(14-8n)\mu^3 \neq 0$.
  This completes the proof of Proposition 4.1 by induction.
\qed

\vskip 0.39 true cm

 R. Kohn  and M. Vogelius first showed that Dirichlet-to-Neumann map associated with the equation $-\mbox{div}\; (a\nabla v)=0$ uniquely determines $a$ and its all order derivatives on $\partial \Omega$ in \cite{KV1,KV2,KV3}, and J. Lee and G. Uhlmann showed that for the equation $\frac{1}{\sqrt{|g|}}\sum_{j,k=1}^n \frac{\partial}{\partial x_j}(\sqrt{|g|} g^{jk} \frac{\partial v}{\partial x_k})=0$, the Dirichlet-to-Neumann map uniquely determines $\{g_{jk}(x)\big|_{\partial \Omega}\}$, which are a very important results for global uniqueness theorems for various inverse boundary value problems (see \cite{LU}, \cite{SU2}, \cite{SU3}). Our new method has used the exact complicated calculations for
 the full symbol of the pseudodifferential operators $\Xi_g$, and then get our result by discussing the $(n,n)$-entry of $p_{1-l}(x,\xi)$ for $l=0,1,2,3,\cdots$, where  $p_{1-l}(x,\xi)$ is homogeneous in $\xi$ of degree $1-l$ for $|\xi|\ge 1$ and $\sigma(\Xi_g)\sim
 \sum_{l\ge 0}p_{1-l}(x,\xi)$ (see also \cite{Liu1} and \cite{Liu2} for more details).

\vskip 0.39 true cm

 \noindent{\bf Lemma 4.2.} \ {\it Let $\bar \Omega$ be a compact real-analytic $n$-dimensional Riemannian manifold with real-analytic boundary, $n\ge 2$, and let $g$ and $\tilde{g}$ be real analytic metrics on $\bar \Omega$. Assume that the Lam\'{e} constants $\mu>0$ and $\lambda$ satisfy $\lambda+\mu> 0$ and $(n-1) \lambda^3+(4n-2) \lambda^2\mu +(n+5)\lambda \mu^2 +(14-8n)\mu^3 \neq 0$. If $\Xi_g=\Xi_{\tilde g}$, then there exists a neighborhood $\mathfrak{U}$ of $\partial \Omega$ in $\bar \Omega$ and a real-analytic map $\varrho_0: \mathfrak{U}\to \bar \Omega$ such that $\varrho_0\big|_{\partial \Omega}=\mbox{identity}$ and $g=\varrho_0^* \tilde {g}$.}

 \vskip 0.23 true cm

    \noindent  {\it Proof.} \  Let $\mathfrak{V}$ be some connected open set in the half-space $\{ x_n\ge 0\}\subset {\Bbb R}^n$ containing the origin. For any $x_0\in \partial \Omega$ we define a real-analytic diffeomorphism $\varsigma_{x_0}: \mathfrak{V}\to {\mathfrak{U}}_{x_0}$ (respectively, ${\tilde\varsigma}_{x_0}: \mathfrak{V}\to {\tilde{\mathfrak{U}}}_{x_0}$), where  $\{ x_1, \cdots, x_{n}\}$ (respectively, $\{ {\tilde x}_1, \cdots, {\tilde x}_{n}\}$) denote the corresponding boundary normal coordinates for $g$ (respectively, $\tilde g$) defined on a connected neighborhood ${\mathfrak{U}}_{x_0}$ (respectively, ${\tilde{\mathfrak{U}}}_{x_0}$) of $x_0$ (see \cite{LU}).     It follows from Proposition 4.1 that the metric $\varsigma_{x_0}^* g$ and ${\tilde{\varsigma}}_{x_0}^* \tilde{g}$ are real-analytic metric on $\mathfrak{V}$, where Taylor series at origin by explicit formulas involving the symbol of $\Xi_g$ in $\{x_\alpha\}$ coordinates. Clearly, these two metric must be identical on $\mathfrak{V}$. Set $ \varrho_{x_0} ={\tilde{\varsigma}}_{x_0} \circ \varsigma_{x_0}^{-1} :{\mathfrak{U}}_{x_0}\to {\tilde{\mathfrak{U}}}_{x_0}$.  Then, $\varrho_{x_0}$
is a real-analytic diffeomorphism which fixes the portion of $\partial \Omega$ lying in ${\mathfrak{U}}_{x_0}$ and satisfies ${\varrho}^*_{x_0}\tilde {g}= g$.
Therefore, by such a construction for each $x_0\in \partial \Omega$, we get a real analytic map $\varrho_0: \mathfrak{U}\to \bar \Omega$ such that $\varrho_0\big|_{\partial \Omega}=\mbox{identity}$ and $\varrho_0^*\tilde{g} =g$. \ \ \  \  \qed

\vskip 0.2 true cm

 A Riemannian manifold-with-boundary $(\bar \Omega, g)$ is strongly convex if, for every pair of points $p, q\in \bar \Omega$, there is a unique minimal geodesic $\gamma$ jointing $p$ and $q$ whose interior lies in $\Omega$.
The following Lemma, which is due to J. Lee and G. Uhlmann \cite{LU} and S. D. Myers\cite{Myers}, will be used late. It provides a possibility of extending a local isometry to a global one:

\vskip 0.19 true cm

 \noindent{\bf Lemma 4.3 (Lee-Uhlmann-Myers}, see p.$\,$1107 of \cite{LU} or \cite{Myers}). \ {\it Let $\bar \Omega$ be a compact, connected, real-analytic $n$-manifold with real-analytic boundary, $n\ge 2$, and assume that $\pi(\bar \Omega, \partial \Omega)$ is trivial. Let $g$ and $\tilde{g}$ be real analytic metric on $\bar \Omega$, and suppose on some neighborhood $\mathfrak{U}$ of $\partial \Omega$
in $\bar \Omega$ we are given a real-analytic map $\varrho_0: \mathfrak{U}\to \bar \Omega$ such that $\varrho_0\big|_{\partial \Omega}=\,$identity and $\varrho_0^* \tilde g=g$.
 Assume that one of the following conditions holds:

 (a) \ \ $\bar \Omega$ is strongly convex with respect to both $g$ and $\tilde {g}$;

\noindent or

 (b) \ \  either $g$ or $\tilde{g}$ extends to a complete real-analytic metric on a non-compact real-analytic manifold $\tilde{\Omega}$ without boundary containing $\bar \Omega$.

 Then $\varrho_0$ extends to a real-analytic diffeomorphism $\varrho:\bar \Omega\to \bar \Omega$ such that $g=\varrho^* \tilde {g}$.}

\vskip 0.32 true cm

  \noindent  {\it Proof of Theorem 1.1.} \  From Proposition 4.1 and Lemma 4.2, we see that there exists a neighborhood ${\mathfrak{U}}$ of $\partial \Omega$ in $\bar \Omega$ and a real analytic map $\rho_0: {\mathfrak{U}}\to \bar \Omega$ such that $\rho_0\big|_{\partial \Omega}=\mbox{identity}$ and $g=\rho_0^* {\tilde g}$. It follows from Lee-Uhlmann-Mayers' result (i.e., Lemma 4.3) that the map $\rho_0$ can be extends to a real analytic diffeomorphism $\rho: \bar \Omega\to \bar \Omega$ such that $g=\rho^* {\tilde g}$.     \qed

\vskip 0.22 true cm

\noindent{\bf Remark 4.4.} \   i) \  The conclusion of Theorem 1.1 holds for any $n\ge 2$ if $(\mu, \lambda)\in \{(\mu,\lambda)\in {\mathbb{R}}^2\,\big|\, \mu>0,\, \lambda+\mu> 0,\; (n-1) \lambda^3+(4n-2) \lambda^2\mu +(n+5)\lambda \mu^2 +(14-8n)\mu^3 \ne 0\}$. But for the Dirichlet-to-Neumann map associated with Laplacian, the corresponding result holds only for $n\ge 3$ because  in two-dimensional case, Lee-Uhlmann's result says that $\varrho^*\tilde{g}$ is a conformal multiple of $g$ (see \cite{LU}).

ii) \  In particular, if $\Omega$ is simply connected in Theorem 5.3, then $\pi_1(\bar \Omega, \partial \Omega)=0$ (see p.$\,$162 of \cite{Whi}).

\vskip 0.22 true cm

\noindent{\bf Remark 4.5.} \ Using our new method, we can also prove that Theorem 1.1 and Corollary 1.2 hold when the Lam\'{e} constants $\mu>0$ and $\lambda$ become smooth functions on $\Omega$ with $\mu>0$ and $\lambda+\mu> 0$.

\vskip 1.29 true cm

\section{Spectral invariants of elastic Dirichlet-to-Neumann map}

\vskip 0.45 true cm

In this section, we will compute the heat (spectral) invariants of the elastic Dirichlet-to-Neumann map. Because we have explicitly obtained the full symbol of the Dirichlet-to-Neumann map $\Xi_g$ associated with elastic Lam\'{e} operator, we
can now calculate the local expressions for the first $n-1$ heat invariants of $\Xi_g$ by some new techniques together with symbol calculus.

\vskip 0.18 true cm

 We need the following lemma late:

\vskip 0.2 true cm

\noindent{\bf Lemma 5.1.} \ {\it For any constants $\tilde{a}\ne 0$ and $\tilde{b}$, we have
\begin{eqnarray}\label{19.4.14;3}\;\;\; \left[\tilde{a}\delta_{jk}\! +\!\tilde{b} \xi_{j}\xi_k\right]_{(n\!-\!1)\times (n\!-\!1)}^{-1} \!=\! \left[ \frac{\delta_{jk}}{\tilde{a}}\!-\! \frac{\tilde{b} }{\tilde{a}(\tilde{a}\!+\!\tilde{b}|\xi'|^2)}\xi_j\xi_k\right]_{(n\!-\!1)\times(n\!-\!1)},\end{eqnarray}
where $\xi'=(\xi_1, \cdots ,\xi_{n-1})\in {\Bbb R}^{n-1}$ and $|\xi'|=\sqrt{\sum_{l=1}^{n-1}\xi_l^2}$.}
\vskip 0.25 true cm

\noindent  {\it Proof.} \ Note that \begin{eqnarray*}
\left[ \xi_j\xi_k\right]_{(n-1)\times (n-1)} \left[ \xi_j\xi_k\right]_{(n-1)\times (n-1)} =|\xi'|^2\left[ \xi_j\xi_k\right]_{(n-1)\times (n-1)}. \end{eqnarray*} This implies that the inverse of the matrix $\left[\tilde{a}\delta_{jk} +\tilde{b} \xi_{jk}\right]_{(n-1)\times (n-1)}$ has form
\begin{eqnarray*} \left[\tilde{a}\delta_{jk} +\tilde{b} \xi_{j}\xi_k\right]_{(n-1)\times (n-1)}^{-1}= \left[c\delta_{jk} +d\xi_j\xi_k\right]_{(n-1)\times (n-1)},\end{eqnarray*}
where $c$ and $d$ are two undetermined constants.
From \begin{eqnarray*} I_{n-1}&=& \left[\tilde{a}\delta_{jk} +\tilde{b} \xi_{jk}\right]_{(n-1)\times (n-1)}\left[c\delta_{jk} +d\xi_j\xi_k\right]_{(n-1)\times (n-1)}\\
&=&\tilde{a}c I_{n-1} +( \tilde{a} d+\tilde{b} c) \left[ \xi_j\xi_k\right]_{(n-1)\times (n-1)} +\tilde{b} d |\xi'|^2 \left[ \xi_j\xi_k\right]_{(n-1)\times (n-1)},\end{eqnarray*}
we get  $$\tilde{a}c=1, \; \, \tilde{a} d+ \tilde{b} c +\tilde{b} d |\xi'|^2 =0,$$ so that
\begin{eqnarray*}c=\frac{1}{\tilde{a}},  \quad d= -{\tilde{b}}c ({\tilde{a}} +\tilde{b} |\xi'|^2)^{-1}.\end{eqnarray*}
This yields (\ref{19.4.14;3}).  \qed

\vskip 0.39 true cm

\noindent  {\it Proof of Theorem 1.3.} \ \ \  Step 1. \  Let $\Psi(\tau)$ be a two-sided parametrix for $\Xi_g -\tau I$, i.e., $\Psi(\tau)$ is a pseudodifferential operator of order $-1$ with parameter $\tau$ for which
 \begin{eqnarray} \label{9.14-1}  (\Xi_g-\tau I)\,\Psi (\tau) = I \quad\mbox{mod}\;\; OPS^{-\infty},\\
 \Psi (\tau)\, (\Xi_g-\tau I)= I \quad \mbox{mod} \;\; OPS^{-\infty},\nonumber\end{eqnarray}
 where $I$ denotes the identity operator. Because $\Xi_g$ is an elliptic pseudodifferential operator of order $1$, it follows (see \cite{Ta2}) that such an $\Psi(\tau)$ is also a pseudodifferential operator of order $-1$ with symbol \begin{eqnarray*} \psi(x',\xi', \tau)\sim \psi_{-1}(x', \xi', \tau) +\psi_{-2}(x', \xi', \tau) +\cdots + \psi_{-1-m}(x', \xi', \tau) +\cdots.
  \end{eqnarray*}
   Let the complex parameter $\tau$ have homogeneity $1$, and let the $\psi_{-1-m}(x', \xi', \tau)$ be homogeneous of order $-1-m$ in the variables $(\xi', \tau)$. This infinite
   sum defines $\psi$ asymptotically. The symbol of the composition of the operator defined by $\Psi$ is given by (see \cite{Gr}, \cite{Ta2}, \cite{Tre} or p.$\,$71 of \cite{Ho3})
   \begin{eqnarray*}  \iota( (\Xi_g-\tau I)\Psi) \sim \sum\limits_{\vartheta} \frac{(-i)^{|\vartheta|}}{\vartheta!} \big(\partial_{\xi'}^\vartheta (\iota(\Xi_g-\tau I))\big)\big(\partial^\vartheta_{x'} \psi\big),  \end{eqnarray*}
   where $\vartheta=(\vartheta_1, \cdots, \vartheta_n)$ is a multi-index, $\iota(A)$ is the full symbol of a pseudodifferential operator $A$.
   Set $\varsigma_1=p_1 (x', \xi') -\tau I$, $\,\varsigma_0= p_0(x',\xi'), \cdots$, $\varsigma_k(x', \xi')=p_{k}(x',\xi'), \cdots$. Decompose this sum into order of homogeneity:
   \begin{eqnarray} \label{9.14-2} \sigma\big((\Xi_g-\tau I)\Psi\big)\sim \sum\limits_{\substack{-m=k-|\vartheta| -1-l\\0\le l\le m, \, -1 -m\le -1, \, k\le 1 }} \frac{(-i)^{|\vartheta|}}{\vartheta!}
   (\partial^\vartheta_{\xi'} \varsigma_k) (\partial_{x'}^\vartheta \psi_{-1-l}\big).\end{eqnarray}
  From (\ref{9.14-1}) we get the equations:
\begin{eqnarray*} I\!&=\!&\sum\limits_{\substack{-m=k-|\vartheta| -1-l\\0\le l=m,\, -1 -m\le -1, \, k\le 1}}\! \frac{(-i)^{|\vartheta|}}{\vartheta!}
   (\partial^\vartheta_{\xi'} \varsigma_k) (\partial_{x'}^\vartheta \psi_{-1-l}\big) \\
   \!&=\!& (p_1(x', \xi', \tau) -\tau I) \,\psi_{-1}(x',\xi'),\\
   0 \!&=\!&\sum\limits_{ \substack{-m=k-|\vartheta| -1-l\\0\le l < m,\, -1 -m< -1, \, k\le 1 }} \frac{(-i)^{|\vartheta|}}{\vartheta!}
   (\partial^\vartheta_{\xi'} \varsigma_k) (\partial_{x'}^\vartheta \psi_{-1-l}\big)\nonumber\\
   &=& (p_1(x',\xi')-\tau I) \,\psi_{-1-m}(x',\xi',\tau) \\
   && + \sum\limits_{ \substack{-m=k-|\vartheta| -1-l\\0\le l<m,\, m\ge 1, \, k\le 1 }}
 \frac{(-i)^{|\vartheta|}}{\vartheta!} \, \big(\partial_{\xi'}^\vartheta \varsigma_k\big) \big(\partial_{x'}^\vartheta \psi_{-1-l}(x',\xi',\tau)\big).\nonumber\end{eqnarray*}
  These equations define the $\psi_{-1-m}$ inductively, as follows:
 \begin{eqnarray} \label{19.6.16-1}
  \!\!\!\!&&\! \!\!\!\! \!\psi_{-1} (x', \xi, \tau) =( p_1(x',\xi')-\tau I)^{-1}, \end{eqnarray}
  \begin{eqnarray} \label{19.6.16-2}
  \!\!\!\!\! && \!\!\! \!\!\!\psi_{\!-\!1\!-\!m} (x'\!, \xi'\!, \tau)\!=\!- \!(p_1(x'\!, \xi')\! -\!\tau I)^{-1}\! \!\!\sum\limits_{\substack{-m=k\!-\!|\vartheta| \!-\!1\!-\!l\\0\le l<m, m\ge 1, k\le 1}}\!\!\!
 \frac{(\!-\!i)^{|\vartheta|}}{\vartheta!} \big(\partial^\vartheta_{\xi'} \varsigma_k\big) \big( \partial^\vartheta_{x'} \psi_{-\!1\!-\!l}\big),  \;\, m\ge 1. \end{eqnarray}
 For example, we can write out the first three terms for $\psi_{-1}$, $\psi_{-2}$ and $\psi_{-3}$.:
 \begin{eqnarray} \label{18/11/26--1} &&\psi_{-1} (x', \xi', \tau)= (p_1-\tau I)^{-1}, \\
 \label{18/11/26--2} &&\psi_{-2} (x', \xi', \tau) =-(p_1-\tau I)^{-1} \bigg(p_0 \psi_{-1} -i \sum\limits_{\alpha}\frac{\partial p_1}{\partial {\xi}_\alpha} \, \frac{\partial \psi_{-1}}{\partial {x}_\alpha} \bigg); \\
\label{18/11/26--3} && \psi_{-3} (x', \xi', \tau)= -(p_1 -\tau I)^{-1} \bigg[ p_0 \psi_{-2} +p_{-1} \psi_{-1} \\ &&\; \quad \qquad  -i  \sum\limits_{\gamma} \bigg(\frac{\partial p_1}{\partial {\xi}_\gamma} \, \frac{\partial \psi_{-2}}{\partial {x}_\gamma}
 \!+\!\frac{\partial p_0}{\partial {\xi}_\gamma}  \frac{\partial \psi_{-1}}{\partial {x}_\gamma} \bigg)
  \!-\! \frac{1}{2} \sum\limits_{\alpha,\beta}\!\frac{\partial^2 p_1}{\partial \xi_\alpha \xi_\beta}\, \frac{\partial^2 \psi_{-1}}{\partial x_\alpha \partial x_\beta}\bigg].\nonumber\end{eqnarray}
From the theory of elliptic equations (see  \cite{Mo1}, \cite{Mo2}, \cite{Mo3}, \cite{Pa}, \cite{St}), we see that the elastic Dirichlet-to-Neumann map $\Xi_g$ associated with Lam\'{e} operator can generate a strongly continuous semigroup $e^{-t\Xi_g}$ in a suitable space of vector-valued functions defined on $\partial \Omega$. Furthermore,
   there exists a matrix-valued function ${K} (t, x', y')$, which is called the parabolic kernel, such that (see p.$\,$4 of \cite{Frie})
                \begin{eqnarray*}  e^{-t\Xi_g}\,{w}_0(x')=\int_{\partial \Omega} {K}(t, x',y') {w}_0(y')dS(y'), \quad \,
        {{w}}_0\in (H^1(\partial \Omega))^n, \end{eqnarray*}
where $(H^1(\partial \Omega))^n= H^1(\partial \Omega) \times \cdots \times H^1(\partial \Omega)$.  Let $\{{{u}}_k\}_{k=1}^\infty$ be orthnormal eigenvectors of the elastic Dirichlet-to-Neumann map $\Xi_g$ corresponding to eigenvalues $\{\tau_k\}_{k=1}^\infty$, then the parabolic kernel  ${{K}}(t, x', y')=e^{-t \Xi_g} \delta(x'-y')$ is given by \begin{eqnarray} \label{18/12/18} {{K}}(t,x',y') =\sum_{k=1}^\infty e^{-t \tau_k} {{u}}_k(x')\otimes {{u}}_k(y').\end{eqnarray}
This implies that the integral of the trace of ${{K}}(t,x',y')$ is actually a spectral invariants:
\begin{eqnarray} \label{1-0a-2}\int_{\partial\Omega} \mbox{Tr}({{K}}(t,x',x')) dS(x')=\sum_{k=1}^\infty e^{-t \tau_k}.\end{eqnarray}
On the other hand,  the strongly continuous semigroup $e^{-t\Xi_g}$ can also be represented as  \begin{eqnarray*} e^{-t\Xi_g} =\frac{i}{2\pi } \int_{\mathcal{C}} e^{-t\lambda} (\Xi_g -\tau I)^{-1} d\tau,\end{eqnarray*}
where $\mathcal{C}$ is a suitable curve in the complex plane in the positive direction around the spectrum of $\Xi_g$ (i.e., a contour around the positive real axis).
 It follows that \begin{eqnarray}\label{18/12/2}&& {{K}}  (t,x',y') \!=\! e^{-t\Xi_g}\big(\delta(x'-y')I_{n-1}\big)  \\
&\!\!=\!&\!\!\! \frac{1}{(2\pi)^{n-1}} \!\int_{T^*(\partial \Omega)} \!\!e^{i\langle x'-y', \xi'\rangle}\! \bigg\{\!\frac{i}{2\pi}\! \int_{\mathcal{C}} \!e^{-t\tau}\! \big(\iota((\Xi_g\!-\!\tau I)^{-1})\big) d\tau\!\bigg\} d\xi'\nonumber\\
 &\!\!=\!&\!\! \! \frac{1}{(2\pi)^{n-1}} \int_{{\Bbb R}^{n-1}} e^{i\langle x'-y', \xi'\rangle} \bigg\{\frac{i}{2\pi} \int_{\mathcal{C}} e^{-t\tau} \begin{small} \bigg(\end{small}
 \psi_{-1} (x',\xi',\tau) \nonumber\\
 \!\!&&\! \!+\psi_{-2} (x',\xi',\tau) +\psi_{-3} (x',\xi',\tau) +\cdots \begin{small}\bigg) \end{small} d\tau\bigg\} d\xi',\nonumber\end{eqnarray}
 where $T^*(\partial \Omega)$ is the cotangent space at $x'$, so that \begin{eqnarray} \label{18/12/2-2}  && \mbox{Tr}\,({{K}}  (t,x',x')) \!\!=\!\! \frac{1}{(2\pi)^{n-1}}\!\int_{{\Bbb R}^{n-1}} \!\bigg\{\!\frac{i}{2\pi}\! \int_{\mathcal{C}}\! e^{-t\tau}\!  \bigg(\!
 \mbox{Tr}\,\big(\psi_{-1} (x',\xi',\tau)\big)\\
 && \quad \;\;+ \mbox{Tr}\,\big(\psi_{-2} (x',\xi',\tau)\big)\! +\! \mbox{Tr}\,\big(\psi_{-3} (x',\xi',\tau)\big) \!+\!\cdots \bigg) d\tau\!\bigg\} d\xi'.\nonumber\end{eqnarray}
  Thus \begin{eqnarray} && \sum_{k=1}^\infty \!\!e^{-t \tau_k} \!\!=\!\!\int_{\partial \Omega} \! \frac{1}{(2\pi)^{n\!-\!1}}\! \int_{{\Bbb R}^{n\!-\!1}}\! \bigg\{\!\frac{i}{2\pi} \!\int_{\mathcal{C}}\! e^{-t\tau} \! \bigg(\!
 \mbox{Tr}\,\big(\psi_{-1} (x',\xi',\tau)\big) \\
 &&\quad  + \,\mbox{Tr}\,\big(\psi_{-2} (x',\xi',\tau)\big)\!+\! \mbox{Tr}\,\big(\psi_{-3} (x',\xi',\tau)\big) \!+\!\cdots \!\bigg) d\tau\!\bigg\} d\xi'\, dS(x'),\nonumber\end{eqnarray}
We will calculate the asymptotic expansion of the trace of the semigroup $e^{-t\Xi}$ as $t\to 0^+$. More precisely, we will figure out
 \begin{eqnarray} \label{18/12/2-5} \widetilde{a_m}(x') \!\!=\!\! \frac{1}{(2\pi)^{n\!-\!1}}\! \int_{{\Bbb R}^{n\!-\!1}}\! \frac{i}{2\pi} \!\int_{\mathcal{C}}\! e^{-t\tau}
 \mbox{Tr}\,\big(\psi_{\!-\!1\!-\!m} (x'\!,\xi',\tau)\big) d\tau\!\bigg\} d\xi'\quad \mbox{for}\;\,  0\le m\le n\!-\!1.\end{eqnarray}

\vskip 0.26 true cm

 Step 2. \   In order to make the eigenvalues are nonnegative, we have taken the principal symbol of $p_1$ is positive. By  (\ref{19.3.28-13}), we have
\begin{eqnarray}&& \label{18/11/24} \end{eqnarray} \begin{eqnarray}
\!\!&& \!\!\!\!\!\! \!\!  p_1(x', \xi')-\tau I_n \qquad \qquad \qquad \qquad \nonumber \\
\!\!&&\!\!\!\!\!\!\!\!  \quad = {\begin{small} \left[ \begin{BMAT}(@, 1pt, 1pt){c.c}{c.c}\bigg[s_1\mu \sqrt{\sum\limits_{\alpha,\beta} g^{\alpha\beta} \xi_\alpha \xi_\beta}\, \delta_{jk} + s_2 \mu \frac{\sum_\alpha g^{j\alpha}\xi_\alpha\xi_k}{\sqrt{\sum\limits_{\alpha,\beta} g^{\alpha\beta}\xi_\alpha\xi_\beta}}-\tau\delta_{jk}\bigg] &
 \bigg[\mu (s_4-i) \sum_{\alpha} g^{j\alpha}\xi_\alpha\bigg] \\
 \bigg[\big(s_5 (\lambda+2\mu)-i \lambda\big)\xi_k\bigg] &  (s_1+s_3)(\lambda+2\mu) \sqrt{\sum\limits_{\alpha,\beta}g^{\alpha\beta} \xi_\alpha\xi_\beta}-\tau \end{BMAT}\right].\end{small}}\nonumber \end{eqnarray}
 Let us denote
 \begin{eqnarray*} p_1-\tau I_n:=\begin{small}\left[\begin{BMAT}(@, 3pt, 3pt){c.c}{c.c} \left[ p_1^{jk}-\tau \delta_{jk}\right]_{(n-1)\times (n-1)} & \left[ p_1^{jn}\right]_{(n-1)\times 1}\\
 \left[ p_1^{nk}\right]_{1\times (n-1)} & p_1^{nn}-\tau\end{BMAT}\right],\end{small}\end{eqnarray*}
 where \begin{eqnarray*} &&\!\!\!\! \!\!\!\!\!p_1^{jk}\! -\! \tau \delta_{jk}\!=\! s_1\mu \sqrt{\sum\limits_{\alpha,\beta}\! g^{\alpha\beta} \xi_\alpha \xi_\beta} \;\delta_{jk}\! +\! s_2 \mu \frac{\sum_\alpha g^{j\alpha}\xi_\alpha\xi_k}{\sqrt{\sum\limits_{\alpha,\beta}\!g^{\alpha\beta}\xi_\alpha\xi_\beta}}\!-\!\tau\delta_{jk},\;\;\, 1\le j,k\le n\!-\!1,\\
&&\!\! \!\!\!\!\! \!\!p_1^{jn}=  \mu (s_4-i) \sum_\alpha g^{j\alpha}\xi_\alpha, \,\quad \, 1\le j\le n-1\\
&&\!\!\!\!\!\!\!\!\! p_1^{nk}=\big(s_5 (\lambda+2\mu)-i \lambda\big)\xi_k,\,\quad \, 1\le k\le n-1\\
&&\!\!\!\! \!\!\!\!\! p_1^{nn}-\tau=(s_1+s_3)(\lambda+2\mu) \sqrt{\sum\limits_{\alpha,\beta}g^{\alpha\beta} \xi_\alpha\xi_\beta}-\tau.\end{eqnarray*}
 It is easy to verify that  \begin{eqnarray*} (p_1(x', \xi')-\tau I_n)^{-1}:= \begin{small}\left[\begin{BMAT}(@, 3pt, 3pt){c.c}{c.c} \sigma & \eta \\
 \zeta &  \theta \end{BMAT}\right]\end{small},
 \end{eqnarray*}
where \begin{eqnarray*} \!\sigma\! &=&\! \left(\left[ p_1^{jk} - \tau \delta_{jk}\right]_{(n-1)\times (n-1)} -\frac{1}{p_1^{nn}-\tau} \left[\, p_1^{jn}\,\right]_{(n-1)\times 1} \left[\, p_1^{nk}\,\right]_{1\times (n-1)}\right)^{-1},\\
\zeta \!&=&\! -\frac{1}{p_1^{nn}-\tau} \left[p_1^{nk}\right]_{1\times (n-1)}   \left(\left[ p_1^{jk} - \tau \delta_{jk}\right]_{(n-1)\times (n-1)} \right.\\
&&\left.-\frac{1}{p_1^{nn}-\tau} \left[ \,p_1^{jn}\,\right]_{(n-1)\times 1} \left[\, p_1^{nk}\,\right]_{1\times (n-1)}\right)^{-1},\\
\eta \!&=&\! -\frac{1}{p_1^{nn}-\tau}  \left(\left[ p_1^{jk} - \tau \delta_{jk}\right]_{(n-1)\times (n-1)}\right.\\
 && \left. -\frac{1}{p_1^{nn}-\tau} \left[ \, p_1^{jn}\,\right]_{(n-1)\times 1} \left[ \,p_1^{nk}\,\right]_{1\times (n-1)}\right)^{-1}\left[\,p_1^{jn}\,\right]_{(n-1)\times 1},\\
\theta \!&=& \!\frac{1}{p_1^{nn}-\tau}  \left\{ 1+ \frac{\left[\,p_1^{nk}\,\right]_{1\times (n-1)}}{(p_1^{nn}-\tau)}\left(\left[ p_1^{jk} - \tau \delta_{jk}\right]_{(n-1)\times (n-1)}\right.\right.\\
 && \left.\left.-\frac{1}{p_1^{nn}-\tau} \left[\, p_1^{jn}\,\right]_{(n-1)\times 1} \left[ \,p_1^{nk}\,\right]_{1\times (n-1)}\right)^{-1}\right\}.
\end{eqnarray*}
In particular, at the origin $x_0$ of the boundary normal coordinates (see (\ref{18/7/14/1})) we get
\begin{eqnarray} \label{19.4.14=1} \end{eqnarray} \begin{eqnarray} \sigma\!=\! \left[(s_1\mu |\xi'|\!-\!\tau)\delta_{jk}\! +\!\left(\frac{s_2\mu }{|\xi'|}\!-\! \frac{\mu(s_4-i) \big(s_5( \lambda\!+\!2\mu)-i\lambda\big)}{(s_1\!+\!s_3)(\lambda\!+\!2\mu)|\xi'|-\tau}\right)\xi_j\xi_k\right]_{(n\!-\!1)\times (n\!-\!1)}^{-1}.\nonumber\end{eqnarray}
Applying Lemma 5.1 to the expression (\ref{19.4.14=1}) of $\sigma$, we have
\begin{eqnarray}\label{19.4.14-,6} \sigma= \frac{1}{s_1\mu|\xi'|-\tau} \big[ \delta_{jk} - \omega \,\xi_j\xi_k\big]_{(n-1)\times(n-1)},\end{eqnarray}
  where \begin{eqnarray*}\omega\!=\! \frac{\mu s_2 \big( (s_1+s_3) (\lambda+2\mu) |\xi'|-\tau\big)-|\xi'| \mu (s_4-i)\big(s_5(\lambda+2\mu)-i\lambda\big)}{\big((s_1\!+\!s_3)(\lambda\!+\!2\mu)|\xi'|\!-\!\tau\big)\big( \mu (s_1\!+\!s_2)|\xi'|\!-\!\tau\big)|\xi'| \!-\! \mu (s_4\!-i) \big( s_5(\lambda\!+\!2\mu)\! -i\lambda\big)|\xi'|^3}.\end{eqnarray*}
Obviously, $\omega$ can be rewritten as
\begin{align}\label{19.6.13.2} \omega &\!=\! \frac{\mu s_2 \big( (s_1\!+\!s_3) (\lambda+2\mu) |\xi'|-\tau\big)-|\xi'| \mu (s_4-i)\big(s_5(\lambda+2\mu)-i\lambda\big)}{|\xi'| (\tau \!-\!\tau_+)(\tau\!-\!\tau_-)}\\
&  \! =\! \frac{ \mu |\xi'|\big[ s_2(s_1+s_3)(\lambda+2\mu) - (s_4-i)\big( s_5(\lambda+2\mu)-i\lambda\big)\big] -\mu s_2 \tau_{+}}{ |\xi'|(\tau\!-\!\tau_{+})(\tau\!-\!\tau_{-})} \!-\!\frac{\mu s_2}{|\xi'|(\tau\!-\!\tau_{-})}
.\nonumber\end{align}
 where $\tau_+=2\mu |\xi'|$, $\, \tau_-= \frac{2\mu(\lambda+\mu)|\xi'|}{\lambda+3\mu}$.
Therefore,
\begin{align} \label{18/11/22-1} &\psi_{-1}(x_0,\xi',\tau):=\big(p_1 (x_0, \xi') -\tau I_n\big)^{-1}=\begin{small} \left[\begin{BMAT}(@, 9pt, 9pt){c.c}{c.c} \sigma& \eta\\
\zeta & \theta \end{BMAT}\right]\end{small}\\
&\!\!\!=\! \!{\tiny \left[\!\begin{BMAT}(@, 0.1pt, 0.1pt){c.c}{c.c} \sigma & -\frac{(s_4-i)\mu}{(s_1+s_3)(\lambda+2\mu)|\xi'|-\tau} \sigma \left[\xi_j\right]_{(n\!-\!1)\!\times\! 1}\\
 \!-\frac{s_5 (\lambda\!+\!2\mu)-i\lambda}{ (s_1\!+\!s_3)(\lambda\!+\!2\mu)|\xi'|\!-\!\tau} \left[\xi_k\right]_{1\!\times\!(n\!-\!1)} \sigma &
 \frac{1}{(s_1\!+\!s_3)(\lambda\!+\!2\mu)|\xi'|\!-\!\tau} \!+\!\frac{\mu (s_4\!-\!i) \big(s_5 (\lambda\!+\!2\mu)\!-\!i\lambda\big)} {\big((s_1\!+\!s_3)(\lambda\!+\!2\mu)|\xi'|\!-\!\tau\big)^2} \left[\xi_k\right]_{1\!\times \!(n\!-\!1)} \sigma \left[\xi_j\right]_{(n\!-\!1)\!\times \!1}
  \end{BMAT}\!\right]}\nonumber\quad \;\, \\
  &\!\!=\! \!{\tiny \left[\!\begin{BMAT}(@, 0.1pt, 0.1pt){c.c}{c.c}\!
 -\frac{1}{\tau \!-\!\mu |\xi'|s_1}\left[ \delta_{jk}\!-\!\omega \xi_j\xi_k \right] _{( n\!-\!1) \times ( n\!-\!1 )}&  -\frac{\mu ( s_4\!-\!i ) ( 1\!-\!\omega |\xi'|^2 )}{\left( \tau \!-\!\mu |\xi'|s_1 \right) \left( \tau \!-\!|\xi'| \left( s_1\!+\!s_3 \right) \left( \lambda \!+\!\text{2}\mu \right) \!\right)}\left[ \xi _j \right]_{\left( n\!-\!1 \right) \times 1}\\
 \!\!\frac{\left( -( \lambda\! +\!\text{2}\mu )s_5\!+\!\lambda i \right) \left( 1\!-\!\omega |\xi'|^2 \right)}{( \tau\! -\!\mu |\xi'| s_1 ) ( \tau \!-\!|\xi'|( s_1\!+\!s_3 ) ( \lambda\! +\!\text{2}\mu ) )}\left[ \xi _k \right] _{1\!\times \!( n\!-\!1 )}&  \frac{
 1}{ (s_1\!+\!s_3)(\lambda\!+\!2\mu)|\xi|\!-\!\tau}\!+\! \frac{\mu \left( ( \lambda \!+\!2\mu ) s_5\!-\!\lambda i \right)\! ( s_4\!-\!i )  ( 1\!-\!\omega |\xi'|^2 ) |\xi'| ^2}{\left(  s_1\mu |\xi'| \!-\!\tau \right) \left( \tau \!-\!|\xi'| ( s_1\!+\!s_3 ) ( \lambda \!+\!{2}\mu ) \right) ^2}\\
\end{BMAT} \!\!\right]\!}\nonumber
  \end{align}
It follows that  \begin{eqnarray*}\!\!\! &&\!\!\!\!\!\!\! \!\!\! \label{19.4.14-,7} \mbox{Tr\,}(\sigma)\!=\! \frac{1}{ s_1\mu|\xi'|\!-\!\tau}  \bigg\{
n-1\\
&& \!\!-\frac{s_2 \mu \big((s_1\!+\!s_3) (\lambda\!+\!2\mu) |\xi'|\! -\!\tau\big) |\xi'|- \mu (s_4 \!-\!i) \big( s_5 (\lambda\!+\!2\mu) -i\lambda\big)|\xi'|^2}{\big((s_1\!+\!s_3) (\lambda\!+\!2\mu) |\xi'| \!-\!\tau\big) \big( (s_1 \!+\!s_2)\mu|\xi'|\!-\!\tau\big)\!-\!\mu (s_4\!-\!i)\big( s_5(\lambda\!+\!2\mu)\!-\!i\lambda\big)|\xi'|^2}\!\bigg\}\end{eqnarray*}
and
\begin{eqnarray*} &&\!\!\!\!\!\!\!\!\!\!\!\!\!\!\!\theta\!=\! \frac{1}{(s_1+s_3)(\lambda+2\mu)|\xi'|-\tau} \left\{ 1+ \frac{\big( s_5 (\lambda+2\mu)-i\lambda\big) \mu (s_4-i) |\xi'|^2}{
(s_1\mu|\xi'|-\tau) \big( (s_1+s_3) (\lambda+2\mu)|\xi'|-\tau\big) }\bigg( 1 \right.\\
&&\!\left.\!\!\!\!\!-\frac{s_2 \mu \big((s_1\!+\!s_3) (\lambda\!+\!2\mu) |\xi'| \!-\!\tau\big) |\xi'|\!-\! \mu (s_4 -i) \big( s_5 (\lambda\!+\!2\mu) -i\lambda\big)|\xi'|^2}{\big((s_1\!+\!s_3) (\lambda\!+\!2\mu) |\xi'| \!-\!\tau\big) \big( (s_1 \!+\!s_2)\mu|\xi'|\!-\!\tau\big)\!-\!\mu (s_4\!-\!i)\big( s_5(\lambda\!+\!2\mu)\!-\!i\lambda\big)|\xi'|^2}  |\xi'|\!\bigg)\!\right\}.\end{eqnarray*}
Thus we have
\begin{align} \label{19.4.14;8} & \!\!\!\!\mbox{Tr}\big(\psi_{-1}( x_0,\xi',\tau)\big)=\mbox{Tr}\, (\sigma)+\theta
\\
&=-\frac{n-1}{\tau -\mu | \xi ' |s_1}-\frac{1}{\tau -| \xi' |( \lambda +\text{2}\mu ) ( s_1+s_3 )}-\frac{\mu | \xi ' |s_2}{\left( \tau -\tau _+ \right) ( \tau -\tau _- )}\nonumber
\\
&\;\;\,\;+\frac{\mu | \xi' |^2\big( s_2\left( ( \lambda +\text{2}\mu ) ( s_1+s_3 ) -\mu s_1 \right) +\left( -s_5( \lambda +\text{2}\mu ) +\lambda {i} \right) ( s_4-{i} ) \big)}{( \tau \!-\!\tau _+ ) ( \tau \!-\!\tau _- ) \left( \tau \!-\!\mu | \xi' |s_1 \right)}\nonumber
\\
&\;\;\,\;+\frac{\mu | \xi' |^2\left( -s_5( \lambda +{2}\mu ) +\lambda {i} \right) ( s_4-{i} )}{\left( \tau \!-\!| \xi' |( \lambda \!+\!\text{2}\mu ) ( s_1\!+\!s_3 ) \right)^2\left( \tau \!-\!\mu | \xi' |s_1 \right)}\!+\!\frac{\mu ^2| \xi ' |^3s_2\big(\! -\!s_5( \lambda \!+\!\text{2}\mu ) \!+\!\lambda {i} \big) ( s_4\!-{i} )}{( \tau \!-\!\tau _+ ) ( \tau \!-\!\tau _- ) \left( \tau \!-\!| \xi' |( \lambda \!+\!\text{2}\mu ) ( s_1\!+\!s_3 ) \right) ^2}\nonumber
\\
&\;\;\,\;-\!\frac{\mu ^2| \xi' |^4\!\left\{\! s_2\big( ( \lambda \!+\!\text{2}\mu ) ( s_1\!+\!s_3 )\! -\!\mu s_1 \big) \!+\!\left( \!-\!s_5( \lambda \!+\!\text{2}\mu ) \!+\!\lambda {i} \right) ( s_4\!-\!{i} ) \!\right\}\!\! \big(\!\! -\!s_5( \lambda \!+\!\text{2}\mu ) \!+\!\lambda {i} \big) \!( s_4\!-\!{i} )}{( \tau \!-\!\tau _+ ) ( \tau \!-\!\tau _- ) \big( \tau \!-\!| \xi ' |( \lambda \!+\!\text{2}\mu ) ( s_1\!+\!s_3 ) \big) ^2\left( \tau \!-\mu | \xi ' |s_1 \!\right)}\nonumber\\
& \,\;:=II_1^{(0)} +II_2^{(0)} +II_3^{(0)}+II_4^{(0)} +II_5+II_6^{(0)}+II_7^{(0)}, \nonumber\end{align}
where \begin{align*} & II_1^{(0)}=-\frac{n-1}{\tau -\mu | \xi ' |},\\
& II_2^{(0)} = \frac{\lambda +\text{3}\mu}{\text{4}\mu ^2| \xi ' |-\lambda \tau -\text{3}\mu \tau +\text{2}\lambda \mu | \xi ' |}\\
& II_3^{(0)} =\frac{\mu | \xi ' |( \lambda +\mu )}{( \tau -\text{2}\mu | \xi' | ) \left( \text{2}\mu^2| \xi ' |-\lambda \tau -\text{3}\mu \tau +\text{2}\lambda \mu | \xi '| \right)},
\\
&II_4^{(0)}=\frac{\mu ^2| \xi ' |^2( \lambda -\mu )}{( \tau -\mu | \xi ' | ) ( \text{2}\mu| \xi '|-\tau ) \left( \text{2}\mu ^2| \xi '|-\lambda\tau -\text{3}\mu \tau +\text{2}\lambda \mu| \xi ' | \right)},\\
&II_5^{(0)}= \frac{\text{4}\mu ^4| \xi ' |^2}{( \mu | \xi ' |-\tau ) \left( \text{4}\mu ^2| \xi ' |-\lambda \tau -\text{3}\mu \tau +\text{2}\lambda \mu | \xi ' | \right) ^2},\\
&II_6^{(0)} =\frac{-\text{4}\mu ^5| \xi '|^3( \lambda +\mu )}{( \text{2}\mu | \xi ' |\!-\!\tau ) \!\left(\! \text{2}\mu^2| \xi ' |\!-\!\lambda \tau \!-\!\text{3}\mu \tau \!+\!\text{2}\lambda \mu | \xi ' | \right) \!\left(\! \text{4}\mu ^2| \xi ' |\!-\!\lambda \tau \!-\!\text{3}\mu \tau \!+\!\text{2}\lambda \mu | \xi ' | \right)^2},\\
&II_7^{(0)}\!\! =\!\frac{-\text{4}\mu ^6| \xi' |^4( \lambda -\mu )}{( \mu | \xi' |\!-\!\tau )\! \left(\! \text{2}\mu | \xi ' |\!-\!\tau \right)\! \left(\! \text{2}\mu ^2| \xi'|\!-\!\lambda \tau \!-\!\text{3}\mu \tau \!+\!\text{2}\lambda \mu | \xi ' |\right)\!\left( \!\text{4}\mu ^2| \xi ' |\!-\!\lambda \tau \!-\!\text{3}\mu \tau\! +\!\text{2}\lambda \mu | \xi ' | \right)^2}.
\end{align*}
It is easily seen that  $\tau_+, \,\tau_-$, $\mu |\xi'|$ and $(\lambda+2\mu)(s_1+s_3)|\xi'|$) are distinct numbers for any $|\xi'|\ne 0$. By applying residue theorem of complex variable function (see \cite{Ahl}), we find from (\ref{19.4.14;8}) that
\begin{align*}
&\frac{1}{2\pi i}\int_{\mathcal{C}}{e^{-t\tau}}(II_1^{(0)})\,d\tau =-( n-1)\text{e}^{-\mu |\xi'|t},
\\
&\frac{1}{2\pi i}\int_{\mathcal{C}}{e^{-t\tau}}(II_2^{(0)})\,d\tau =-\text{e}^{-\frac{\text{2}\mu \left( \lambda +\text{2}\mu \right)|\xi'|t}{\lambda +\text{3}\mu}},
\\
&\frac{1}{2\pi i}\int_{\mathcal{C}}{e^{-t\tau}}(II_3^{(0)})\,d\tau =\frac{\lambda +\mu}{\text{4}\mu}\text{e}^{-\frac{\text{2}\mu \left( \lambda +\mu \right)|\xi'|t}{\lambda +\text{3}\mu}}-\frac{\lambda +\mu}{4\mu }\text{e}^{-\text{2}\mu |\xi'|t},
\\
&\frac{1}{2\pi i}\int_{\mathcal{C}}{e^{-t\tau}}(II_4^{(0)})\,d\tau =-\frac{\lambda +\text{3}\mu}{\text{4}\mu }\text{e}^{-\frac{\text{2}\mu \left( \lambda +\mu \right)|\xi'|t}{\lambda +\text{3}\mu}}+\frac{\lambda -\mu}{4\mu }\text{e}^{-\text{2}\mu |\xi'|t}+\text{e}^{-\mu |\xi'|t},
\\
&\frac{1}{2\pi i}\int_{\mathcal{C}}{e^{-t\tau}}(II_5^{(0)})\,d\tau =\frac{4\mu ^2\left( \lambda +\text{3}\mu +\mu ^2|\xi'|t+\lambda \mu |\xi'|t \right)}{\left( \lambda +\text{3}\mu \right) ( \lambda +\mu ) ^2}\text{e}^{-\frac{\text{2}\mu \left( \lambda +\text{2}\mu \right)|\xi'|t}{\lambda +\text{3}\mu}}\\
& \qquad \qquad \qquad \qquad \qquad \,-\frac{\text{4}\mu ^2}{( \lambda +\mu ) ^2}\text{e}^{-\mu |\xi'|t},
\\
&\frac{1}{2\pi i}\int_{\mathcal{C}}{e^{-t\tau}}(II_6^{(0)})\,d\tau =-\frac{\mu |\xi'|t( \lambda +\mu )}{\left( \lambda +\text{3}\mu \right)}\text{e}^{-\frac{\text{2}\mu \left( \lambda +\text{2}\mu \right)|\xi'|t}{\lambda +\text{3}\mu}}+\frac{\lambda +\mu}{4\mu \,}\text{e}^{-\frac{\text{2}\mu \left( \lambda +\mu \right)|\xi'|t}{\lambda +\text{3}\mu}}\\
& \qquad \qquad \qquad\, \qquad \qquad \, -\frac{\lambda +\mu}{4\mu \,}\text{e}^{-{2}\mu |\xi'|t},
\\
&\frac{1}{2\pi i}\int_{\mathcal{C}}{e^{-t\tau}}(II_7^{(0)})\,d\tau =\frac{\left( \lambda -\mu \right)\left( \lambda +\text{3}\mu +\mu ^2|\xi'|t+\lambda \mu |\xi'|t \right)}{( \lambda +\mu ) ^2}\text{e}^{-\frac{\text{2}\mu \left( \lambda +\text{2}\mu \right)|\xi'|t}{\lambda +\text{3}\mu}}\\
&\qquad \qquad \qquad \qquad \qquad \,  -\frac{\lambda +\text{3}\mu}{\text{4}\mu }\text{e}^{-\frac{\text{2}\mu \left( \lambda +\mu \right)|\xi'|t}{\lambda +\text{3}\mu}} +\frac{\lambda -\mu}{4\mu \,}\text{e}^{-\text{2}\mu |\xi'|t}+\frac{\text{4}\mu ^2}{\left( \lambda +\mu \right) ^2}\text{e}^{-\mu |\xi'|t}.\end{align*}
Thus
\begin{align*}
\frac{i}{2\pi}\int_{\mathcal{C}}{e^{-t\tau}}\left(\mbox{Tr}\,\big( ( p_1-\tau I_n)^{-1} \big)\right) \,d\tau& =-\frac{1}{2\pi i}\int_{\mathcal{C}}{e^{-t\tau}}  \bigg(\sum_{l=1}^7 II_{l}^{(0)}\bigg) \, d\tau
\\
&=\text{e}^{-\text{2}\mu |\xi'|t}+\text{e}^{-\frac{\text{2}\mu\left( \lambda +\mu \right) |\xi'|t}{\lambda +\text{3}\mu}}+\left( n-2 \right) \text{e}^{-\mu |\xi'|t}.
\end{align*}
It follows from \cite{Liu2} that for any positive constant $C$,
\begin{align} \label{19.6.15-9,1}
\int_{\mathbb{R}^{n-1}}{e^{-C\sqrt{\sum_{j=1}^{n-1}{\xi_{j}^{2}}}}} \,d\xi _1\,d\xi _2\cdots \,d\xi _{n-1}
= \frac{\text{vol}( \mathbb{S}^{n-2} ) \varGamma ( n-1)  }{C^{n-1}}  \quad \mbox{for}\;\; n\geqslant \text{2.}
\end{align}
 where $\text{vol}\left( \mathbb{S}^{n-2} \right) = \frac{2\pi^{(n-1)/2}}{\Gamma(\frac{n-1}{2})} $ is the area of $(n-2)$-dimensional unit sphere ${\Bbb S}^{n-2}$ in ${\Bbb R}^{n-1}$, and $\Gamma$ is the gamma function.
Therefore we have  \begin{align}\label{19.6.14-1}
\widetilde{a_0}(x')=&\frac{1}{(2\pi)^{n-1}}\!\int_{{\Bbb R}^{n-1}}\!{\left( \!\frac{{i}}{2\pi}\!\int_{\mathcal{C}}{e^{-t\tau}} \,\mbox{Tr} ( p_1\!-\!\tau I_n )^{-1} \,d\tau \right) d\xi _1\,d\xi _2\cdots d\xi _{n-1}}\\
=& \frac{1}{(2\pi)^{n-1}}\!\int_{{\Bbb R}^{n-1}}\!\left( \text{e}^{-\text{2}\mu |\xi'|t}\!+\!\text{e}^{-\frac{\text{2}\mu \left( \lambda \!+\!\mu \right)|\xi'|t}{\lambda \!+\!\text{3}\mu}}\!+\!\left( n-2 \right) \text{e}^{-\mu |\xi'|t}\right) d\xi_1\cdots d\xi_{n-1} \nonumber\\
=&\, t^{1-n} \,\frac{\left( n-2 \right) !\,\text{vol}( \mathbb{S}^{n-2} ) }{(2\pi)^{n-1}\mu ^{n-1}}\bigg( \frac{1+\big( \frac{\lambda +\text{3\,}\mu}{\lambda +\mu} \big) ^{n-1}}{2^{n-1}}+\left( n-2 \right) \bigg):=t^{1-n}\, a_0(x').\nonumber\end{align}

\vskip 0.26 true cm

 Step 3. \  At the origin $x_0=(0,\cdots,0)$ of the boundary normal coordinates, by (\ref{19.3.23-1}) and  (\ref{18/7/14/1}) we have
  \begin{flalign*}\label{19.4.24-1}&E_1= \big(b_0 q_1\big)\big|_{x_0} + \frac{\partial q_1}{\partial x_n}\big|_{x_0}-c_1\qquad \qquad \qquad \qquad \qquad \qquad \\
  &\!=\!{\tiny \left[\!\begin{BMAT}(@, 1pt, 1pt){c.c}{c.c}  \left[\frac{1}{2}\! \sum_{\alpha}\! \frac{\partial g_{\alpha\alpha}}{\partial x_n}  \,\delta_{jk}\!+\!\frac{\partial g_{jk}}{\partial x_n}\right]_{(n\!-\!1)\times (n\!-\!1)} &  \left[0\right]_{(n\!-\!1)\times 1}\\
  \left[0\right]_{1\times (n-1)} & \frac{1}{2} \sum_{\alpha} \frac{\partial g_{\alpha \alpha}}{\partial x_n}\end{BMAT} \right]\!
    \left[\begin{BMAT}(@, 1pt, 1pt){c.c}{c.c}  \left[s_1 |\xi'| \!+\!s_2 \frac{\xi_j\xi_k}{|\xi'|} \right]_{(n\!-\!1)\times (n\!-\!1)}  &  \left[s_4\,\xi_j\right]_{(n\!-\!1)\times 1}\\
  \left[s_5 \,\xi_k \right]_{1\times (n-1)} & (s_1+s_3)|\xi'|\end{BMAT} \right]
  } \\
  &\!+ \!  {\tiny \left[\begin{BMAT}(@, 1pt, 1pt){c.c}{c.c} \!\!\left[\!\frac{s_1 }{2|\xi'|}\! \sum_{\alpha,\beta} \! \frac{\partial g^{\alpha\beta}}{\partial x_n} \xi_\alpha \xi_\beta \delta_{jk}
   \! +\!\frac{s_2}{|\xi'|}\sum_{\alpha}\! \frac{\partial g^{j\alpha}}{\partial x_n} \xi_\alpha \xi_k \!-\!\frac{s_2}{2|\xi'|^3} \sum_{\alpha, \beta}
    \!\frac{\partial g^{\alpha\beta}}{\partial x_n} \xi_\alpha \xi_\beta  \xi_j\xi_k \!\right] & \left[s_4 \sum_{\alpha} \!\frac{\partial g^{j\alpha}}{\partial x_n}  \xi_\alpha  \right]_{(n\!-\!1)\times 1}\\
  \left[0\right]_{1\times (n-1)} & \frac{s_1\!+\!s_2}{2|\xi'|} \sum_{\alpha,\beta}\! \frac{\partial g^{\alpha \beta}}{\partial x_n} \xi_\alpha \xi_\beta\end{BMAT} \right]
  } \\
  &\!- \! \begin{small} \left[\begin{BMAT}(@, 1pt, 1pt){c.c}{c.c} \begin{small} \left[0 \right]_{(n-1)\times (n-1)} \end{small} &  \left[\frac{i(\lambda+\mu)}{2\mu} \sum_{\alpha} \frac{\partial g_{\alpha\alpha}}{\partial x_n} \, \xi_j +i \sum\limits_{\alpha} \frac{\partial g_{j\alpha}}{\partial x_n}\xi_\alpha\right]_{(n-1)\times 1}\\
  \left[-\frac{i\mu}{\lambda+2\mu} \sum\limits_{\alpha} \frac{\partial g_{k\alpha}}{\partial x_n}\xi_\alpha\right]_{1\times (n-1)} & 0\end{BMAT} \right]
  \end{small}.  \end{flalign*}
  Since $\sum_{\gamma} g^{\alpha \gamma} g_{\gamma \beta}=\delta_{\alpha\beta}$ for all $1\le \alpha, \beta\le n-1$, we get
  $$\sum_{\gamma} \big( \frac{\partial g^{\alpha \gamma}}{\partial x_n} \, g_{\gamma\beta} +g^{\alpha \gamma}\frac{\partial g_{\gamma \beta}}{\partial x_n} \big)\big|_{x_0} = 0,$$ i.e., $$\frac{\partial g^{\alpha\beta}}{\partial x_n}\bigg|_{x_0} =-\frac{\partial g_{\alpha \beta}}{\partial x_n}\bigg|_{x_0}\quad \mbox{for all}\;\; 1\le \alpha, \beta\le n-1.$$
Thus, applying this  and (\ref{18/7/14/1}) once more,  we get
\begin{eqnarray} && \!\!\!\!\! \label{2021.6.23-1} \end{eqnarray} \begin{eqnarray}
 \!\!\!\!&&\!\!\!\!\!\!\!\!\! E_1 =s_1\big( \sum_{\alpha} \kappa_\alpha |\xi'| -\sum_{\alpha} \frac{\kappa_\alpha \xi_\alpha^2}{|\xi'|}\big) I_n\nonumber\\
\!\!\!\!&&\!\!\!\!\! +\!
 {\begin{small} \left[\!\begin{BMAT}(@, 0pt, 0pt){c.c}{c.c} \left[\!2s_1 |\xi'| \kappa_j\delta_{jk} \!+\! s_2 \big( \sum_{\alpha} \kappa_\alpha \frac{\xi_j\xi_k}{|\xi'|}  \!+\!\sum_{\alpha} \kappa_\alpha \xi_\alpha^2\, \frac{\xi_j\xi_k}{|\xi'|^3}\big) \right]  &  \left[\! \big( s_4 \!-\!\frac{i(\lambda\!+\!\mu)}{\mu} \big) \big(\sum_{\alpha} \kappa_\alpha\big) \xi_j \!-\!2 \, i\,\kappa_j\xi_j\!\right]_{(n\!-\!1)\times 1}\\
  \left[s_5 \big(\sum_{\alpha} \kappa_\alpha \big) \xi_k \!+\!\frac{2i\mu}{\lambda\!+\!2\mu} \kappa_k \xi_k\right]_{1\times (n\!-\!1)} & s_3 \big(\sum_\alpha \! \kappa_\alpha |\xi'| -\sum_{\alpha} \frac{\kappa_\alpha \xi_\alpha^2}{|\xi'|} \big) \end{BMAT}\! \right]\end{small}},\nonumber \end{eqnarray}
  where $\kappa_\alpha =\frac{1}{2}\frac{g_{\alpha\alpha}}{\partial x_n}\big|_{x_0}$, $1\le \alpha \le n-1$ is the principal curvature of $\partial \Omega$ at the origin $x_0$.
At the origin (see (\ref{18/7/14/1})) of the boundary normal coordinates, (\ref{19.3.28-1}) becomes
\begin{eqnarray}\label{19.4.24-6}\end{eqnarray} \begin{eqnarray}
 &&\!\! \!\!\!\!\!\!\!\!\!\! X = {\tilde{s}}_1\frac{1}{|\xi'|} E +  \begin{small} \left[\begin{BMAT}(@, 1pt, 1pt){c.c}{c.c} \begin{small} \left[{\tilde{s}}_2\frac{\xi_j\xi_k}{|\xi'|^3} \right]_{(n-1)\times (n-1)} \end{small} &  \left[ {\tilde{s}}_4 \frac{\xi_j}{|\xi'|^2}\right]_{(n-1)\times 1}\\
  \left[{\tilde{s}}_5 \frac{\xi_k}{|\xi'|^2}  \right]_{1\times (n-1)} & {\tilde{s}}_3 \frac{1}{|\xi'|}  \end{BMAT} \right]  \end{small}E \nonumber\\
   && \,+ E \begin{small} \left[\begin{BMAT}(@, 1pt, 1pt){c.c}{c.c} \begin{small} \left[{\tilde{s}}_6\frac{\xi_j\xi_k}{|\xi'|^3} \right]_{(n-1)\times (n-1)} \end{small} &  \left[ {\tilde{s}}_8 \frac{\xi_j}{|\xi'|^2}\right]_{(n-1)\times 1}\\
  \left[{\tilde{s}}_9 \frac{\xi_k}{|\xi'|^2}  \right]_{1\times (n-1)} & {\tilde{s}}_7 \frac{1}{|\xi'|}  \end{BMAT} \right]  \end{small} \nonumber\\
 && \, + \begin{small} \left[\!\begin{BMAT}(@, 1pt, 1pt){c.c}{c.c} \begin{small} \left[{\tilde{s}}_{10}\frac{\xi_j\xi_k}{|\xi'|^3} \right]_{(n\!-\!1)\times (n\!-\!1)} \end{small} &  \left[ {\tilde{s}}_{18} \frac{\xi_j}{|\xi'|^2}\right]_{(n\!-\!1)\times 1}\\
  \left[{\tilde{s}}_{19} \frac{\xi_k}{|\xi'|^2}  \right]_{1\times (n\!-\!1)} & {\tilde{s}}_{11} \frac{1}{|\xi'|}  \end{BMAT} \right]  \end{small}
E \begin{small} \left[\begin{BMAT}(@, 1pt, 1pt){c.c}{c.c} \begin{small} \left[\frac{\xi_j\xi_k}{|\xi'|^2}  \right]_{(n\!-\!1)\times (n\!-\!1)} \end{small} &  \left[ 0\right]_{(n\!-\!1)\times 1}\\
  \left[0  \right]_{1\times (n\!-\!1)} & -1  \end{BMAT} \right]  \end{small}\nonumber
\\  &&\,  + \begin{small} \left[\begin{BMAT}(@, 1pt, 1pt){c.c}{c.c} \begin{small} \left[{\tilde{s}}_{14}\frac{\xi_j\xi_k}{|\xi'|^4} \right]_{(n-1)\times (n-1)} \end{small} &  \left[ {\tilde{s}}_{22} \frac{\xi_j}{|\xi'|^3}\right]_{(n-1)\times 1}\\
  \left[{\tilde{s}}_{23} \frac{\xi_k}{|\xi'|^3}  \right]_{1\times (n-1)} & {\tilde{s}}_{15} \frac{1}{|\xi'|^2}  \end{BMAT} \right]  \end{small} E
 \begin{small} \left[\begin{BMAT}(@, 1pt, 1pt){c.c}{c.c} \begin{small} \left[0\right]_{(n-1)\times (n-1)} \end{small} &  \left[ \xi_j\right]_{(n-1)\times 1}\\
  \left[\xi_k  \right]_{1\times (n-1)} & 0  \end{BMAT} \right]  \end{small}.\nonumber \end{eqnarray}
  Now, by replacing $E$ and $X$ in (\ref{19.4.24-6}) by $E_1$ and $q_0$ at the origin $x_0$ of the boundary normal coordinates , respectively, we get \begin{eqnarray} \label{19.4.25} \end{eqnarray} \begin{eqnarray} &&  q_0= s_1 {\tilde{s}}_1\! \left(\! \!\big(\!\sum_{\alpha}\kappa_\alpha\!\big) \!-\!\sum_{\alpha}\!\frac{ \kappa_\alpha \xi_\alpha^2}{|\xi'|^2}\!\right)I_n \!+\! \begin{small} \left[\begin{BMAT}(@, 1pt, 1pt){c.c}{c.c} \begin{small} \!\left[q_0^{jk} \right]_{(n\!-\!1)\times (n\!-\!1)} \end{small} &  \left[ q_0^{jn}\right]_{(n\!-\!1)\times 1}\\
  \left[q_0^{nk}  \right]_{1\times (n\!-\!1)} & q_0^{nn}  \end{BMAT}\! \right]  \end{small}\nonumber \end{eqnarray}
where
 \begin{eqnarray} \!\!\!\!\!&&\!\!\!\!\!\label{19.6.13-6}\end{eqnarray}\begin{eqnarray}
 &&\!\!\!\!  q_0^{jk}  = 2s_1{\tilde{s}}_1 \kappa_j \delta_{jk} + \big(2 s_1{\tilde{s}}_2 + \frac{2i\mu}{\lambda+2\mu} {\tilde{s}}_4 \big)\frac{\xi_j \kappa_k \xi_k}{|\xi'|^2}\nonumber\\
 &&  + (2s_1 {\tilde{s}}_6 +2(s_4-i){\tilde{s}}_9 -2s_4 {\tilde{s}}_8) \frac{\kappa_j \xi_j \xi_k}{|\xi'|^2} +  \bigg\{s_1\big({\tilde{s}}_2 +{\tilde{s}}_6 + {\tilde{s}}_{10}+{\tilde{s}}_{22}\big)\nonumber\\
 &&\, +s_2 \big({\tilde{s}}_1 +{\tilde{s}}_2 + {\tilde{s}}_{6}+{\tilde{s}}_{10}\big)
 +s_3{\tilde{s}}_{22}
 +\big(s_4-\frac{i(\lambda+\mu)}{\mu} \big)  \big({\tilde{s}}_{8} +{\tilde{s}}_{14}\big)\nonumber
\\
 && \!\!\!\! \quad + s_5 \big({\tilde{s}}_4 +{\tilde{s}}_{18}\big)
 \bigg\}\big(\sum_{\alpha} \kappa_\alpha \big)
 \,\frac{\xi_j\xi_k}{|\xi'|^2} + \bigg\{ s_2{\tilde{s}}_2 +(s_1+s_2) {\tilde{s}}_{10}+\frac{2i\mu}{\lambda+2\mu} {\tilde{s}}_{18} \nonumber\\
  &&\,\,\, -2i {\tilde{s}}_{14}    -s_1 \big({\tilde{s}}_2 +{\tilde{s}}_6 + {\tilde{s}}_{22}\big)+ s_2 \big({\tilde{s}}_1  + {\tilde{s}}_{6}\big)-s_3 {\tilde{s}}_{22} \bigg\}\big(\sum_{\alpha} \kappa_\alpha \xi_\alpha^2\big)
\frac{\xi_j\xi_k}{|\xi'|^4},\nonumber\\
&&  \!\!     q_0^{jn}  =\!  \bigg\{s_1\big({\tilde{s}}_4 \!+\!{\tilde{s}}_9 \!+\! {\tilde{s}}_{14}\!-\!{\tilde{s}}_{18}\big)
 + s_2 \big({\tilde{s}}_9 +{\tilde{s}}_{14} \big) + s_3 \big( {\tilde{s}}_{4}-  {\tilde{s}}_{18}\big)\!\nonumber\\
 &&\;\; +\! \big(s_4-\frac{i(\lambda+\mu)}{\mu} \big)  \big({\tilde{s}}_{1} \!+\! {\tilde{s}}_{2}\!+\!{\tilde{s}}_{7} \!-\!{\tilde{s}}_{10}\big)  +s_5{\tilde{s}}_{22}\bigg\} \big(\sum_{\alpha} \kappa_\alpha \big)\nonumber
 \,\frac{\xi_j}{|\xi'|}
 \\
  &&\;\; + \bigg\{ (s_1+s_2){\tilde{s}}_{14} + \frac{2i\mu}{\lambda+2\mu} {\tilde{s}}_{22}  -s_1 ({\tilde{s}}_{4} +{\tilde{s}}_{9}  - {\tilde{s}}_{18})\nonumber\\
\!\!\!\! && \;\;\,   + s_3 \big(-{\tilde{s}}_{4} +{\tilde{s}}_{18}\big)  + s_2 {\tilde{s}}_{9}
   \!-\! 2i s_2 +2(i-2s_4) {\tilde{s}}_{10}  \! \bigg\}\big(\sum_{\alpha} \kappa_\alpha\xi_\alpha^2  \!\big)\frac{\xi_j}{ |\xi'|^3}
  \nonumber\\
  && \;\;\; +      \bigg\{\!-2i({\tilde{s}}_{1}+{\tilde{s}}_{7}) \!-\! 2 s_2 {\tilde{s}}_{9}\!+\! 2(s_1+s_2){\tilde{s}}_{8}\!\bigg\}
    \frac{\kappa_j \xi_j}{|\xi'|},\nonumber\\
&&  \!\!     q_0^{nk}  = \bigg\{ \frac{2i\mu}{\lambda+2\mu} ({\tilde {s}}_1+ {\tilde{s}}_3) +2 s_1 {\tilde{s}}_5\bigg\} \frac{\kappa_k\xi_k}{|\xi'|} +  \bigg\{s_1\big({\tilde{s}}_5 +{\tilde{s}}_8 + {\tilde{s}}_{15}+{\tilde{s}}_{19}\big)
\nonumber\\
&&\quad
+   s_2 \big({\tilde{s}}_5 +{\tilde{s}}_{19} \big)+ s_3  \big({\tilde{s}}_8 +{\tilde{s}}_{15} \big) +\big(s_4-\frac{i(\lambda+\mu)}{\mu} \big)  {\tilde{s}}_{23} \nonumber \\
&&  + s_5 \big( {\tilde{s}}_{1}+  {\tilde{s}}_{3}+
 {\tilde{s}}_{6}+  {\tilde{s}}_{11}\big)
\bigg\} \big(\sum_{\alpha} \kappa_\alpha \big)
 \,\frac{\xi_k}{|\xi'|} + \bigg\{ (s_1+s_2) {\tilde{s}}_{19} +s_2 {\tilde{s}}_{5} -2i {\tilde{s}}_{23}\nonumber\\
 &&
    -s_1  \big({\tilde{s}}_{5} +{\tilde{s}}_{8}+ {\tilde{s}}_{15} \big)- s_3 \big(  {\tilde{s}}_{8} +{\tilde{s}}_{15} \big) +\frac{2i\mu}{\lambda+2\mu}({\tilde{s}}_{6}+{\tilde{s}}_{11}) \bigg\}\big(\sum_\alpha \kappa_\alpha \xi_\alpha^2 \big) \frac{\xi_k}{|\xi'|^3},\nonumber\\
&&  \!\!\! \! q_0^{nn}  =  \bigg\{s_1\big({\tilde{s}}_3 +{\tilde{s}}_7 - {\tilde{s}}_{11}+{\tilde{s}}_{23}\big)+ s_2  {\tilde{s}}_{23}
+   s_3 \big({\tilde{s}}_1 +{\tilde{s}}_{3}+ {\tilde{s}}_{7} -{\tilde{s}}_{11}  \big) \nonumber\\
 && + \big(s_4-\frac{i(\lambda+\mu)}{\mu} \big) \big(  {\tilde{s}}_{5} - {\tilde{s}}_{19} \big)  + s_5 \big( {\tilde{s}}_{9}+  {\tilde{s}}_{15}\big)
\bigg\} \big(\sum_{\alpha} \kappa_\alpha \big)\nonumber\\
&&
 + \bigg\{ -2i ({\tilde{s}}_{5}-{\tilde{s}}_{19})+ \frac{2i\mu}{\lambda
  +2\mu} \big({\tilde{s}}_{8}+{\tilde{s}}_{15}\big) -s_1 \big({\tilde{s}}_{3} +{\tilde{s}}_{7}- {\tilde{s}}_{11}\big)
   \nonumber\\
    &&    +(s_1+s_2){\tilde{s}}_{23}+  s_3 \big(- {\tilde{s}}_{1} -{\tilde{s}}_{3} -{\tilde{s}}_{7} +{\tilde{s}}_{11} \big)
   \bigg\}\big(\sum_\alpha \kappa_\alpha\xi_\alpha^2 \big) \frac{1}{|\xi'|^2}.\nonumber \end{eqnarray}
Thus \begin{eqnarray} \label{19.4.30-3} &&
\quad  \;p_0=\! \left[\begin{BMAT}(@, 1pt, 1pt){c.c}{c.c}
 \left[\mu \,\delta_{jk}\right]_{(n-1)\times (n-1)} &  \left[0\right]_{(n-1)\times 1}\\
 \left[0\right]_{1\times (n-1)} &  \lambda +2\mu\end{BMAT} \right]
 q_0\!-\!\left[\begin{BMAT}(@, 1pt, 1pt){c.c}{c.c} \left[0\right]_{(n-1)\times (n-1)} &  \left[0\right]_{(n-1)\times 1}\\
 \left[0\right]_{1\times (n-1)} &  \lambda \,\sum_{\alpha =1}^{n-1}{\kappa _{\alpha}}\end{BMAT}\!\right]\qquad \qquad \qquad  \\
 && \;\;\quad \quad = \!\begin{small} \left[\!\begin{BMAT}(@, 1pt, 1pt){c.c}{c.c}
 \left[\!\mu q_0^{jk}\right]_{(n-1)\times (n-1)}  & \mu  \left[q_0^{jn}\right]_{(n-1)\times 1} \\ \left[(\lambda+2\mu)q_0^{nk}\right]_{1\times (n-1)} & (\lambda+2\mu)q_0^{nn} -  \lambda \,\sum_{\alpha =1}^{n-1}{\kappa _{\alpha}} \end{BMAT}
\! \right] \end{small},\nonumber\end{eqnarray}
 where $q_0^{jk}, q_0^{jn}, q_0^{nk}, q_0^{nn}$ are given by (\ref{19.6.13-6}).
Since $(p_1(x',\xi')-\tau) \psi_{-1} =I$, we have that for $1\le \alpha \le n-1$,
\begin{eqnarray*} \frac{\partial (p_1 (x', \xi')-\tau I)}{\partial x_\alpha} \, \psi_{-1} +\big(p_1(x', \xi') -\tau I\big) \frac{\partial \psi_{-1}}{\partial x_\alpha} =0.\end{eqnarray*}
Combining this and $\frac{\partial (p_1(x',\xi') -\tau I)}{\partial x_\alpha} \big|_{x=x_0}= 0$ (by (\ref{18/11/24})), we get
\begin{eqnarray*} \frac{\partial \psi_{-1}}{\partial x'_\alpha}(x_0,\xi')=0 \quad \, \mbox{for } 1\le \alpha \le n-1.\end{eqnarray*}
It follows from these and (\ref{18/11/26--2}) that
 \begin{eqnarray} \label{9.15-1} \psi_{-2} (x_0,\xi') = -\psi_{-1} (x_0,\xi') \,  p_0(x'_0, \xi') \psi_{-1} (x_0, \xi').\end{eqnarray}
By direct  calculations, we get
\begin{eqnarray} \label{2019.4.24-3} \end{eqnarray}
\begin{align*}
&    \!\!\!\! \mbox{Tr}\left( \psi _{-1}\,p_0\,\psi _{-1} \right) =\mbox{Tr}\, \left( \psi_{-1}^2\,p_0 \right)
\nonumber \\
&= \!\bigg( \!(n+1)\mu s_1{\tilde{s}}_1 \sum\limits_{\alpha} \kappa_\alpha -(n\!-\!1) \mu s_1{\tilde{s}}_1 \frac{\sum_{\alpha} \kappa_\alpha \xi_\alpha^2 }{|\xi'|^2} \!+\!M_2\!+\!M_4\!\bigg) \frac{1}{(s_1\mu|\xi'|-\tau)^2}\\
&   \;\; \;  \!+\!  \big( s_5(\lambda\!+\!2\mu)\!-\!i\lambda\big) \mu (s_4 \!-\!i)(M_1|\xi'|^2 \!+\!M_3)\frac{(1\!-\!\omega |\xi'|^2)^2}{ \big( (s_1\!+\!s_3)(\lambda\!+\!2\mu)|\xi'|\!-\!\tau\big)^2 (s_1 \mu|\xi'|\!-\!\tau)^2}   \\
& \; \!\;\; +\!   \frac{M_1}{ \big( (s_1\!+\!s_3)(\lambda+2\mu)|\xi'|\!-\!\tau\big)^2 } \\
& \;\;\;\!+\!2M_1\big( s_5(\lambda\!+\!2\mu)\!-\!i\lambda\big) \mu (s_4 \!-\!i)|\xi'|^2 \frac{ 1-\omega |\xi'|^2}{ \big( (s_1\!+\!s_3)(\lambda\!+\!2\mu)|\xi'|\!-\!\tau\big)^3 (s_1 \mu|\xi'|\!-\!\tau)}   \\
 & \;\;\;+  M_1\big( s_5(\lambda\!+\!2\mu)\!-\!i\lambda\big)^2 \mu^2 (s_4 \!-\!i)^2|\xi'|^4 \frac{(1\!-\!\omega |\xi'|^2)^2}{ \big( (s_1\!+\!s_3)(\lambda\!+\!2\mu)|\xi'|\!-\!\tau\big)^4 (s_1 \mu|\xi'|\!-\!\tau)^2}\\
&  \;\;\;     + M_3 |\xi'|^2 \frac{\omega^2}{(s_1 \mu|\xi'|-\tau)^2} - 2M_3 \frac{\omega}{ (s_1\mu |\xi'|-\tau)^2} \\
& \;\;\;+ (M_5+M_6) \bigg\{ \frac{(1-\omega|\xi'|^2)^2}{\big((s_1+s_3)(\lambda+2\mu)|\xi'|-\tau\big)(s_1 \mu|\xi'|-\tau)^2} \\
& \;\;\;     + \frac{1-\omega|\xi'|^2}{\big((s_1+s_3)(\lambda+2\mu)|\xi'|-\tau)^2(s_1 \mu|\xi'|-\tau)}\\
& \;\;\;+ \frac{\big( s_5(\lambda+2\mu)-i\lambda\big)\mu (s_4-i) |\xi'|^2(1-\omega|\xi'|^2)^2}{\big((s_1+s_3)(\lambda+2\mu)|\xi'|-\tau\big)^3(s_1 \mu|\xi'|-\tau)^2} \bigg\}
,\end{align*}
 where $\omega$ is given by (\ref{19.6.13.2}) and
\begin{align*} &\!\!\! \!M_1:=\! \sum\limits_{\alpha}\kappa_\alpha \bigg\{\! -\lambda+(\lambda+2\mu)\bigg( (s_1+s_3)({\tilde{s}}_{1}+ {\tilde{s}}_{3}+ {\tilde{s}}_{7} -{\tilde{s}}_{11})+(s_1+s_2) {\tilde{s}}_{23}  \\
&\;\;\; +  \big(s_4- \frac{i(\lambda+\mu)}{\mu}\big) ({\tilde{s}}_{5}-{\tilde{s}}_{19})
+\! s_5 ({\tilde{s}}_{9}+{\tilde{s}}_{15})\!\bigg)\!\bigg\}  \\
& \;\;\;- \frac{\sum_{\alpha }\kappa_\alpha \xi_\alpha^2}{|\xi'|^2}\bigg\{ (\lambda+2\mu)\bigg((s_1+s_3)( {\tilde{s}}_{1}+{\tilde{s}}_{3}+{\tilde{s}}_{7}-{\tilde{s}}_{11}) \\
&\;\;\; +\!(s_4\!+\!i) ({\tilde{s}}_{5}\!-\!{\tilde{s}}_{19}) \!-\!\frac{i\mu}{\lambda\!+\!2\mu} ({\tilde{s}}_{8}\!+\!{\tilde{s}}_{15})\!\bigg)\!\bigg\}\\
&      =-\,\frac{\lambda ^3+\text{4\,}\lambda ^2\,\mu +\lambda \,\mu ^2-\text{8\,}\mu ^3}{\left( \lambda +\text{3\,}\mu \right) ^2}\sum_{\alpha}{\kappa _{\alpha}}- \frac{\mu \left( \text{3\,}\lambda ^2+\text{12\,}\lambda \,\mu +\text{13\,}\mu ^2\right) }{2\left( \lambda +\text{3\,}\mu \right)^2}\,\frac{\sum_{\alpha}{\kappa _{\alpha}\xi_{\alpha}^2}}{|\xi'|^2},    \\
&  \!\!\!   M_2:= \sum\limits_{\alpha} \kappa_\alpha\, \mu\bigg( s_1({\tilde{s}}_{2}+ {\tilde{s}}_{6}+{\tilde{s}}_{10}+{\tilde{s}}_{22})
+ s_2({\tilde{s}}_{1}+ {\tilde{s}}_{2}+{\tilde{s}}_{6}+{\tilde{s}}_{10})+ s_3 {\tilde{s}}_{22} \\
& \;\;\;\;\;\;+\big(s_4 -\frac{i(\lambda+\mu)}{\mu} \big)({\tilde{s}}_{8}+{\tilde{s}}_{14})  +
s_5( {\tilde{s}}_{4}+{\tilde{s}}_{18})  \bigg)\\
& \;\;\;\;\;\;- \frac{\sum_{\alpha }\kappa_\alpha \xi_\alpha^2}{|\xi'|^2}\mu
\bigg(\! s_1({\tilde{s}}_{2}\!+\! {\tilde{s}}_{6}\!+\!{\tilde{s}}_{10}\!+\!{\tilde{s}}_{22})\! +\!s_2({\tilde{s}}_{1}\!+\! {\tilde{s}}_{2}\!+\!{\tilde{s}}_{6}\!+\!{\tilde{s}}_{10})\! +\!2s_4 ( {\tilde{s}}_{14}\!+\!{\tilde{s}}_{8}) \!+\! s_3  {\tilde{s}}_{22}\!\bigg)\\
&   \;\;  \; \,\,  =-\frac{\left(\!\lambda\! +\!\mu \!\right) ^2}{2\left(\! \lambda \!+\!\text{3\,}\mu \right) ^2} \left(\,\left( {2\,}\lambda \!+\!\text{5\,}\mu \right)\sum_{\alpha}{\kappa _{\alpha}}\! -\!\mu \frac{\sum_{\alpha}{\kappa _{\alpha}\xi _{\alpha}^{2}}}{|\xi'|^2} \right) ,    \\
&\!\!\!   M_3:= \sum\limits_{\alpha} \kappa_\alpha |\xi'|^2  \mu\bigg( (s_1+s_2)({\tilde{s}}_{1}+ {\tilde{s}}_{2}+{\tilde{s}}_{6}+{\tilde{s}}_{10})
+ (s_1+s_3){\tilde{s}}_{22} +s_5({\tilde{s}}_{4}+ {\tilde{s}}_{18})\nonumber\\
& \;\;\;\; +\big(s_4\!-\!\frac{i(\lambda\!+\!\mu)}{\mu} \big) ({\tilde{s}}_{8}\!+\!{\tilde{s}}_{14})
\!+\! \sum\limits_{\alpha}\!\kappa_\alpha \xi_\alpha^2\mu \!\bigg(\!( {\tilde{s}}_{1}\!+\!{\tilde{s}}_{2})({\tilde{s}}_{1}\!+\! {\tilde{s}}_{2}\!+\!{\tilde{s}}_{6}\!+\!{\tilde{s}}_{10})
\!-\!(s_1\!+\!s_3) {\tilde{s}}_{22} \nonumber\\
&\;\;\;\;-
  2s_4({\tilde{s}}_{14}\!+\! {\tilde{s}}_{8}) \!+\!\frac{2i\mu}{\lambda+2\mu} ({\tilde{s}}_{4}\!+\!{\tilde{s}}_{18}) \!+\!2(s_4\!-\!i)({\tilde{s}}_{9}\!+\! {\tilde{s}}_{14})\!\bigg)\\
   &  =-\frac{\left( \lambda +\text{2\,}\mu \right)}{\left( \lambda +\text{3\,}\mu \right) ^2} \left(\left( \lambda ^2+\text{2\,}\lambda \,\mu -\mu ^2 \right)|\xi'|^2 \sum_{\alpha}{\kappa _{\alpha}}\, +\text{2\,}\lambda \,\mu \sum_{\alpha}{\kappa _{\alpha}\xi _{\alpha}^{2}} \right),       \\
   &    \!\!\!   M_4 := \frac{\sum_{\alpha} \kappa_\alpha \xi_\alpha^2}{|\xi'|^2}2\mu\bigg( (s_1+s_2) ({\tilde{s}}_{2}+ {\tilde{s}}_{6}+{\tilde{s}}_{10})+s_2{\tilde{s}}_1\\
    & \;\;\;\;\,+ \frac{i\mu}{\lambda+2\mu} ({\tilde{s}}_{4}+{\tilde{s}}_{18})+(s_4 -i)({\tilde{s}}_9 +{\tilde{s}}_{14})\bigg)\\
    & =-\,\frac{\mu \,\left( \lambda +\mu \right) \,\left( \text{3\,}\lambda +5\mu \right)}{\left( \lambda +\text{3\,}\mu \right) ^2}\,\frac{\sum_{\alpha}{\kappa _{\alpha}\xi _{\alpha}^{2}}}{\left| {\xi '} \right|^2},       \\
   & \!\!\!  M_5:= \mu(s_4-i) (\lambda+2\mu) \bigg\{ - \sum\limits_{\alpha} \kappa_\alpha |\xi'|  \bigg( s_1({\tilde{s}}_{5}+ {\tilde{s}}_{8}+{\tilde{s}}_{15}+{\tilde{s}}_{19})
+ s_2({\tilde{s}}_{5}+ {\tilde{s}}_{19}) \\
&\;\,\;\;\; +s_5({\tilde{s}}_{1}+{\tilde{s}}_{3}+ {\tilde{s}}_{6}+ {\tilde{s}}_{11})+  s_3 ({\tilde{s}}_{8}+ {\tilde{s}}_{15})+ \big(s_4- \frac{i(\lambda+\mu)}{\mu} \big){\tilde{s}}_{23} \bigg)\\
&\;\,\;\;\; + \frac{\sum_{\alpha}\!\kappa_\alpha \xi_\alpha^2}{|\xi'|} \bigg(\!
  s_1(- {\tilde{s}}_{5}\!+\!{\tilde{s}}_{8}+{\tilde{s}}_{15}\!-\! {\tilde{s}}_{19})\!-s_2({\tilde{s}}_{5}\!+\!{\tilde{s}}_{19})
 \qquad \;\,\\
&  \;\;\; \;\, +\!s_3( {\tilde{s}}_{8}\!+\!{\tilde{s}}_{15})- \frac{2i\mu}{\lambda+2\mu}({\tilde{s}}_{1}\!+\! {\tilde{s}}_{3} \!+\!
{\tilde{s}}_{6}+{\tilde{s}}_{11}) +2 i {\tilde{s}}_{23}\bigg)\bigg\}\\
  &  =\!\frac{-\text{2\,}\,\mu ^2}{\! \left( \lambda \!+\!\text{3\,}\mu \right) ^3}\!\left( \!\left( \lambda \!+\!\mu \right) \left( \lambda \!+\!\text{2\,}\mu \right) ^2 \left| {\xi '} \right|\,\sum_{\alpha}{\kappa _{\alpha}}\!+\!\mu \left( \text{3\,}\lambda ^2\!+\!\text{11\,}\lambda \,\mu \!+\!\text{12\,}\mu ^2 \right)\frac{\sum_{\alpha}{\kappa _{\alpha}\xi _{\alpha}^{2}}}{|\xi'|} \right),  \\
     & \!\!\!  M_6:= \mu\big(s_5(\lambda\!+\!2\mu) \!-\!i\lambda\big)  \bigg\{ \!\!-\! \sum\limits_{\alpha} \kappa_\alpha |\xi'|   \bigg( s_1({\tilde{s}}_4\!+\! {\tilde{s}}_{9}\!+\!{\tilde{s}}_{14}\!-\!{\tilde{s}}_{18})
\!+\! s_2({\tilde{s}}_{9}\!+\! {\tilde{s}}_{14}) \\
& \;\;\;\;+s_3({\tilde{s}}_{4}-{\tilde{s}}_{18})  +  s_5{\tilde{s}}_{22}  + \big( s_4-\frac{i(\lambda+\mu)}{\mu}\big) ({\tilde{s}}_{1}+ {\tilde{s}}_{2}+{\tilde{s}}_{7}-{\tilde{s}}_{10} )   \bigg)\\
&\;\;\; \;+ \frac{\sum_{\alpha}\!\kappa_\alpha \xi_\alpha^2}{|\xi'|} \bigg(\!
  s_1( {\tilde{s}}_{4}\!+\!{\tilde{s}}_{9} \!+\!{\tilde{s}}_{14}\!-\!{\tilde{s}}_{18})\!+\!s_2({\tilde{s}}_{9}\!+\!{\tilde{s}}_{14})
\!+\!4s_4 {\tilde{s}}_{10}\\
&\; \;\;\;+\!2i ({\tilde{s}}_{1}\!+\!{\tilde{s}}_{2}\!+\!{\tilde{s}}_{7}\!-\!{\tilde{s}}_{10}) \!+\!s_3 ({\tilde{s}}_{4}\!-\! {\tilde{s}}_{18}) \!-\! 2 (s_1\!+\!s_2) ({\tilde{s}}_{8}\!+\! {\tilde{s}}_{14}) \!-\!\frac{2i\mu}{\lambda\!+\!2\mu}{\tilde{s}}_{22}\bigg)\bigg\}\\
&          =\frac{-\text{2\,}\mu ^2}{\left( \lambda +\text{3\,}\mu \right) ^4}\left( \left( \lambda +\mu \right)\left( \lambda +\text{3\,}\mu \right) \left( \lambda +\text{2\,}\mu \right)^2\left|{\xi '} \right|\,\sum_{\alpha}{\kappa _{\alpha}}\right.\\
& \;\;\;\left.\;\;+\,\mu \left( \text{2\,}\lambda ^2+\text{7\,}\lambda \,\mu +\text{7\,}\mu ^2 \right)\left( \lambda +\text{5\,}\mu \right) \frac{\sum_{\alpha}{\kappa _{\alpha}\xi _{\alpha}^{2}} }{|\xi'|}\right).   \end{align*}
It follows from residue theorem of complex variable function (see \cite{Ahl}) that
\begin{align*}
&    \frac{1}{2\pi i}\int_{\mathcal{C}}{ \frac{1}{\left(s_1 \mu \left| {\xi '} \right|-\tau \right) ^2}\,e^{-t\tau}\,d\tau}=-te^{-t\mu \left| {\xi '} \right|},\\
&  \frac{1}{2\pi i}\int_{\mathcal{C}}{\frac{(1-\omega |\xi'|^2)^2}{ \big( (s_1+s_3)(\lambda+2\mu)|\xi'|-\tau\big)^2 (s_1 \mu|\xi'|-\tau)^2}\,e^{-t\tau}\,d\tau}\\
& \;\;\;\;\,\,=
-\frac{\left( \lambda +{3}\mu \right) ^2\left({2}\left| {\xi '} \right|t\mu^2+{3}\mu +\lambda \right)}{32\mu ^6\left| {\xi '} \right|^3}{e}^{-{2}\mu \left| {\xi '} \right|t}\\
& \;\;\;\, \,\quad \,\,  +\frac{\left( \lambda +{3}\mu \right) ^2\left( -{2}\left| {\xi '} \right|t\mu ^2+{3}\mu +\lambda \right)}{32\mu ^6\left| {\xi '} \right|^3}{e}^{-\frac{{2}\mu \left| {\xi '} \right|t\left( \lambda +\mu \right)}{\lambda +{3}\mu}},\\
&\frac{1}{2\pi i}\int_{\mathcal{C}}{  \frac{1}{ \big( (s_1+s_3)(\lambda+2\mu)|\xi'|-\tau\big)^2 } \,e^{-t\tau}\,d\tau}=-t{e}^{-\frac{2\mu \left| {\xi '} \right|t\left( \lambda +2\mu \right)}{\lambda +3\mu}},\\
&    \frac{1}{2\pi i}\int_{\mathcal{C}}{\frac{ 1-\omega |\xi'|^2}{ \big( (s_1\!+\!s_3)(\lambda\!+\!2\mu)|\xi'|\!-\!\tau\big)^3 (s_1 \mu|\xi'|\!-\!\tau)}\,e^{-t\tau}\,d\tau}=\frac{\left( \lambda +3\mu \right)^3}{16\mu ^6\,\left| \xi' \right|^3 }e^{-2\mu \left| \xi' \right|t}\\
& \;\;\;\;\,\,\quad \quad-\frac{\left( \lambda +3\mu \right)^3}{16 \mu^6\left| \xi' \right|^3}e^{-\frac{2\mu \left| \xi' \right|t\left( \lambda +\mu \right)}{\lambda +3\mu}} +\frac{t\left( \lambda +3\mu \right) ^2}{4\mu ^4\left| \xi' \right|^2 }e^{-\frac{2\mu \left| \xi' \right|t\left( \lambda +2\mu \right)}{\lambda +3\mu}},\\
&      \frac{1}{2\pi i}\int_{\mathcal{C}}{\frac{(1-\omega |\xi'|^2)^2}{ \big( (s_1+s_3)(\lambda+2\mu)|\xi'|-\tau\big)^4 (s_1 \mu|\xi'|-\tau)^2}\,e^{-t\tau}\,d\tau}\\
& \;\;\;\;\,\,=-\frac{\left( \lambda +3\mu \right) ^4\left( 2\left| {\xi '} \right|t\mu ^2+9\mu +3\lambda \right)}{128\mu ^{10}\left|{\xi '} \right|^5}\, e^{-2\mu \left| {\xi '} \right|t}\\
& \;\;\;\;\quad\,\,    +\frac{\left( \lambda \!+\!3\mu \right) ^4\left( -2\left| \xi' \right|t\mu ^2\!+\!9\mu \!+\!3\lambda \right)}{128\mu ^{10}\left| \xi' \right|^5}\,{e}^{-\frac{2\mu \,\left| {\xi '} \right|\,t\,\left( \lambda \!+\!\mu \right)}{\lambda \!+\!\text{3\,}\mu}} \! -\!\frac{t\left( \lambda \!+\!3\mu \right) ^4}{16\mu ^8\left| \xi' \right|^4}e^{-\frac{\text{2\,}\mu \,\left| \xi' \right|t\left( \lambda \!+\!2\mu \right)}{\lambda +3\mu}},\\
 &     \frac{1}{2\pi i}\int_{\mathcal{C}}{ \frac{\omega^2}{(s_1 \mu|\xi'|-\tau)^2}\,e^{-t\tau}\,d\tau}=
 \frac{2\lambda +2\mu +\mu ^2| \xi ' |t-\lambda \mu | \xi '|t}{\mu \left| \xi' \right|^5\left( \lambda -\mu \right)}e^{-\mu \left| \xi' \right|t}\\
 &\quad \,\,\;\;\;\; +\frac{-2\left| \xi' \right|t\mu^2-5\mu +\lambda}{8\mu^2\left| \xi' \right|^5}\text{e}^{-2\mu \left| \xi' \right|t}\qquad \qquad
\\
&   \;\;  \;  \,\,   \quad  \;  -\frac{\lambda ^2+\text{2\,}\left| \xi' \right|t\lambda \mu^2+10\lambda\mu -2\left| \xi' \right|t\mu^3+21\mu^2}{8\mu ^2\left| \xi' \right|^5\left( \lambda -\mu \right) }\,e^{-\frac{2\mu \left|\xi' \right|t\left( \lambda +\mu \right)}{\lambda +3\mu}},\\
&        \frac{1}{2\pi i}\int_{\mathcal{C}}{  \frac{\omega}{ (s_1\mu |\xi'|\!-\!\tau)^2} e^{-t\tau}\,d\tau}=\frac{\lambda\! +\!\mu \!+\!\mu ^2| \xi' |t\!-\!\lambda\mu | \xi' |\,t}{\mu|\xi'|^3\left( \lambda \!-\!\mu \right)}e^{-\mu | \xi' |\,t}\!  \\
 & \;\;\;\;\,\,\quad-\!\frac{1}{2\mu | \xi |^3}\,e^{-2\mu \left| \xi' \right|t}-\frac{\lambda \!+\!3\mu}{2| \xi' |^3\left( \lambda\! -\!\mu \right) \mu }e^{-\frac{2\mu | \xi'|t\left( \lambda \!+\!\mu \right)}{\lambda +3\mu}},\\
 &  \frac{1}{2\pi i}\int_{\mathcal{C}}{  \frac{(1-\omega|\xi'|^2)^2}{\big((s_1+s_3)(\lambda+2\mu)|\xi'|-\tau\big)(s_1 \mu|\xi'|-\tau)^2} e^{-t\tau}\,d\tau}\\
 &\;\;\; \;\,\,=
 \frac{t\left( \lambda +3\mu \right)}{8\mu^2| \xi' |}e^{-2\mu | \xi'|t}-\frac{t\left( \lambda +3\mu \right)}{8\mu^2| \xi' |}e^{-\frac{2\mu | \xi'|t\left( \lambda +\mu \right)}{\lambda +3\mu}},\\
 &     \frac{1}{2\pi i}\int_{\mathcal{C}}{\frac{1-\omega|\xi'|^2}{\big((s_1+s_3)(\lambda+2\mu)|\xi'|-\tau)^2(s_1 \mu|\xi'|-\tau)} e^{-t\tau}\,d\tau}\\
 &\; \;\;\;\,\,=-\frac{\left( \lambda +3\mu \right) ^2}{8\mu ^4\left| \xi' \right|^2}\,e^{-2\mu \left| \xi' \right|t}-\frac{\left( \lambda +3\mu \right)^2}{8\mu ^4\left|\xi' \right|^2} e^{-\frac{2\mu \left| \xi' \right|t\left( \lambda +\mu \right)}{\lambda +3\mu}} +\frac{\left( \lambda +3\mu \right) ^2}{4\mu ^4\left| \xi' \right|^2}e^{-\frac{2\mu \left| \xi' \right|t\left( \lambda +2\mu \right)}{\lambda +3\mu}},\\
 &     \frac{1}{2\pi i}\int_{\mathcal{C}}{\frac{(1-\omega|\xi'|^2)^2}{\big((s_1+s_3)(\lambda+2\mu)|\xi'|-\tau\big)^3(s_1 \mu|\xi'|-\tau)^2} e^{-t\tau}\,d\tau}\\
 &\;\;\;\;\, \,= \frac{\left( \lambda +3\mu \right) ^3\left( \left| \xi' \right|t\mu^2+3\mu +\lambda \right)}{32\mu^8\left| \xi' \right|^4 }e^{-2\mu \left| \xi' \right|t}\\
&\;\;\;\;\quad \,\, +\frac{\left( \lambda \!+\!3\mu \right) ^3\left( -\left| \xi' \right|t\mu^2\!+\!3\mu \!+\!\lambda \right)}{32\mu ^8\left| \xi' \right|^4 }e^{-\frac{2\mu \left| \xi ' \right|t\left( \lambda +\mu \right)}{\lambda \!+\!3\mu}} \! -\!\frac{\left( \lambda \!+\!3\mu \right)^4}{16\,\mu ^8\left| \xi' \right|^4}\,e^{-\frac{2\mu \left| \xi' \right|t\left( \lambda +2\mu \right)}{\lambda +3\mu}}.
 \end{align*}
Inserting these, (\ref{19.4.6-8}) and (\ref{2019.10.24-2}) into (\ref{2019.4.24-3})), and then calculating the complex integral $\int_{\mathcal{C}} \mbox{Tr}\left( \psi _{-1}\,p_0\,\psi _{-1} \right)e^{-t\tau} \, d\tau,$ we obtain
\begin{eqnarray} \label{2019.10.26-1} &&\frac{{i}}{2\pi}\int_{\mathcal{C}}{\big(\mbox{Tr}\left( -\psi _{-1}p_0\psi _{-1} \right)\big)e^{-t\tau}}d\tau= \frac{1}{2\pi i}\int_{\mathcal{C}} \left( \mbox{Tr}\left( \psi _{-1}\,p_0\,\psi _{-1} \right)\right) e^{-t\tau} \, d\tau\qquad \qquad \qquad \quad \qquad\;\\
&&\quad \quad \quad  =\!  -\frac{1}{2} n \mu\,  t\,e^{-\mu | \xi' |t}\left(\! \sum_{\alpha}{\kappa_{\alpha}}\!-\!\frac{\sum_{\alpha}{\kappa_{\alpha}\xi_{\alpha}^{2}}}{|\xi'|^2}\!\right)\!\nonumber\\
 &&\qquad \quad \quad  -\mu \,t\, e^{-2\mu \left| \xi' \right|t}\!\bigg\{\! \!\sum_{\alpha}{\kappa_{\alpha}} \!+ \!\frac{ \sum_{\alpha} \kappa_\alpha \xi_\alpha^2} {|\xi'|^2}\bigg(\! \frac  {\lambda ^3\!+\!23\lambda ^2\mu \!+\!87\lambda \mu ^2\!+\!97\mu^3}{4\left( \lambda \!+\!3\mu \right) ^3}\!\bigg) \!\bigg\} \nonumber\\
&&\quad \quad \quad \quad    +t e^{-\frac{2\mu \left( \lambda \!+\!\mu \right)| \xi' |t}{\lambda\! +\!3\mu}}\bigg\{\!\sum_{\alpha}\!{\kappa_{\alpha}}\left(\frac{2\lambda^3\!+\!9\lambda^2\mu \!+\!10\lambda \mu^2\!-\!\mu^3}{(\lambda+3\mu)^2}\right) \nonumber \\
&& \qquad \quad\quad  + \frac{\sum_{\alpha}\!{\kappa_{\alpha}\xi_{\alpha}^{2}}}{|\xi'|^2} \left( \frac{\mu( 15\lambda ^3\!+\!105\lambda^2\mu \!+\!233\lambda \mu^2\!+\!175\mu^3) }{4(\lambda+3\mu)^3} \right)\! \bigg\}.\nonumber
\end{eqnarray}
It follows from \cite{Liu2} that
\begin{eqnarray*}
&\int_{\mathbb{R}^{n-1}}{\begin{array}{c}\!\!
 \left( \sum_{j=1}^{n-1}{\xi _{j}^{2}} \right) ^{\frac{m}{2}}\xi _{k}^{\alpha}e^{-\sqrt{\sum_{j=1}^{n-1}{\xi_{j}^{2}}}}d\xi _1\,d\xi _2\cdots \,d\xi _{n-1}\\
\end{array}}\\
& \;\;\;\; \qquad\quad \; \;\; =\!\begin{cases}
 \varGamma \left( n\!-\!1\!+\!m \right) \,\text{vol}( \mathbb{S}^{n\!-\!2} ) \; &\mbox{for}\;\,\alpha =0\\
 \text{0                              }\; &\mbox{for}\,\;\alpha =\text{1,}\\
\end{cases}\;\, n\geqslant \text{2,}\nonumber
\\
&\int_{\mathbb{R}^{n-1}}{\begin{array}{c}\!\!
 \left( \sum_{j=1}^{n-1}{\xi _{j}^{2}} \right) ^{\frac{m-2}{2}}\xi _k\xi _le^{-\sqrt{\sum_{j=1}^{n-1}{\xi_{j}^{2}}}}d\xi _1\,d\xi _2\cdots d\xi _{n-1}\\
\end{array}}\\
& \;\;\;\;\qquad\quad \;  \;\; =\!\begin{cases}
 \frac{\varGamma \left( n-1+m \right)\, \text{vol}( \mathbb{S}^{n\!-\!2} )}{n-1}\;\,&\mbox{for}\; k=l\\
 \text{0                              }\; &\mbox{for}\;\, k\ne l,\\
\end{cases}\;\;n\geqslant \text{3},\qquad\nonumber
\end{eqnarray*}
so that, for any positive real number $b$, \begin{eqnarray} \label{19.6.15-9} \left. \begin{array}{ll} & \frac{1}{(2\pi)^{n-1}}\int_{{\Bbb R}^{n-1}} e^{-b|\xi|t} d\xi_1\cdots d\xi_{n-1} = \frac{\Gamma(n-1) \,\text{vol}( \mathbb{S}^{n\!-\!2} ) }{(2\pi)^{n-1} (bt)^{n-1}},\\
&\frac{1}{(2\pi)^{n-1}}\int_{{\Bbb R}^{n-1}} \frac{ \xi_\alpha^2}{|\xi'|^2} e^{-b|\xi|t} d\xi_1\cdots d\xi_{n-1} = \frac{\Gamma(n-1) \, \text{vol} (\mathbb{S}^{n\!-\!2} ) }{(2\pi)^{n-1} (n-1) (bt)^{n-1}}.\end{array}\right. \end{eqnarray}
 Combining  (\ref{18/12/2-5}), (\ref{9.15-1}), (\ref{2019.10.26-1}) and (\ref{19.6.15-9}), we find that  \begin{align*} &
\widetilde{a_1}(x')= \frac{1}{(2\pi)^{n-1}}\int_{{\Bbb R}^{n-1}}\bigg( \frac{{i}}{2\pi}\int_C e^{-t\tau}\,\mbox{Tr}( -\psi _{-1}\,p_0\,\psi _{-1} )\,d\tau \bigg) d\xi _1\cdots \, d\xi _{n-1}\quad \\
& =\frac{(n-2)! \big(\text{vol}\,( \mathbb{S}^{n-2} )\big)\big(\sum_\alpha{\kappa _{\alpha}}\big) }{(2\pi)^{n-1}t^{n-2}}\bigg\{
\frac{n(2-n)}{2(n-1)\mu^{n-2}} \\
 & \;\;\;\;- \bigg( 1 + \frac{(\lambda^3\! +\!23 \lambda^2 \mu \!+\!87 \lambda \mu^2 \!+\!97 \mu^3)}{ 4(n\!-\!1)(\lambda\!+\!3\mu)^3} \bigg) \frac{\mu}{(2\mu)^{n-1}} \! +\! \bigg( \frac{2\lambda^3 \!+\!9 \lambda^2\mu\!+\! 10\lambda\mu^2 \!-\!\mu^3}{(\lambda+3\mu)^2}\\
 &\;\;\;\; + \frac{\mu (15\lambda^3 \!+\!105\lambda^2 \mu \!+\!233\lambda\mu^2 \!+\!175 \mu^3)}{4(n-1)(\lambda\!+\!3\mu)^3} \bigg) \bigg(\frac{\lambda+3\mu}{2\mu(\lambda\!+\!\mu)}\bigg)^{n-1}\bigg\}\!:=\!t^{2-n} a_1(x').\end{align*}

 \vskip 0.26 true cm

Step 4. \   For general $m\ge 1$, replacing the matrix $E$ in (\ref{19.3.28-1}) by $E_{-m}$ of (\ref{19.6.1-2}) at the original $x_0$ in local boundary normal coordinates, we can figure out $q_{-m}(x_0,\xi')$, so that $p_{-m}(x_0,\xi')$
is gotten by (\ref{19.3.28-15}). Next, applying (\ref{19.6.16-1})--(\ref{19.6.16-2}) we get $\psi_{-1-m}(x_0, \xi',\tau)$, and hence (\ref{18/12/2-5}) gives \begin{eqnarray*} && \label{018/12/2-5} \widetilde{a_m}(x_0)= \frac{1}{(2\pi)^{n-1}} \int_{{\Bbb R}^{n-1}} \frac{i}{2\pi} \int_{\mathcal{C}} e^{-t\tau}\,
 \mbox{Tr}\,\big(\psi_{-1-m} (x_0,\xi',\tau)\big) d\tau\bigg\} d\xi'\\
 && \qquad \qquad \qquad \; \quad \;\mbox{for}\;\;  0\le m\le n-1.\nonumber\end{eqnarray*}
 Because $x_0$ is any point on $\partial \Omega$, we can obtain all $\widetilde{a_m}(x')$ for $0\le m\le n-1$.  \qed

\vskip 0.38 true cm

 \noindent{\bf Remark 5.2.} \  {\it  \
By applying the Tauberian theorem (see, for example, Theorem 15.3 of p.$\,$30 of \cite{Kor} or p.$\,$107 of \cite{Ta2}) for the first term on the right side of (\ref{6.0.1}), we immediately get the Weyl-type law for the counting function $N(\tau):=
 \#\{k\big|\tau_k \le \tau\} $ of elastic Steklov eigenvalues:
 \begin{eqnarray}  \label{019.6.29-1,,1} && N(\tau)\sim \frac{\big(\text{vol}\,( \mathbb{S}^{n-2} )\big)}{(n-1)\,\mu ^{n-1}}\begin{small}\bigg\{ \frac{1\!+\!\big( \frac{\lambda \!+\!{3}\mu}{\lambda +\mu} \big) ^{n-1}}{2^{n-1}}
 \!+\!\left( n-2 \right)\! \bigg\}\end{small}\big(\text{vol}\,(\partial \Omega) \big)\,\tau^{n-1}\\
 &&\qquad \quad \; \;\; + o(\tau^{n-1}) \quad \mbox{as}\;\;
 \tau\to +\infty.\nonumber\end{eqnarray}}

\vskip 0.38 true cm

 \noindent{\bf Remark 5.3.} \  {\it  \  i) \ \  It can be verified that, when $n\ge 3$, the third coefficient $a_2(x',n)$ in (\ref{19.6.15-1}) has the form
 \begin{eqnarray} \label{19.6.16-8} &&\!\!\!\!\!\!\!\!\!
{a_2}(x')= ( n\!-\!2 )!\,\big(\text{vol}\,( \mathbb{S}^{n\!-\!2} )\big) \!\bigg\{\! h_1 \big(\!\sum_{\alpha=1}^{n\!-\!1} \kappa_\alpha\!\big) \!+\! h_2 \bigg(\!\frac{1}{(n\!-\!1)(n\!-\!2)} \!\sum_{\underset {\alpha\ne \beta}{1\le \alpha,\beta\le n\!-\!1}}\! \kappa_\alpha \kappa_\beta\!\bigg)\\
&&\qquad\;\,\; \!+h_3 R_\Omega + h_4 R_{\partial \Omega}+ h_5 {\mbox{Ric}}_{{}_{\Omega}} +h_6 {\mbox{Ric}}_{{}_{\partial \Omega}} \bigg\},\nonumber\end{eqnarray}
  where $h_l, h_2, h_3, h_4, h_5, h_6$  are constants depending only on $n$, $\mu$ and $\lambda$, and $R_\Omega$ (respectively, ${\mbox{Ric}}_{{}_\Omega}$) is the scalar (respectively, Ricci) curvature of the domain $\Omega$, and $R_{\partial \Omega}$  (respectively, ${\mbox{Ric}}_{{}_{\partial \Omega}}$) is the scalar (respectively, Ricci) curvature of the boundary $\partial \Omega$.}

{\it  ii) \ \  Besides the above obtained $a_0$ and $a_1$,   by our new method we can also get all coefficients $a_l$, $2\le l\le m$, for the asymptotic expansion
   $\sum_{k=1}^\infty e^{-t \tau_k}\sim \sum\limits_{m=0}^{n-1} a_m t^{m+1-n} +o(1)$ as $t\to 0^+$}.

\vskip 1.08 true cm

\centerline {\bf  Acknowledgments}

\vskip 0.39 true cm

The author would like to thank Professors Gunther Uhlmann, Ari Laptev, Robert V. Kohn, Mark S. Ashbaugh, Fang-Hua Lin and John M. Lee for their great support and many useful comments and discussions. This research was supported by NNSF of China (11671033/A010802) \ and  NNSF
of China \ (11171023/A010801).

  \vskip 1.68 true cm

\vskip 0.32 true cm

\end{document}